\documentclass[12pt,a4]{article}
\usepackage{amsmath}
\usepackage{amsfonts}
\usepackage{bm}
\usepackage{cases}
\usepackage{graphicx}
\usepackage{stmaryrd}
\usepackage{amssymb}
\usepackage{ascmac}
\usepackage{comment}
\usepackage{latexsym}
\usepackage{ascmac}
\usepackage{mathrsfs}
\usepackage{bm}			
\usepackage{comment}
\usepackage{latexsym}

\allowdisplaybreaks[4]

\newtheorem{thm}{Theorem}[section]
\newtheorem{lem}[thm]{Lemma}
\newtheorem{prop}[thm]{Proposition}

\numberwithin{equation}{section}

\makeatletter
\newcommand{\opnorm}{\@ifstar\@opnorms\@opnorm}
\newcommand{\@opnorms}[1]{%
  \left|\mkern-1.5mu\left|\mkern-1.5mu\left|
   #1
  \right|\mkern-1.5mu\right|\mkern-1.5mu\right|
}
\newcommand{\@opnorm}[2][]{%
  \mathopen{#1|\mkern-1.5mu#1|\mkern-1.5mu#1|}
  #2
  \mathclose{#1|\mkern-1.5mu#1|\mkern-1.5mu#1|}
}
\makeatother

\newcommand{\wt}{\widetilde}

\newcommand{\del}{\partial}
\newcommand{\ov}{\overline}
\newcommand{\delx}{\partial_{x}}
\newcommand{\delt}{\partial_{t}}
\newcommand{\ts}{\textstyle}
\newcommand{\real}{{\rm Re}\,}

\renewcommand{\div}{\mbox{\rm div}\,}
\newcommand{\lang}{{\langle}}
\newcommand{\rang}{{\rangle}}

\newcommand{\re}{{\rm Re}} 
\newcommand{\im}{{\rm Im}}
\newcommand{\trans}{{}^\top}
\renewcommand{\Pr}{{\rm Pr}}

\newcommand{\Ker}{{\rm Ker}\,}

\newcommand{\bv}{\boldsymbol{v}}
\newcommand{\bw}{\boldsymbol{w}}
\newcommand{\bu}{\boldsymbol{u}}
\newcommand{\bW}{\boldsymbol{W}}
\newcommand{\bU}{\boldsymbol{U}}
\renewcommand{\b}{\boldsymbol}

\newcommand{\ep}{\varepsilon}

\renewcommand{\trans}{{}^{\top}}
\renewcommand{\lang}{{\langle}}
\renewcommand{\rang}{{\rangle}}

\newcommand{\ds}{\displaystyle}

\begin{document}

\title{\large \bf Singular limit in Hopf bifurcation for doubly diffusive convection equations II: bifurcation and stability  
%Hopf bifurcation in the artificial compressible system for doubly diffusive convection
}
\author{\normalsize Chun-Hsiung Hsia${}^1$, Yoshiyuki Kagei${}^2$, Takaaki Nishida${}^3$ \\[1ex] 
\normalsize and Yuka Teramoto${}^4$}
\date{}

\footnotetext[1]{Department of Mathematics, National Taiwan University, Taipei  10617, Taiwan}
\footnotetext[2]{Department of Mathematics, Tokyo Institute of Technology, Tokyo   152-8551, Japan}
\footnotetext[3]{Department of Advanced Mathematical Sciences, Kyoto University, Kyoto  606-8317, Japan}
\footnotetext[4]{Research Institute for Interdisciplinary Science, Okayama University, Okayama 700-8530, Japan}

\maketitle

\begin{abstract}
\noindent
A singular perturbation problem from the artificial compressible system to the incompressible system is considered for a doubly diffusive convection when a Hopf bifurcation from the motionless state occurs in the incompressible system. It is proved that the Hopf bifurcation also occurs in the artificial compressible system for small singular perturbation parameter, called the artificial Mach number. The time periodic solution branch of the artificial compressible system is shown to converge to the corresponding bifurcating branch of the incompressible system in the singular limit of vanishing artificial Mach number. 
\end{abstract}

\section{Introduction}\label{Introduction}
This paper studies a singular limit problem for Hopf bifurcation in the artificial compressible system for thermal convection equations in the presence of the diffusion of solute concentration. The system of equations under consideration is written as   
\begin{align}
\ep^2\delt \phi + \div\bw &= 0, \label{1.1} \\ 
\delt\bw-\Pr\Delta\bw + \Pr\nabla \phi -\Pr\mathcal{R}_1\theta\boldsymbol{e}_2 + \Pr\mathcal{R}_2\psi\boldsymbol{e}_2 + \bw\cdot\nabla\bw &= \boldsymbol{0}, \label{1.2} \\
\delt\theta - \Delta\theta - \mathcal{R}_1\bw\cdot\boldsymbol{e}_2 + \bw\cdot\nabla\theta &= 0, \label{1.3}\\
\delt \psi - d \Delta \psi - \mathcal{R}_2\bw\cdot\boldsymbol{e}_2 + \bw\cdot\nabla \psi &= 0. \label{1.4}
\end{align}
Here $\phi = \phi(x, t)$, $\bw=\trans(w^1(x,t), w^2(x,t))$, $\theta = \theta(x, t)$ and $\psi = \psi(x, t)$ denote the perturbation of the pressure, velocity field, temperature and solute concentration, respectively, at position $x \in \mathbb{R}^{2}$ and time $ t \in \mathbb{R}$, from their values of the motionless state in a  thermal convection; $\mathcal{R}_{1}$, $\mathcal{R}_{2}$, $\Pr$ and $d$ are non-dimensional positive parameters; $\mathcal{R}_{1}^{2}$ and $\mathcal{R}_{2}^{2}$ are the thermal and salinity Rayleigh numbers, respectively; $\Pr$ and $d$ are called the Prandtl and Lewis numbers, respectively. The parameter $\ep > 0$ is called the artificial Mach number.  
  
The system \eqref{1.1}--\eqref{1.4} is considered on the two dimensional infinite layer $\{x =(x_{1}, x_{2}) \in \mathbb{R}^{2}; x_{1} \in \mathbb{R}, \, 0 < x_{2} < 1\}$ under the following boundary condition on the boundary $\{x_2=0,1\}$:  
\begin{equation}\label{1.5}
\mbox{$\ds \frac{\del w^1}{\del x_2} = w^2 = \theta = \psi = 0$ on $\{x_2 = 0,1\}$}
\end{equation}
and the periodic boundary condition in $x_{1}$ for $\phi$, $\bw$, $\theta$ and $\psi$ with period $\frac{2\pi}{\alpha}$, where $\alpha > 0$ . The system \eqref{1.1}--\eqref{1.4} is thus considered in the domain 
\[
\Omega = \mathbb{T}_{\frac{2\pi}{\alpha} } \times (0, 1)
\]
under the boundary condition \eqref{1.5}. Here $\mathbb{T}_{a}$ denotes $\mathbb{T}_a=\mathbb{R}/a\mathbb{Z}$. We write $u=\trans(\phi, \bw, \theta, \psi)$ for the unknown functions of \eqref{1.1}--\eqref{1.4}. Here and in what follows, the superscript $\trans \, \cdot \,$ denotes the transposition. 

When $\ep = 0$ one obtains the thermal convection equations for a viscous incompressible fluid under the Oberbeck-Boussinesq approximation. A rough explanation of the structure of $\eqref{1.1}$--\eqref{1.4} with $\ep = 0$ is as follows. The terms including $\mathcal{R}_{1}$ have a symmetric structure which may cause instabilities against the dissipativity by the Laplacians when $\mathcal{R}_{1}$ increases; the terms including $\mathcal{R}_{2}$ have a skew-symmetric structure which may cause oscillatory behavior. Indeed, when $\ep = 0$, it was proved in \cite{Bona-Hsia-Ma-Wang} that, for some range of $\mathcal{R}_{2}$, there exists a critical number $\mathcal{R}_{1,c}$ such that if $\mathcal{R}_{1} < \mathcal{R}_{1,c}$ then the motionless state $u = 0$ is asymptotically stable, while if $\mathcal{R}_{1} > \mathcal{R}_{1,c}$ then $u = 0$ is unstable and nontrivial time periodic solutions bifurcate from $u = 0$. 

On the other hand, when $\ep > 0$, the system \eqref{1.1}--\eqref{1.4} is a hyperbolic-parabolic system and the limit $\ep \to 0$ is a singular limit from a hyperbolic-parabolic system to a parabolic system. One of mathematical questions is thus whether the artificial compressible system \eqref{1.1}--\eqref{1.4} gives a good approximation of the incompressible system in the singular limit $\ep \to 0$. 

The aim of this paper is to consider whether the artificial compressible system approximates well the incompressible system when a Hopf bifurcation, i.e., a time periodic bifurcation, occurs in the incompressible system. In this paper we shall investigate the Hopf bifurcation problem of the singularly perturbed system \eqref{1.1}--\eqref{1.4} with $0 < \ep \ll 1$ for $\mathcal{R}_{1}$ near a critical value $\mathcal{R}_{1,c}$, where a Hopf bifurcation occurs from the motionless state $u = 0$ in the incompressible system $\eqref{1.1}$--\eqref{1.4} with $\ep = 0$ for $\mathcal{R}_{1} \sim \mathcal{R}_{1,c}$ as was shown in \cite{Bona-Hsia-Ma-Wang}. 

The artificial compressible system \eqref{1.1}--\eqref{1.4} with $\ep > 0$ was introduced by Chorin (\cite{Chorin1, Chorin2, Chorin3}) and Temam (\cite{Temam1, Temam2}) to avoid the computational difficulties in the incompressible system due to the constraint $\div\bw =0$. By using the artificial compressible system \eqref{1.1}--\eqref{1.4} with small $\ep>0$, Chorin \cite{Chorin1, Chorin2, Chorin3} computed stationary convective bifurcating patterns near the onset of convection when $\mathcal{R}_{2} = 0$, which suggests the artificial compressible system \eqref{1.1}--\eqref{1.4} would give a good approximation of the incompressible system. Temam (\cite{Temam1, Temam2, Temam3}) introduced the artificial compressible system with an additional stabilizing nonlinear term (in the context of \eqref{1.1}--\eqref{1.4}, the corresponding stabilizing terms are to be $+\frac{1}{2}(\div \bw)\bw$, $+\frac{1}{2}(\div \bw)\theta$ and $+\frac{1}{2}(\div \bw)\psi$ which are added on the left-hand side of \eqref{1.2}, \eqref{1.3} and \eqref{1.4}, respectively); and for such a system with stabilizing nonlinear term, the global in time solutions of the initial boundary value problem for the artificial compressible system on a two-dimensional bounded domain $D$ converge in the limit $\ep \to 0$ to the one for the incompressible system on time interval $(0, T)$ for all $T>0$. Convergence result was also established in the framework of weak solutions on three-dimensional bounded domains. Donatelli \cite{Donatelli1, Donatelli2} and Donatelli and Marcati \cite{Donatelli-Marcati1, Donatelli-Marcati2} proved similar convergence results in the case of unbounded domains on any finite time interval by using the wave equation structure of the pressure and the dispersive estimates. 

In this paper we shall show that a time periodic bifurcation occurs in \eqref{1.1}--\eqref{1.4} with $0 < \ep \ll 1$ for $\mathcal{R}_{1}$ near the criticality $\mathcal{R}_{1,c}$ and the bifurcating solution branch converges as $\ep \to 0$ to the time priodic bifurcating solution branch of the incompressible system obtained in \cite{Bona-Hsia-Ma-Wang}.  

To prove these results, we shall employ the Lyapunov-Schmidt method in a time periodic function space $\mathcal{X}_{a}$. For this purpose, we need to investigate the spectral properties of the operator $\delt + L^{\ep}_{\mathcal{R}_{1}}$ in $\mathcal{X}_{a}$ for $\mathcal{R}_{1}$ near the criticality $\mathcal{R}_{1,c}$, where $L^{\ep}_{\mathcal{R}_{1}}$ denotes the linearized operator around $u = 0$ for \eqref{1.1}--\eqref{1.4}, which was studied in \cite{Hsia-Kagei-Nishida-Teramoto1} in detail. Based on the spectral properties of the operator $\delt + L^{\ep}_{\mathcal{R}_{1}}$, we shall show that \eqref{1.1}--\eqref{1.4} has a nontrivial bifurcating time periodic solution $u^{\ep}$ for sufficiently small $\ep > 0$ and that $u^{\ep}$ converges to the bifurcating time periodic solution of the incompressible system $\eqref{1.1}$--\eqref{1.4} with $\ep = 0$ as $\ep \to 0$. 

One of the key points in the proof is to use the norms with $\ep$-weights, which enables us to establish the uniform estimates in the nonlinear problem as well as in the linearized problem.  Another key is the fact that the limiting solution of the incompressible system and the eigenfunctions for the critical eigenvalues of $L^{\ep}_{\mathcal{R}_{1} }$ for $\mathcal{R}_{1} \sim \mathcal{R}_{1,c}$ are smooth, which compensate the loss of uniform estimates in $\ep$. There is one more thing to be mentioned. At criticality $\mathcal{R}_{1} = \mathcal{R}^{\ep}_{1,c}$, the linearized problem has nontrivial time periodic solutions and the bifurcating solutions of the nonlinear problem are obtained as a perturbation of a time periodic solution of the linearized problem. The period of the bifurcating solution of the nonlinear problem is also a perturbation of the period of the time periodic solution of the linearized problem. This perturbation of the period is transformed into the equation as a perturbation of the time derivative of the unknown. In the case of the incompressible problem (i.e., in a classical setting), this perturbation term can be regraded as a regular perturbation. In contrast to the classical case, in the case of the artificial compressible problem, this perturbation term is not a regular perturbation. To overcome this difficulty, we put it into the principal part of the linearized operator and establish uniform estimates for small $\ep$. 

In the proof, we make use of the two-dimensional aspect of the problem to estimate the nonlinearity uniformly in small $\ep$. It is not straightforward to extend the argument to the three-dimensional problem.  

We shall also show that the bifurcating time periodic solutions of \eqref{1.1}--\eqref{1.4} for small $\ep > 0$ is stable under perturbations with the same symmetries as those of the time periodic solutions, if the corresponding bifurcating time periodic solutions of the incompressible system is stable. This, in particular, implies that, by using the artificial compressible system with small $\ep >0 $, one can numerically compute bifurcating time periodic solutions which are close to the bifurcating ones of the incompressible system near the onset of convection. See \cite{Kagei-Nishida, Kagei-Nishida-Teramoto, Teramoto} for the stability of stationary bifurcating solutions.    

This paper is organized as follows. In section 2 we introduce notation used in this paper. In section 3 we state the results on the Hopf bifurcation for the incompressible problem obtained in \cite{Bona-Hsia-Ma-Wang}. 
In section 4 we state the result on the occurrence of the Hopf bifurcation in the artificial compressible system with small $\ep$ and also the result on the singular limit $\ep \to 0$ of the time periodic bifurcating solutions. 
In section \ref{Spectrum of B} we summarize the results on the spectral properties of the linearized problem obtained in \cite{Hsia-Kagei-Nishida-Teramoto1}. 
Section \ref{Proof of thm4.2} is devoted to a proof of the existence and singular limit.  
In section \ref{Stability} we state the result on the stability of the bifurcating solutions and give a proof of the stability.

\section{Preliminaries}

In this section we introduce notation used in this paper. Let $\Omega=\mathbb{T}_{\frac{2\pi}{\alpha}}\times(0,1)$. 
We denote by $L^p(\Omega)$ the usual $L^p$ space on $\Omega$ with norm $\|\cdot\|_p$. 
We also denote by $\b L^p(\Omega)$ the space of all $L^p$ vector fields on $\Omega$ with norm $\|\cdot\|_p$. 
The inner product of $L^2(\Omega)$ is denoted by $(\cdot, \cdot)_{L^2}$. We also denote  by $(\cdot,\cdot)_{L^2}$ the inner product of $\b L^2(\Omega)$. 

The inner product and the norm of $\b L^2(\Omega)\times L^2(\Omega)\times L^2(\Omega)$ are defined by 
\[
(\bu_1,\bu_2)=\Pr^{-1}(\bw_1,\bw_2)_{L^2}+(\theta_1,\theta_2)_{L^2}+(\psi_1,\psi_2)_{L^2} 
\]
for $\bu_j=\trans(\bw_j,\theta_j,\psi_j)\in \b L^2(\Omega)\times L^2(\Omega)\times L^2(\Omega)$
\[
\|\bu\|_2=\sqrt{(\bu,\bu)}=\sqrt{\Pr^{-1}\|\bw\|_2^2+\|\theta\|_2^2+\|\psi\|_2^2} 
\]
for  $\bu=\trans(\bw,\theta,\psi)\in \b L^2(\Omega)\times L^2(\Omega)\times L^2(\Omega)$, respectively. 

The symbols $H^k(\Omega)$ and $\b H^k(\Omega)$ stand for the $k$ th order $L^2$ Sobolev spaces for scalar functions and for vector fields, respectively. 

We define the spaces of $L^{2}$ and $H^{k}$ functions with vanishing mean value on $\Omega$ by     
\[
L_*^2(\Omega)=\{\phi \in L^2(\Omega); \int_{\Omega}\phi(x) \,dx=0\}
\]
and 
\[ 
H_*^k(\Omega)=H^k(\Omega)\cap L_*^2(\Omega), 
\]
respectively. 

We shall employ the following function spaces with symmetries: 
\begin{align*}
& L_{sym}^2(\Omega) = \{\theta\in L^2(\Omega); \mbox{ \rm $\theta$ is even in $x_1$} \},
\\
& L_{*,sym}^2(\Omega) =L_{sym}^2(\Omega) \cap L_*^2(\Omega), 
\\
& H^k_{*,sym}(\Omega) =H^k(\Omega)\cap L_{*,sym}^2(\Omega), 
\\
& \b L^2_{sym}(\Omega) =\{\bv = \trans(v^1,v^2)\in \b L^2(\Omega); \mbox{ \rm $v^1$ is odd in $x_1$, $v^2$ is even in $x_1$} \}. 
\end{align*}

It is known (\cite{Galdi, Sohr, Temam3}) that $\b L^2_{sym}(\Omega)$ admits the Helmholtz decomposition: 
\[
\mbox{$\b L^2_{sym}(\Omega)=\b L^2_{\sigma,sym}(\Omega)\oplus \b G_{sym}^2(\Omega)$ \quad (orthogonal decomposition)}, 
\]
where   
\begin{align*}
& \b L^2_{\sigma,sym}(\Omega)
=\{\bw \in\b L_{sym}^2(\Omega);\,\div\bw=0\mbox{ in }\Omega,\,\bw\cdot\b n|_{\del\Omega}=0 \}, 
\\
& \b G_{sym}^2(\Omega)=\{\nabla \phi; \phi \in H_{*,sym}^1(\Omega)\}.
\end{align*}
We denote by $\mathbb{P}_\sigma$ the orthogonal projection on $\b L^2_{\sigma,sym}(\Omega)$. 

We define the function spaces $\b X$, $\b Y$, $\mathbb{X}_{\sigma}$, $\mathbb{Y}_{\sigma}$, $X$ and $Y$ by  
\begin{align*}
\b X &= \b L^2_{sym}(\Omega) \times L_{sym}^2(\Omega)\times L_{sym}^2(\Omega), 
\\
\b Y & =[\b H_{b}^2(\Omega)\times (H^2\cap H^1_0)(\Omega)\times (H^2\cap H^1_0)(\Omega))]\cap \b X, 
\\ 
\mathbb{X}_\sigma &= \b L^2_{\sigma,sym}(\Omega) \times L_{sym}^2(\Omega)\times L_{sym}^2(\Omega), 
\\ 
\mathbb{Y}_\sigma & =[\b H_{b}^2(\Omega)\times (H^2\cap H^1_0)(\Omega)\times (H^2\cap H^1_0)(\Omega))]\cap \mathbb{X}_\sigma, 
\\
X & =H_{*,sym}^1(\Omega)\times \b L^2_{sym}(\Omega) \times L_{sym}^2(\Omega)\times L_{sym}^2(\Omega),  
\\
Y & =[H_{*,sym}^1(\Omega)\times \b H_{b}^2(\Omega)\times (H^2\cap H^1_0)(\Omega)\times (H^2\cap H^1_0)(\Omega))]\cap X,
\end{align*}
respectively, where 
\[
\b H_b^2(\Omega) = \left\{\bw = \trans(w^1,w^2)\in \b H^2(\Omega);\,\frac{\del w^1}{\del x_2} = w^2 = 0\mbox{ on }\{x_2 = 0,1\}\right\}
\]
\[
H_0^1(\Omega) = \{\theta\in H^1(\Omega);\,\theta|_{x_2 = 0,1} = 0\}.
\]
We will also employ the space $\b X^1$ defined by 
\[
\b X^1=[\b H_b^1(\Omega)\times H_0^1(\Omega)\times H_0^1(\Omega)]\cap \b X 
\]
with norm
\[
\|\bu\|_{\b X^1}=\sqrt{\Pr^{-1}\|\nabla\bw\|_2^2+\|\nabla\theta\|_2^2+\|\nabla\psi\|_2^2}
\]
for $\bu=\trans(\bw,\theta,\psi) \in \b X^1$. Here 
\[
\b H_b^1(\Omega) = \left\{\bw = \trans(w^1,w^2)\in \b H^1(\Omega);\, w^2|_{x_2 = 0, 1} = 0\right\}.
\]
Observe that the Poincar\'e inequality $\|\bw\|_2\le C\|\nabla\bw\|_2$ holds for $\bw=\trans(w^1,w^2) \in \b H_b^1(\Omega)\cap \b L^2_{sym}(\Omega)$, which follows from the fact $\int_\Omega w^1\,dx=0$ due to the oddness of $w^1$ in $x_1$ and the fact $w^2|_{x_2 = 0, 1} = 0$. 

We next introduce $\|\bu\|_{(\b X^1)^*}$. We expand $\bu=\trans(\bw,\theta,\psi) \in \b X$ with $\bw=\trans(w^1,w^2)$ in the Fourier series as 
\begin{align*}
w^1 & =\sum_{j\ge1, k\ge0}w_{jk}^1\sin{\alpha j x_1}\cos{k\pi x_2}, 
\\
w^2 & =\sum_{j\ge0, k\ge1}w_{jk}^2\cos{\alpha j x_1}\sin{k\pi x_2}, 
\\
\theta & =\sum_{j\ge0, k\ge1}\theta_{jk}\cos{\alpha j x_1}\sin{k\pi x_2}, 
\\
\psi & =\sum_{j\ge0, k\ge1}\psi_{jk}\cos{\alpha j x_1}\sin{k\pi x_2},
\end{align*}
and define $\|\bu\|_{(\b X^1)^*}$ by 
\[
\|u\|_{(\b X^1)^*}=\left[\sum_{j\ge0,k\ge0}(\alpha^2 j^2+k^2\pi^2)^{-1}\left(\Pr^{-1}\left\{(w_{jk}^1)^2+(w_{jk}^2)^2\right\}+\theta_{jk}^2+\psi_{jk}^2 \right) \right]^{\frac{1}{2}},
\]
where $w_{0k}^1=w_{j0}^2=\theta_{j0}=\psi_{j0}=0$ for $k,j\ge0$. 
It is easily verified that 
\[
\| \bu \|_{(\b X^1)^*} \le C \| \bu \|_2
\]
and 
\[
|(\bu_1,\bu_2)| \le \|\bu_1\|_{\b X^1}\|\bu_2\|_{(\b X^1)^*}.
\]
We set
\[
X^1 = H_{*,sym}^1(\Omega) \times \b X^1.
\] 

We next introduce function spaces on time intervals. Let $t_{1} < t_{2}$. We define $\mathcal{X}(t_{1}, t_{2})$ and $\mathcal{Y}(t_{1}, t_{2})$ by  
\begin{align*}
\mathcal{X} (t_{1}, t_{2}) & =L^2(t_{1}, t_{2}; X), 
\\
\mathcal{Y} (t_{1}, t_{2})& =L^2(t_{1}, t_{2},T;Y)\cap H^1(t_{1}, t_{2},T; X), 
\end{align*}
with norms $\|u\|_{\mathcal{X}(t_{1}, t_{2})}$ and $\|u\|_{\mathcal{Y}(t_{1}, t_{2})}$, respectively. 

We define the inner product $\lang \bu_1,\bu_2\rang$ by 
\[
\langle\bu_1,\bu_2\rangle = \frac{a}{2\pi}\int_0^{\frac{2\pi}{a}}(\bu_1(t),\bu_2(t))\,dt  
\]
for $\bu_j=\trans(\bw_j,\theta_j,\psi_j)\in L^2(\mathbb{T}_{\frac{2\pi}{a}};\b L^2(\Omega)\times L^2(\Omega)\times L^2(\Omega))$ $(j=1,2)$

Let $\ep$ be a given positive number. For $u_j=\trans(\phi_j,\bu_j)\in L^2(\Omega)\times \b L^2(\Omega)\times L^2(\Omega)\times L^2(\Omega)$ $(j=1,2)$, we define the inner product $(u_1,u_2)_\ep$ and the norm $|||u|||_\ep$ with $\ep$-weights by
\[
(u_1,u_2)_\ep = \ep^2(\phi_1,\phi_2)_{L^2} + (\bu_1,\bu_2)
\]
for $u_j=\trans(\phi_j, \bu_j )\in L^2(\Omega)\times \b L^2(\Omega)\times L^2(\Omega)\times L^2(\Omega)$ $(j = 1,2)$ and  
\[
|||u|||_\ep=\sqrt{(u,u)_\ep}=\sqrt{\ep^2\|\phi\|_2^2+\|\bu\|_2^2}
\]
for $u=\trans(\phi,\bu)\in L^2(\Omega)\times \b L^2(\Omega)\times L^2(\Omega)\times L^2(\Omega)$, 
and, likewise, for $\bu_j=\trans(\bw_j,\theta_j,\psi_j)\in L^2(\mathbb{T}_{\frac{2\pi}{a}};L^2(\Omega)\times \b L^2(\Omega)\times L^2(\Omega)\times L^2(\Omega))$ $(j=1,2)$, the inner product $\langle u_1,u_2\rangle_\ep$ is defined by  
\[
\langle u_1,u_2\rangle_\ep = \frac{a}{2\pi}\int_0^{\frac{2\pi}{a}}(u_1(t),u_2(t))_\ep\,dt. 
\] 
We shall also employ the norms $||| u |||_{\ep, X^1}$, $|||u |||_{\ep, \mathcal{X}(t_{1}, t_{2})}$ and $|||u |||_{\ep, \mathcal{Y}(t_{1}, t_{2})}$ of $u=\trans(\phi, \bu)$ with $\ep$ weights defined by 
\begin{align*}
||| u |||_{\ep,X^1} & = \left\{ |||u|||_\ep^2 + \ep^2 |||\delx u|||_\ep^2 \right\}^{\frac{1}{2}},
\\
||| u |||_{\ep, \mathcal{X}(t_{1}, t_{2})}  & = \left\{\int_{t_{1} }^{t_{2} } \left( \ep^2 \| \phi \|_2^2 + \|\bu\|_{(\b X^1)^*}^2 + \ep^6 \|\delx \phi \|_2^2 +\ep^2 \|\bu\|_2^2\right) \, dt \right\}^{\frac{1}{2}}, 
\\
||| u |||_{\ep, \mathcal{Y}(t_{1}, t_{2})} & = \left\{ \sup_{t_{1} \le t \le t_{2} } ||| u(t) |||_{\ep,X^1}^2 \right. \\
& \quad \quad \left. + \int_{t_{1} }^{t_{2} } \left(|||\delx u |||_\ep^2 + \ep ^2 |||\delt u |||_\ep^2 + \ep^2 \|\delx ^2\bu \|_2^2 + \ep^6 \|\delx\delt \phi\|_2^2 \right) \,dt \right\}^{\frac{1}{2}}. 
\end{align*}

Time periodic function spaces $\mathcal{X}_a$ and $\mathcal{Y}_a$ with period $\frac{2\pi}{a}$ are defined by 
\begin{align*}
\mathcal{X}_a & =L^2(\mathbb{T}_{\frac{2\pi}{a}};X), 
\\
\mathcal{Y}_a & =L^2(\mathbb{T}_{\frac{2\pi}{a}};Y)\cap H^1(\mathbb{T}_{\frac{2\pi}{a}};X), 
\end{align*}
respectively. We denote the norm of $u\in \mathcal{X}_a$ (resp. $u\in \mathcal{Y}_a$) by $\|u\|_{\mathcal{X}_a} = \|u\|_{\mathcal{X}(0, \frac{2\pi}{a})}$ (resp. $\|u\|_{\mathcal{Y}_a} = \|u\|_{\mathcal{Y}(0, \frac{2\pi}{a})}$), and likewise, the norm of $u\in \mathcal{X}_a$ (resp. $u\in \mathcal{Y}_a$) by $||| u |||_{\ep, \mathcal{X}_a} = ||| u |||_{\ep, \mathcal{X}(0, \frac{2\pi}{a})}$ (resp. $||| u |||_{\ep, \mathcal{Y}_a} = ||| u |||_{\ep, \mathcal{Y}(0, \frac{2\pi}{a})}$).

The resolvent set and spectrum of an operator $A$ are denoted by $\rho(A)$ and $\sigma(A)$, respectively. We denote by $\mathfrak{B}(E_{1}, E_{2})$ the space of all bounded linear operators from $E_{1}$ to $E_{2}$.

\section{Hopf bifurcation in the incompressible system}

In this section we summarize the results on the Hopf bifurcation for the incompressible system \eqref{1.1} obtained in \cite{Bona-Hsia-Ma-Wang} and introduce the associated pressure of the velocity field of the bifurcating time periodic solution. 

Setting $\ep = 0$ in \eqref{1.1}--\eqref{1.4}, we obtain the incompressible system: 
\begin{equation}\label{3.4}
\left\{
\begin{array}{rcl}
\div\bw &=& 0,\\
\delt\bw-\Pr\Delta\bw + \Pr\nabla \phi -\Pr\mathcal{R}_1\theta\b{e}_2 + \Pr\mathcal{R}_2\psi\b{e}_2 +\bw\cdot\nabla\bw &=& \b{0},\\
\delt\theta - \Delta\theta - \mathcal{R}_1\bw\cdot\b e_2 + \bw\cdot\nabla\theta &=& 0,\\
\delt \psi - d \Delta \psi - \mathcal{R}_2\bw\cdot\b e_2 + \bw\cdot\nabla \psi &=& 0.
\end{array}
\right.
\end{equation}
The boundary condition on the boundary $\{x_2=0,1\}$ is given as 
\begin{equation}\label{3.5}
\mbox{$\ds \frac{\del w^1}{\del x_2} = w^2 = \theta = \psi = 0$ on $\{x_2 = 0,1\}$.}
\end{equation}
The problem \eqref{3.4}--\eqref{3.5} then has a trivial stationary solution $u=0$ which corresponds to the the motionless state. See \cite{Bona-Hsia-Ma-Wang, Hsia-Kagei-Nishida-Teramoto1} for the derivation of the non-dimensional perturbation equations \eqref{3.4}. 

The result on the Hopf bifurcation for \eqref{3.4}--\eqref{3.5} by Bona, Hsia, Ma and Wang \cite{Bona-Hsia-Ma-Wang} is summarized as follows. There are positive numbers $\mathcal{R}_{2*}$ and $\mathcal{R}_2^*$ such that if $\mathcal{R}_2\in[\mathcal{R}_{2*},\mathcal{R}_2^*)$, $\Pr>1$ and $0<d<1$, then there exists a critical number $\mathcal{R}_{1,c}$ such that the basic state $u=0$ is stable when $\mathcal{R}_1<\mathcal{R}_{1,c}$; and if $\mathcal{R}_1>\mathcal{R}_{1,c}$, then $u=0$ loses its stability and time periodic solutions of the incompressible problem \eqref{3.4}--\eqref{3.5} bifurcate from $u=0$. 

To state the bifurcating result in \cite{Bona-Hsia-Ma-Wang} more precisely, we take $\eta = \mathcal{R}_1-\mathcal{R}_{1,c}$ as a bifurcation parameter. 
The linearized operator $\mathbb{L}_{\mathcal{R}_{1,c}+\eta}$ on $\mathbb{X}_\sigma$ is then given by 
\[
D(\mathbb{L}_{\mathcal{R}_{1,c}+\eta})=\mathbb{Y}_\sigma,  
\]
\[
\mathbb{L}_{\mathcal{R}_{1,c}+\eta} = 
\begin{pmatrix}
-\Pr\mathbb{P}_\sigma\Delta & -\Pr(\mathcal{R}_{1,c} + \eta)\mathbb{P}_\sigma\boldsymbol{e}_2 & \Pr\mathcal{R}_2\mathbb{P}_\sigma\boldsymbol{e}_2 \\
-(\mathcal{R}_{1,c} + \eta)\trans\boldsymbol{e}_2  & -\Delta & 0 \\
-\mathcal{R}_2\trans\boldsymbol{e}_2  & 0 & -d\Delta
\end{pmatrix}
. 
\]

We write problem \eqref{3.4}--\eqref{3.5} in the form 
\begin{equation}\label{3.6}
\delt\bu+\mathbb{L}_{\mathcal{R}_{1,c}}\bu+\eta \mathbb{P}\b K\b u+\mathbb{P}\b N(\bu)=\b 0,  
\end{equation}
where $\bu=\trans(\bw,\theta,\psi)$, 
\[
\mathbb{P}=
\begin{pmatrix}
\mathbb{P}_\sigma & \b0 & \b0\\
\trans\b0 & 1 & 0 \\
\trans\b0 & 0 & 1
\end{pmatrix}, 
\ \ \ 
\b K
=
\begin{pmatrix}
\b O & -\Pr\,\b e_2 & \b0 \\
-\trans\b e_2 & 0 & 0 \\
\trans\b0 & 0 & 0 
\end{pmatrix}, 
\]
and 
\[
\b N(\bu) = \b N(\bu,\bu)
\]
with  
\[
\mbox{
$\ds 
\b N(\bu_1,\bu_2)
=
\begin{pmatrix}
\bw_1\cdot\nabla\bw_2 \\
\bw_1\cdot\nabla\theta_2 \\
\bw_1\cdot\nabla\psi_2
\end{pmatrix}
$ 
for $\bu_j
=
\begin{pmatrix}
\bw_j \\
\theta_j \\
\psi_j
\end{pmatrix}
$ 
$(j=1,2)$.}
\] 

The adjoint operator $\mathbb{L}_{\mathcal{R}_{1,c}+\eta}^*$ on $\mathbb{X}_\sigma$ is given by
\[
\mathbb{L}_{\mathcal{R}_{1,c}+\eta}^* = 
\begin{pmatrix}
-\Pr\mathbb{P}_\sigma\Delta & -\Pr(\mathcal{R}_{1,c} + \eta)\mathbb{P}_\sigma\boldsymbol{e}_2 & -\Pr\mathcal{R}_2\mathbb{P}_\sigma\boldsymbol{e}_2\\
-(\mathcal{R}_{1,c} + \eta)\boldsymbol{e}_2  & -\Delta & 0\\
\mathcal{R}_2\boldsymbol{e}_2 & 0 & -d \Delta
\end{pmatrix}
\]
with domain $D(\mathbb{L}_{\mathcal{R}_{1,c}+\eta}^*) = D(\mathbb{L}_{\mathcal{R}_{1,c}+\eta})=\mathbb{Y}_\sigma$. 

Bona, Hsia, Ma and Wang in \cite{Bona-Hsia-Ma-Wang} proved the following result on the spectrum of $\mathbb{L}_{\mathcal{R}_{1,c}+\eta}$.

\begin{prop}\label{prop3.1} {\rm (\cite{Bona-Hsia-Ma-Wang})} 
{\rm (i)} There exist positive numbers $\mathcal{R}_{2*}$ and $\mathcal{R}_2^*$ such that if $\mathcal{R}_2\in[\mathcal{R}_{2*},\mathcal{R}_2^*)$, $\Pr>1$ and $0<d<1$, then the following assertions hold. 
There exist positive constants $\eta_0$, $b_0$ and $\Lambda_0$ such that if $|\eta|\le\eta_0$ then it holds that 
\[
\Sigma\setminus\{\lambda_+(\eta),\lambda_-(\eta)\}\subset\rho(-\mathbb{L}_{\mathcal{R}_{1,c}+\eta}),
\]
where $\Sigma=\{\lambda\in\mathbb{C};\,\re\lambda\ge-b_0|\im\lambda|^2-\Lambda_0\}$; $\lambda_+(\eta)$ and $\lambda_-(\eta)$ are simple eigenvalues of $-\mathbb{L}_{\mathcal{R}_{1,c}+\eta}$ satisfying $\lambda_-(\eta) = \overline{\lambda_+(\eta)}$ and
\[
\lambda_+(0) = ia, \ \ \ \frac{d\re\lambda_+}{d\eta}(0) >0. 
\]
Here $a$ is a positive constant. 

\vspace{1ex} 
{\rm (ii)} Let $\boldsymbol{u}_{\pm}$ be eigenfunctions for the eigenvalues $\pm ia$ of $-\mathbb{L}_{\mathcal{R}_{1,c}}$ and let $\boldsymbol{u}_{\pm}^*$ be eigenfunctions for the eigenvalues $\mp ia$ of the adjoint operator $-\mathbb{L}_{\mathcal{R}_{1,c}}^*$ satisfying $(\boldsymbol{u}_{j}, \boldsymbol{u}_{k}^*) = \delta_{jk}$, where $j, k \in \{+, - \}$. Then $\b{u} = \overline{\b{u}_{+} }$ and $\b{u}^{*} = \overline{\b{u}^{*}_{+} }$; and the eigenprojections $\b P_\pm$ for the eigenvalues $\pm ia$ of $-\mathbb{L}_{\mathcal{R}_{1,c}}$ are given by
\[
\b P_\pm\boldsymbol{u} = (\boldsymbol{u},\boldsymbol{u}_\pm^*)\boldsymbol{u}_\pm.
\]
Furthermore, it holds that  
\[
\frac{d\re\lambda_+}{d\eta}(0) = -\re(\b K\boldsymbol{u}_+,\boldsymbol{u}_+^*)>0.
\]
\end{prop}

We next define the operators $\mathbb{B}$ and $\mathbb{B}^{*}$ on $L^2(\mathbb{T}_{\frac{2\pi}{a}};\mathbb{X}_\sigma)$ by 
\[
\mathbb{B}=\delt+\mathbb{L}_{\mathcal{R}_{1,c}} \quad  \mbox{ \rm and} \quad  \mathbb{B}^*=-\delt+\mathbb{L}_{\mathcal{R}_{1,c}}^*
\]
with domain $D(\mathbb{B})=D(\mathbb{B}^*)=L^2(\mathbb{T}_{\frac{2\pi}{a}};D(\mathbb{L}_0))\cap H^1(\mathbb{T}_{\frac{2\pi}{a}};\mathbb{X}_\sigma)$. 

One can then see by a direct computation that the functions 
\[
\boldsymbol{z}_{\pm} =  e^{\pm iat}\boldsymbol{u}_{\pm},\,\boldsymbol{z}_{\pm}^* =  e^{\pm iat}\boldsymbol{u}_{\pm}^* 
\]
are eigenfunctions for the eigenvalue $0$ of $\mathbb{B}$ and $\mathbb{B}^{*}$ with properties   
\[
\langle \boldsymbol{z}_{\pm},\boldsymbol{z}_{\pm}^*\rangle = 1,\,\langle\boldsymbol{z}_{\mp},\boldsymbol{z}_{\pm}^*\rangle = 0.
\]
We set 
\[
[\boldsymbol{z}]_\pm = \langle\boldsymbol{z},\boldsymbol{z}_\pm^*\rangle.
\]
The operators $\hat{\mathcal{P}}_{\pm}$ defined by
\[
\hat{\mathcal{P}}_\pm\boldsymbol{z} = [\boldsymbol{z}]_\pm\boldsymbol{z}_\pm \quad (\boldsymbol{z}\in L^2(\mathbb{T}_{\frac{2\pi}{a}};\mathbb{X}_\sigma)) 
\]
are then projections which satisfy $\hat{\mathcal{P}}_j \hat{\mathcal{P}}_k = \delta_{jk} \hat{ \mathcal{P} }_{j}$ for $j,k\in\{+,-\}$. As for the operator $\mathbb{B}$, we then have the following proposition. 

\begin{prop}\label{prop3.2}
$0$ is a semisimple eigenvalue of $\mathbb{B}$ and it holds that 
\[
L^2(\mathbb{T}_{\frac{2\pi}{a}};\mathbb{X}_\sigma) = \Ker (\mathbb{B})\oplus R(\mathbb{B}). 
\]
Set  
\[
\hat{\mathcal{P}}_0 = \hat{\mathcal{P}}_+ + \hat{\mathcal{P}}_-, \ \ \ \hat{\mathcal{Q}}_0 = I - \hat{\mathcal{P}}_0. 
\]
Then $\hat{\mathcal{P}}_0$ is an eigenprojection for the eigenvalue $0$ of $\mathbb{B}$ and $\hat{ \mathcal{\mathcal{Q}} }_0$ is a projection on $R(\mathbb{B})$ along $N(\mathbb{B})$. There holds that $\bu \in R(\mathbb{B})$ if and only if $\hat{\mathcal{P}}_0 \bu = 0$, i.e., $[\b z]_+ = [\b z]_- = 0$. 
\end{prop}

We note that if $\bu$ is a real valued function, then $[\bu]_- = \overline{[\bu]_+}$, and hence,  
\[
\hat{\mathcal{P}}_0 \bu = 2\re([\bu]_+ \b z_+). 
\]

Based on Proposition \ref{prop3.1}, a Hopf bifurcation would occur when $\mathcal{R}_1> \mathcal{R}_{1,c}$. It is expected that \eqref{3.6} has a nontrivial time periodic solution of period $\frac{2\pi}{a(1+\omega)}$ with small $\omega$ for sufficiently small $\eta$. We thus change the variable $t \mapsto (1+\omega)t$. The problem is then reduced to finding a nontrivial time periodic solution of period $\frac{2\pi}{a}$ to the equation 
\begin{equation}\label{3.7}
(1 + \omega) \delt\bu+\mathbb{L}_{\mathcal{R}_{1,c}}\bu+\eta \mathbb{P}\b K\b u+\mathbb{P}\b N(\bu)=\b 0. 
\end{equation}

The following result on the Hopf bifurcation was shown by Bona, Hsia, Ma and Wang in \cite{Bona-Hsia-Ma-Wang}. 

\begin{prop}\label{prop3.3} {\rm (\cite{Bona-Hsia-Ma-Wang})} 
{\rm (i)} There exists a positive constant $\delta_0$ such that \eqref{3.7} has a nontrivial time periodic solution $\bu_\delta \in L^2(\mathbb{T}_{\frac{2\pi}{a}}; \mathbb{Y}_{\sigma} ) \cap H^1(\mathbb{T}_{\frac{2\pi}{a}};\mathbb{X}_\sigma)$ for $\eta = \eta_{\delta}$ and $\omega = \omega_{\delta}$, where $\{\eta_{\delta}, \omega_{\delta}, \bu_\delta\}$ with a parameter $\delta$ takes the form 
\begin{align*}
\eta_{\delta} & = \tilde\eta_0 \delta^2 + \tilde\eta_1(\delta) \delta^3, \\
\omega_{\delta} & = \tilde \omega_0 \delta^2 + \tilde\omega_1(\delta) \delta^3, \\
\bu_\delta & = \delta(\b z_0+\delta\bU_\delta)
\end{align*}
for $|\delta|\le \delta_0$. Here $\tilde \eta_0$ and $\tilde\omega_0$ are constants and $\tilde \eta_0 > 0$; $\tilde \eta_1(\delta)$ and $\tilde \omega_1(\delta)$ are analytic in $\delta$ satisfying $\tilde \eta_1(\delta)=O(1)$ and $\tilde \omega_1(\delta)=O(1)$ as $\delta \to 0$; $\b z_0=\trans(\bw_{0}, \theta_{0}, \psi_{0}) = 2\re( \b{z}_{+} )$; and $\bU_\delta$ is in $R(\mathbb{B})\cap [L^2(\mathbb{T}_{\frac{2\pi}{a}}; \mathbb{Y}_{\sigma} ) \cap H^1(\mathbb{T}_{\frac{2\pi}{a}};\mathbb{X}_\sigma)]$ 
and is analytic in $\delta$. 

\vspace{1ex} 
{\rm (ii)} There exists a neighborhood $\b{O}_{0}$ of $\trans(\eta,\omega,\bu) = \trans(0,0,\b{0})$ in $\mathbb{R}\times\mathbb{R}\times [ L^2(\mathbb{T}_{\frac{2\pi}{a}}; \mathbb{Y}_{\sigma} ) \cap H^1(\mathbb{T}_{\frac{2\pi}{a}};\mathbb{X}_\sigma) ]$ such that if $\trans(\eta,\omega,\bu) \in \b{O}_{0}$ is a solution of \eqref{3.7}, then
\[
\ts 
\trans(\eta,\omega,\bu) \in \{\trans(\eta,\omega, \b{0})\in \b{O}_{0} \}\cup\{\trans(\eta_\delta,\omega_\delta, \bu_\delta(\cdot + \tau) ); |\delta| \le \delta_1, -\frac{\pi}{a} \le \tau < \frac{\pi}{a} \}
\]
where $\{\eta_\delta,\omega_\delta,\bu_\delta\}$ is the solution branch obtained in {\rm (i)}.
\end{prop}

We next introduce the associated pressure $\phi_{\delta}$ of the time periodic bifurcating solution $\bu_{\delta}$ of \eqref{3.7}. It is known \cite{Galdi, Sohr, Temam3} that there exists the associated pressure $\phi_\delta\in L^2(\mathbb{T}_{\frac{2\pi}{a}};H^1_{*,sym}(\Omega))$ of $\bu_\delta=\trans(\bw_\delta,\theta_\delta,\psi_\delta)$, namely, there exists a unique $\phi_\delta\in L^2(\mathbb{T}_{\frac{2\pi}{a}};H^1_{*,sym}(\Omega))$ such that $u_{\delta} = \trans(\phi_{\delta}, \bu_{\delta}) \in \mathcal{Y}_{a}$ is a time periodic solution of 
\begin{equation}\label{3.8} 
(1 + \omega) \delt 
\begin{pmatrix}
0 \\
\bu
\end{pmatrix}
+ L^{\ep}_{\mathcal{R}_{1,c} } u + \eta K u + N(u) = 0 
\end{equation} 
for $\eta = \eta_{\delta}$ and $\omega = \omega_{\delta}$. 
Here $L^{\ep}_{\mathcal{R}_1}$, $K$ and $N(\,\cdot\,)$ are the maps defined as follows;  $L^{\ep}_{\mathcal{R}_1}$ is the linearized operator around $u = 0$ on the space $X$ with domain 
\[
D(L^{\ep}_{\mathcal{R}_1})=Y 
\]
and it is given in the form
\[
L^{\ep}_{\mathcal{R}_1} =
\begin{pmatrix}
0 & \frac{1}{\ep^2}\div & 0 & 0\\
\Pr\nabla & -\Pr\Delta & -\Pr\,\mathcal{R}_{1}\boldsymbol{e}_2 & \Pr\mathcal{R}_2\boldsymbol{e}_2\\
0 & -\mathcal{R}_{1}\trans\boldsymbol{e}_2  & -\Delta & 0\\
0 &-\mathcal{R}_2\trans\boldsymbol{e}_2 & 0 & - d \Delta
\end{pmatrix};  
\]
$K$ is the linear operator given by 
\[
K=
\left(\begin{array}{@{\,}c|ccc@{\,}}
0 & \trans\b0 & 0 & 0 \\ \hline
\b 0 & & &  \\ 
0 & & \b K & \\
0 & & & 
\end{array}\right);  
\] 
and $N( \, \cdot\, )$ is the nonlinear map given by 
\[ 
N(u)= N(u,u),
\]
where 
\[ 
N(u_1,u_2) =
\begin{pmatrix}
0\\
\boldsymbol{N}(\bu_1,\bu_2)
\end{pmatrix}
\]  
for $u_j=\trans(\phi_j,\bu_j)$ $(j=1,2)$. 

The associated pressure $\phi_\delta$ takes the form 
\[
\phi_\delta=\delta(p_0 + \delta\Phi_\delta),
\]
where $p_0 = \re(e^{ia t}\phi_+)$ and $\Phi_\delta\in L^2(\mathbb{T}_{\frac{2\pi}{a}};H^1_{*,sym}(\Omega))$. Here $\phi_+\in H_{*,sym}^1(\Omega)$ is the associated pressure of $\bu_+=\trans(\bw_+,\theta_+,\psi_+)$, i.e., $u_{+} = \trans(\phi_+, \bu_{+})$ satisfies 
\[
ia 
\begin{pmatrix}
0 \\
\bu_{+} 
\end{pmatrix}
+ L^{\ep}_{\mathcal{R}_{1,c}} u_{+} = 0, 
\]
and $\Phi_{\delta}$ is the associated pressure of $\b{U}_{\delta} = \trans(\bW_{\delta}, \Theta_{\delta}, \Psi_{\delta})$, i.e., $U_{\delta} = \trans(\Phi_{\delta}, \b{U}_{\delta})$  satisfies 
\[
\delt 
\begin{pmatrix}
0 \\
\b{U}_{\delta}
\end{pmatrix} 
+ \delta \tilde\omega_{\delta} 
\delt
\begin{pmatrix}
0 \\ 
\b{z}_{0} + \delta \b{U}_{\delta}
\end{pmatrix}
+ L^{\ep}_{\mathcal{R}_{1,c} } U_{\delta} + \delta \tilde\eta_{\delta} K (z_{0} + \delta U_{\delta}) + N(z_{0} + \delta U_{\delta}) = 0,  
\]
where $z_{0} = \trans(p_{0}, \b{z}_{0})$, $\tilde{\eta}_{\delta} = \delta^{-2} \eta_{\delta}$ and $\tilde{\omega}_{\delta} = \delta^{-2} \omega_{\delta}$. Furthermore, $\Phi_\delta$ is in $H^{1}(\mathbb{T}_{\frac{2\pi}{a}};H^1_{*,sym}(\Omega))$ and satisfies $\| \Phi_\delta \|_{H^{1}(\mathbb{T}_{\frac{2\pi}{a}};H^1_{*}(\Omega)) } \le C$ uniformly for $|\delta| \le \delta_{0}$ if $\delta_{0}$ is suitably small. 

The above decomposition of $\phi_{\delta}$ into $p_{0}$ and $\Phi_{\delta}$ parts can be seen by decomposing \eqref{3.7} into $\hat{\mathcal{P} }_{0}$ and $\hat{\mathcal{Q} }_{0}$ parts. To see the boundedness of $\{\Phi_{\delta} \}_{|\delta| \le \delta_{0} }$, we first observe that $\b U_\delta$ is a time periodic solution of  
\[
(1+\delta^{2} \tilde \omega_{\delta} ) \delt \b U_\delta + \mathbb{L}_{\mathcal{R}_{1,c}} \b U_\delta + \delta \tilde \eta_{\delta} \hat{\mathcal{Q}}_0 \mathbb{P}\b K(\b z_0 + \delta \b U_\delta) + \hat{\mathcal{Q}}_0 \mathbb{P} \b N(\b z_0 + \delta \b U_\delta) = \b 0.
\]
Since $\b z_0$ is a smooth function and $\{\b U_\delta\}_{|\delta|\le \delta_0}$ is bounded in $L^2(\mathbb{T}_{\frac{2\pi}{a}}; \mathbb{Y}_{\sigma} ) \cap H^1(\mathbb{T}_{\frac{2\pi}{a}};\mathbb{X}_\sigma)$, we see that $\{\delt \b U_\delta\}_{|\delta|\le \delta_0}$ is bounded in $L^2(\mathbb{T}_{\frac{2\pi}{a}}; \mathbb{Y}_{\sigma} ) \cap H^1(\mathbb{T}_{\frac{2\pi}{a}};\mathbb{X}_\sigma)$ if $\delta_{0}$ is sufficiently small. It then follows that $\{\Phi_\delta \}_{|\delta|\le \delta_0}$ is bounded in $H^{1} (\mathbb{T}_{\frac{2\pi}{a}};H^1_*(\Omega))$. 

In summary we have the following proposition. 

\begin{prop}\label{prop3.4} 
{\rm (i)} There exists a positive constant $\delta_0$ such that \eqref{3.8} has a nontrivial time periodic solution $u_\delta \in \mathcal{Y}_{a}$ for $\eta = \eta_{\delta}$ and $\omega = \omega_{\delta}$, where $\{\eta_{\delta}, \omega_{\delta}, u_\delta\}$ with a parameter $\delta$ takes the form 
\begin{align*}
\eta_{\delta} & = \tilde\eta_0 \delta^2 + \tilde\eta_1(\delta) \delta^3, \\
\omega_{\delta} & = \tilde \omega_0 \delta^2 + \tilde\omega_1(\delta) \delta^3, \\
u_\delta & = \delta(z_0+\delta U_\delta)
\end{align*}
for $|\delta|\le \delta_0$. Here $\eta_{\delta}$ and $\omega_{\delta}$ are the numbers given in {\rm Proposition \ref{prop3.3} } and $U_\delta \in \mathcal{Y}_a$ is analytic in $\delta$. Furthermore, $\| U_\delta \|_{\mathcal{Y}_a} \le C$ uniformly for $|\delta| \le \delta_0$. 

\vspace{1ex} 
{\rm (ii)} There exists a neighborhood $\mathcal{O}_{0}$ of $\trans(\eta,\omega, u) = \trans(0,0, 0)$ in $\mathbb{R}\times\mathbb{R}\times \mathcal{Y}_{a}$ such that if $\trans(\eta,\omega, u) \in \mathcal{O}_{0}$ is a solution of \eqref{3.8}, then
\[
\ts 
\trans(\eta,\omega, u) \in \{ \trans(\eta,\omega, 0) \in \mathcal{O}_{0} \}\cup\{\trans(\eta_\delta,\omega_\delta, u_\delta(\cdot + \tau) ); |\delta| \le \delta_1, -\frac{\pi}{a} \le \tau < \frac{\pi}{a} \}
\]
where $\{\eta_\delta,\omega_\delta, u_\delta\}$ is the solution branch given in {\rm (i)}.
\end{prop}  

We close this section by introducing notation related to Proposition \ref{prop3.4}. The associated pressure $p_{-}$ of $\bu_{-}=\trans(\bw_-,\theta_-,\psi_-)$ is given by $\phi_{-}=\ov{\phi_+}$, and it holds that $u_{-} = \trans(\phi_{-}, \bu_{-})$ satisfies 
\[
-ia 
\begin{pmatrix}
0 \\
\bu_{-} 
\end{pmatrix}
+ L^{\ep *}_{\mathcal{R}_{1,c}} u_{-} = 0.  
\]
where $L^{\ep*}_{\mathcal{R}_{1} }: X \to X$ is the adjoint operator of $L^{\ep}_{\mathcal{R}_{1} }$ and it is defined by 
\[
D(L^{\ep*}_{\mathcal{R}_1})=Y,
\]
\[
L^{\ep*}_{\mathcal{R}_1} =
\begin{pmatrix}
0 & - \frac{1}{\ep^2}\div & 0 & 0\\
- \Pr\nabla & -\Pr\Delta & -\Pr\mathcal{R}_{1}\boldsymbol{e}_2 & - \Pr\mathcal{R}_2\boldsymbol{e}_2 \\
0 & -\mathcal{R}_{1}\trans\boldsymbol{e}_2  & -\Delta & 0\\
0 & \mathcal{R}_2\trans\boldsymbol{e}_2 & 0 & -d \Delta
\end{pmatrix}.  
\]
Similarly, we have the associated pressures $\phi_\pm^{*}$ of $\bu_{\pm}^{*}=\trans(\bw_{\pm}^{*},\theta_{\pm}^{*},\psi_{\pm}^{*})$, and it holds that $u_{\pm}^{*} = \trans(\phi_{\pm}^{*}, \bu_{\pm}^{*})$ satisfies 
\[
\mp ia 
\begin{pmatrix}
0 \\
\bu_{\pm}^{*} 
\end{pmatrix}
+ L^{\ep*}_{\mathcal{R}_{1,c}} u_{\pm}^{*} = 0.  
\]

We define functions $u_{\pm}$, $u^{*}_{\pm}$, $z_{\pm}$ and $z^{*}_{\pm}$ by  
\begin{equation}\label{3.9}
u_{\pm} = \trans(\phi_{\pm},\bu_{\pm}), \quad u_{\pm}^{*}=\trans(\phi_{\pm}^{*},\bu_{\pm}^{*}), \quad z_{\pm} = e^{\pm iat} u_{\pm}, \quad z_{\pm}^{*} = e^{\pm iat} u_{\pm}^{*},  
\end{equation}
respectively. We finally define the operators $\mathcal{P}_{\pm} : \mathcal{X}_{a} \to \mathcal{X}_{a}$ by 
\[
\mathcal{P}_{\pm} u = [\bu ]_{\pm} z_{\pm} 
\]
for $u=\trans(\phi, \bu) \in \mathcal{X}_{a}$, and $\mathcal{P}_{0}, \, \mathcal{Q}_{0} : \mathcal{X}_{a} \to \mathcal{X}_{a}$ by 
\begin{equation}\label{3.10}
\mathcal{P}_{0} = \mathcal{P}_{+} + \mathcal{P}_{-}, \quad \mathcal{Q}_{0} = I - \mathcal{P}_{0}.
\end{equation} 
Observe that if $u = \trans(\phi, \bu)$ is real valued, then $\mathcal{P}_{0} u = 2 \re \left( [\bu]_{+} z_{+} \right)$. 

\section{Hopf bifurcation in the artificial compressible system}\label{Main results} 

In this section we state the result on the occurrence of a Hopf bifurcation bifurcation in the artificial compressible system \eqref{1.1}--\eqref{1.4} under \eqref{1.5} and the result on the singular limit $\ep \to 0$ for the time periodic bifurcating solutions. 

We fix the parameters $\Pr$, $d$ and $\mathcal{R}_{2}$ in such a way that these parameters satisfies the assumption of Proposition \ref{prop3.1}. 

The time periodic problem for the artificial compressible system \eqref{1.1}--\eqref{1.5} is formulated as 
\begin{equation}\label{4.3}
\delt u + L^{\ep}_{ \mathcal{R}_{1}} u + N(u) = 0.
\end{equation}
Here $u = \trans(\phi,\bu)$ with $\bu = \trans(\bw,\theta,\psi)$ is periodic in $t$,  and $L^{\ep}_{ \mathcal{R}_{1,c} }$, $K$ and $N(\cdot)$ are the maps given in \eqref{3.8}. 

As for the spectrum of the linearized operator $-L^{\ep}_{\mathcal{R}_1}$ for $\mathcal{R}_{1}$ close to the criticality $\mathcal{R}_{1,c}$, we obtained the following result in \cite{Hsia-Kagei-Nishida-Teramoto1}. (Cf., \cite{Kagei-Teramoto}.) 

\begin{thm}\label{thm4.1} {\rm (\cite{Hsia-Kagei-Nishida-Teramoto1})} 
{\rm (i)} There exist positive constants $\Lambda_1$, $\ep_1$ and $\eta_1$ such that 
for each $0<\ep\le \ep_1$ there exists a critical value $\mathcal{R}^{\ep}_{1,c}=\mathcal{R}_{1,c}+O(\ep^{2})$ such that  
if $|\eta|\le\eta_1$ with $\eta=\mathcal{R}_{1}-\mathcal{R}^{\ep}_{1,c}$, then
\[
\{\lambda\in\mathbb{C};\,\real\lambda\ge -\Lambda_1\} \setminus \{\lambda^{\ep}_+(\eta),\lambda^{\ep}_-(\eta)\}\subset\rho(-L^{\ep}_{\mathcal{R}_{1,c}^{\ep}+\eta }),
\]
where $\lambda^{\ep}_\pm(\eta)$ are simple eigenvalues of $-L^{\ep}_{\mathcal{R}_{1,c}^{\ep}+\eta}$ satisfying $\lambda^{\ep}_-(\eta)  =  \overline{\lambda^{\ep}_+(\eta)}$,  
\[
\lambda^{\ep}_{\pm}(\eta)  = \lambda_{\pm}(\eta)+O(\ep^2), 
\]
and
\[
\lambda^{\ep}_+(0) = ia^{\ep},
\ \ 
\frac{\partial\re\lambda^{\ep}_+}{\partial\eta}(0) 
= -\re(\b K\boldsymbol{u}_+,\boldsymbol{u}_+^*)+O(\ep^{2})
>0.
\]
Here $a^{\ep}$ is a constant satisfying $a^{\ep}=a+O(\ep^2)$. 

The eigenspaces for $\lambda_\pm^{\ep} (0)$ are spanned by $u^{\ep}_{\pm}$, respectively, where $u^{\ep}_{\pm}$ satisfy $\ov{ u^{\ep}_{-} } = u^{\ep}_{+}$ and $u^{\ep}_{\pm}=u_{\pm}+O(\ep^2)$.

\vspace{1ex}
{\rm (ii)} For $j=\pm$, the eigenprojections $P_j^\ep$ for $\lambda_j^\ep(0)$ $(j=\pm)$ satisfies  
\[
P_j^{\ep}u=(u,u_j^{\ep*})_{\ep} u_j^{\ep},
\]
where $u_{\pm}^{\ep}$ are eigenfunction for the eigenvalues $\lambda_{\pm}^\ep$ of $-L^{\ep}_{\mathcal{R}^{\ep}_{1,c}}$ satisfying $\ov{ u^{\ep}_{-} } = u^{\ep}_{+}$ and $u_{+}^{\ep}=u_{+} + O(\ep^{2})$ in $Y$; $u_{\pm}^{\ep*}$ are the adjoint eigenfunctions for the eigenvalues $\ov{\lambda_j^\ep(0)}$ of $-L^{\ep*}_{\mathcal{R}^{\ep}_{1,c}}$ satisfying $\ov{ u^{\ep*}_{-} } = u^{\ep*}_{+}$, $(u_j^\ep, u_k^{\ep*})_{\ep}=\delta_{jk}$ and $u_{+}^{\ep*}=u_{+}^{*}+O(\ep^{2})$ in $Y$; and $P^{0}_{\pm}$ are projections defined by $P^{0}_{\pm} u = (\bu,\bu_{\pm}^{*}) u_{\pm}$ for $u=\trans(\phi,\bu)$. 
Here $u_{\pm}$ and $u_{\pm}^{*}$ are functions defined in \eqref{3.9}. 
 
Furthermore, if $k\in \mathbb{Z}$ with $k\ge 0$ and $1\le p\le 2$, there exists a positive constant $\tilde{\ep}_{1} = \tilde{\ep}_{1}(k, p, \Omega)$ such that 
\begin{align*}
\|u_{\pm}^{\ep}\|_{H^{k}\times H^{k}} + \|u_{\pm}^{\ep*}\|_{H^{k}\times H^{k}}  & \le C, 
\\ 
\|u_{\pm}^{\ep} - u_{j} \|_{H^{k}\times H^{k}} + \|u_{\pm}^{\ep*} - u_{j}^{*}\|_{H^{k}\times H^{k}}  & \le C \ep^{2}, 
\\ 
\|P_{\pm}^{\ep}u\|_{H^{k}\times H^{k}} & \le C \|u\|_{L^{p}\times L^{p}},
\\
\|(P_{\pm}^{\ep}- P^{0}_{\pm})u\|_{H^{k}\times H^{k}} & \le C\ep^2 \|u\|_{L^{p}\times L^{p}}
\end{align*}
uniformly for $\ep \in (0, \tilde{\ep}_{1}]$.   
\end{thm}

By Theorem \ref{thm4.1}, a Hopf bifurcation is expected to occur when $\mathcal{R}_1$ passes the critical value $\mathcal{R}^{\ep}_{1,c}$. In fact, we have the following bifurcation result. 

To fix the time interval, we change the variables $t$ and $u$ into $\tilde{t}$ and $\tilde{u}$, respectively, by $\frac{a^\ep}{a}t = \tilde{t}$ and $u(x,t) = \tilde{u}(x,\tilde{t})$. Omitting the tildes $\tilde{}$ of $\tilde{t}$ and $\tilde{u}$, we see from \eqref{4.3} that the problem is transformed into   
\begin{equation}\label{4.4}
\frac{a^\ep}{a} \delt u + L^{\ep}_{\mathcal{R}_{1,c}^\ep} u + \eta K u + N(u) = 0.
\end{equation} 

As in the incompressible problem, we look for a nontrivial time periodic solution of \eqref{4.4} whose period is $\frac{2\pi}{a(1+\omega)}$ with small $\omega$ for sufficiently small $\ep$ and $\eta$. We thus again change the variable $t \mapsto (1+\omega)t$. 
The problem is then reduced to finding a nontrivial time periodic solution of period $\frac{2\pi}{a}$ of the following equation 
\begin{equation}\label{4.5}
\frac{a^{\ep}}{a} (1+\omega) \delt u + L^{\ep}_{\mathcal{R}_{1,c}^\ep}u+ \eta K u  + N (u)=0.
\end{equation}

To formulate the time periodic problem in a functional analytic setting, we introduce the operator $B^{\ep}$ on $\mathcal{X}_a$ defined by
\[
B^{\ep} = \frac{a^\ep}{a} \del_t + L^{\ep}_{\mathcal{R}_{1,c}^\ep}
\]
with domain $D(B^{\ep}) = \mathcal{Y}_a$. 
The time periodic problem \eqref{4.5} is then formulated as 
\begin{equation}\label{4.6}
B^{\ep}u + \frac{a^{\ep}}{a} \omega \delt u+ \eta K u + N (u) = 0, \ \ \ u\in \mathcal{Y}_a.
\end{equation}

Let $z_{\pm}^{\ep}=e^{\pm iat}u_\pm^{\ep}$. It is easily verified that $B^{\ep} z_\pm^\ep = 0$. 
Let 
\[
[u]_{\pm, \ep} = \langle u, z_{\pm}^{\ep*} \rangle_\ep 
\]
with $z_{\pm}^{\ep*}=e^{\pm iat}u_\pm^{\ep*}$ and define the operators $\mathcal{P}_{\pm}^{\ep}$ by    
\[
\mathcal{P}_{\pm}^{\ep} u = [u]_{\pm, \ep} z_{\pm}^{\ep}. 
\]
The operators $\mathcal{P}_{\pm}^{\ep}$ are projections onto the spaces $\mathrm{span}\,\{z_\pm^\ep\}$, respectively. 
As was shown in \cite{Hsia-Kagei-Nishida-Teramoto1} (see Theorem \ref{thm5.6} below), $0$ is a semisimple eigenvalue of $B^{\ep}$ and the eigenprojection for the eigenvalue $0$ is given by 
\[
\mathcal{P}^\ep = \mathcal{P}_+^\ep + \mathcal{P}_-^\ep. 
\]

We set $\mathcal{Q}^\ep = I - \mathcal{P}^\ep$. Then 
\[
\mbox{ \rm $u \in R(\mathcal{Q}^\ep)$ if and only if $[u]_{+,\ep} = [u]_{-,\ep} = 0$.} 
\] 
If $u$ is real valued, then $[u]_{+,\ep} = \ov{ [u]_{-,\ep} }$, from which we see that $u \in R(\mathcal{Q}^\ep)$ if and only if $[u]_{+,\ep} = 0$, when $u$ is a real valued function. 

\vspace{2.5ex} 
We state the result on the Hopf bifurcation for the artificial compressible system \eqref{1.1}--\eqref{1.4}. 

\begin{thm}\label{thm4.2}
{\rm (i)} There exist positive numbers $\ep_2$ and $\delta_2$ such that if $0< \ep \le \ep_2$, then there exists a nontrivial solution $u_\delta^\ep\in \mathcal{Y}_a$ of \eqref{4.6} for $\eta = \eta^\ep_{\delta}$ and $\omega = \omega^\ep_{\delta}$, where 
\begin{align*}
\eta^\ep_{\delta} & = \tilde\eta_0^\ep \delta^2 + \tilde\eta_1^\ep (\delta) \delta^3, \\
\omega^\ep_{\delta} & =\tilde \omega_0^\ep \delta^2 + \tilde\omega_1^\ep (\delta) \delta^3, \\
u_\delta^\ep & =\delta(z_0^{\ep} + \delta U_\delta^\ep)  
\end{align*}
for $|\delta|\le \delta_2$. Here $z_0^{\ep}=\re(e^{iat} u_+^\ep)$ and $[U_\delta^\ep]_{\ep, +}=0$; and $\tilde{\eta}^{\ep}_{0}$, $\tilde{\omega}^{\ep}_{0}$, $\tilde{\eta}^{\ep}_{1}(\delta)$, $\tilde{\omega}^{\ep}_{1}(\delta)$ and $U_\delta^\ep$ satisfy $|\tilde{\eta}^{\ep}_{0}| + |\tilde{\omega}^{\ep}_{0}| \le C$ and $|\tilde{\eta}^{\ep}_{1}(\delta)| + |\tilde{\omega}^{\ep}_{1}(\delta)| \le C$; and, furthermore, $U^{\ep}_{+}$ satisfies $\delt U^{\ep}_{\delta} \in \mathcal{Y}_{a}$ with $[\delt U^{\ep}_{\delta} ]_{+, \ep} = 0$ and $||| U_\delta^\ep|||_{\ep, \mathcal{Y}_{a} } + ||| \delt U_\delta^\ep|||_{\ep, \mathcal{Y}_{a} } \le C$ uniformly in $0 < \ep \le \ep_{2}$ and $|\delta| \le \delta_{2}$. 

\vspace{1ex}
{\rm (ii)} Let $u_\delta=\trans(\phi_\delta,\bu_\delta)=\trans(\delta(p_0 + \delta\Phi_\delta), \delta(\b z_0 + \delta\b U_\delta))$ be the time periodic bifurcating solution of \eqref{3.8} for $\eta= \eta_{\delta}$ and $\omega=\omega_{\delta}$ given in Proposition {\rm \ref{prop3.4}}. 
Then 
\begin{align*}
& \tilde\eta_0^\ep = \tilde\eta_0 + O(\ep),  \ \ \tilde\eta_1^\ep (\delta) =  \tilde\eta_1(\delta) + O(\ep), \\  
& \tilde\omega_0^\ep = \tilde\omega_0 + O(\ep), \ \ \tilde\omega_1^\ep (\delta) =  \tilde\omega_1 (\delta) + O(\ep)
\end{align*}
and  
\[
|||U_\delta^\ep - U_\delta |||_{\ep, \mathcal{Y}_a} \le C \ep 
\]
uniformly for $0< \ep \le \ep_2$ and $|\delta|\le \delta_2$, where $\tilde\eta_0$, $ \tilde\omega_0$, $\tilde\eta_1(\delta)$ and $\tilde\omega_1(\delta)$ are the numbers given in Proposition {\rm \ref{prop3.4}}.  
In particular, 
\begin{align*}
& |\eta^\ep_{\delta} - \eta_{\delta}| + |\omega^\ep_{\delta} - \omega_{\delta}| \le C \ep \delta^2, \\
& \|\phi_\delta^\ep-\phi_\delta\|_{C(\mathbb{T}_{\frac{2\pi}{a}};L^2)} + \|\phi_\delta^\ep-\phi_\delta\|_{L^2(\mathbb{T}_{\frac{2\pi}{a}};H^1)} \le C |\delta|
\\
& \|\b {u}_\delta^\ep-\b{u}_\delta\|_{C(\mathbb{T}_{\frac{2\pi}{a}};\b X)} + \|\b {u}_\delta^\ep-\b{u}_\delta\|_{L^2(\mathbb{T}_{\frac{2\pi}{a}};\b X^1)} \le C\ep |\delta|
\\
& \|\del_x^2(\b{u}_\delta^\ep-\b{u}_\delta)\|_{L^2(\mathbb{T}_{\frac{2\pi}{a}};\b X)} + \|\delt(\b{u}_\delta^\ep-\b{u}_\delta)\|_{L^2(\mathbb{T}_{\frac{2\pi}{a}};\b X)} \le C |\delta|
\end{align*}
uniformly for $0< \ep \le \ep_2$ and $|\delta|\le \delta_2$. 

\vspace{1ex}
{\rm (iii)} There exists a neighborhood $\mathcal{O}$ of $\trans(\eta,\omega,u) = \trans(0,0,0)$ in $\mathbb{R}\times\mathbb{R}\times\mathcal{Y}_a$ such that if $\trans(\eta,\omega,u)\in\mathcal{O}$ is a solution of \eqref{4.6}, then
\[
\ts
\trans(\eta,\omega,u) \in \{\trans(\eta,\omega,0)\in\mathcal{O}\}\cup\{\trans(\eta^{\ep}_{\delta} ,\omega^{\ep}_{\delta}, u^{\ep}_{\delta}(\cdot + \tau ) ); |\delta| \le \delta_1, -\frac{\pi}{a} \le \tau < \frac{\pi}{a} \}
\]
where $\{\eta^{\ep}_{\delta},\omega^{\ep}_{\delta},u^{\ep}_{\delta} \}$ is the solution branch obtained in {\rm (i)}.
\end{thm}

The stability of the bifurcating solution $u^{\ep}_{\delta}$ is also an important question; we shall investigate it in section \ref{Stability}. 

\vspace{1ex}
To prove Theorem \ref{thm4.2}, we employ the spectral properties of $ - B^{\ep}$ which was investigated in \cite{Hsia-Kagei-Nishida-Teramoto1}. Some of them will be summarized in section \ref{Spectrum of B} which will then be used in sections \ref{Proof of thm4.2} and \ref{Stability}. 

In what follows we assume that $\mathcal{R}_{1}$ is in the interval $[\frac{1}{2}\mathcal{R}_{1,c}, \frac{3}{2}\mathcal{R}_{1,c}]$. 

\section{Spectrum of $B^{\ep}$}\label{Spectrum of B}
 
In this section we summarize the results on the spectrum of $B^{\ep}$ obtained in \cite{Hsia-Kagei-Nishida-Teramoto1}. 

We begin with the following lemma.

\begin{lem}\label{lem5.1} {\rm (\cite[Lemma 5.1]{Hsia-Kagei-Nishida-Teramoto1})} 
Let $\beta$ and $T$ be positive constants. Suppose that $u_0\in X^1$ and let $F=\trans(f,\b g,h,k)\in  L^2(0,T; X)$. 
Then there is a unique solution $u(t)\in L^2(0, T;Y) \cap H^1(0, T;X)$ of the problem
\begin{align}\label{5.1}
\beta \del_t u + L^{\ep}_{\mathcal{R}_{1,c}^\ep}u & =  F, \\ \label{5.2}
u(0)  & = u_0
\end{align}
with estimate
\begin{align*}
& \beta ||| u(t) |||_{\ep, X^1}^2 
+ \int_0^t \left(|||\delx u |||_\ep^2 + \beta^{2} \ep ^2 |||\delt u |||_\ep^2 + \ep^2 \|\delx ^2\bu \|_2^2 + \beta^{2} \ep^6 \|\delx\delt \phi\|_2^2 \right) \,ds \\
& \ \le C e^{\frac{C}{\beta} t} \left\{ ||| u_0 |||_{\ep, X^1}^2  
+ \int_0^t \left( \ep^2\| f \|_2^2 + |\re(\b F, \bu)| +\ep^2 \|\b F\|_2^2 + \ep^6 \|\delx f\|_2^2 \right) \, ds \right\}
\end{align*}
uniformly for $0 \le t \le T$ and $0 <\ep \le \ep_{1}$ and $\beta$. 
Furthermore, if $\beta \ge \beta_{1}$, then the estimate 
\[
||| u |||_{\ep, \mathcal{Y}(0, T) } \le C \left\{ ||| u_0 |||_{\ep, X^1} + ||| F |||_{\ep, \mathcal{X}(0, T) } \right \}
\]
holds uniformly for $0 <\ep \le \ep_{1}$ with a constant depending only on $\beta_{1}$. 
\end{lem}

We denote by $\mathscr{V}^{\ep}_{\omega} (t)$ the restriction of the semigroup $ e^{-\frac{a t}{ a^{\ep} (1 + \omega)} L^{\ep}_{ \mathcal{R}_{1,c}^\ep } }$ to the space $X^{1}$: 
\[
\mathscr{V}^{\ep}_{\omega} (t) =  e^{-\frac{a t}{ a^{\ep} (1 + \omega)} L^{\ep}_{ \mathcal{R}_{1,c}^\ep } } \,|_{X^{1} }. 
\] 
As for the decay property of $\mathscr{V}^{\ep}_{\omega} (t)$, we have the following result.  

We define the projections $P^{\ep}$ and $Q^\ep$ by  
\[
\mbox{ \rm $P^{\ep} = P^{\ep}_{+} + P^{\ep}_{-}$ and $Q^\ep = I - P^\ep$, } .
\] 
respectively. 
Since $\pm ia^{\ep}$ are simple eigenvalues of $ - L^{\ep}_{\mathcal{R}^{\ep}_{1,c}}$, we have   
\[
P^{\ep}_{\pm} \mathscr{V}^{\ep}_{\omega}(t) = \mathscr{V}^{\ep}_{\omega}(t) P^{\ep}_{\pm}  = e^{ \pm \frac{ iat }{ 1 + \omega} } P^{\ep}_{\pm}. 
\]

The $Q^{\ep}$ part of $\mathscr{V}^{\ep}_{\omega}(t)$ has the following decay property. 

\begin{lem}\label{lem5.5} {\rm (\cite[Lemma 4.3]{Hsia-Kagei-Nishida-Teramoto1})} 
There exist positive constants $\ep_{1}$ and $\kappa_1$ such that the estimate 
\[
|||\mathscr{V}^{\ep}_{0}(t) Q^\ep u_0 |||_{\ep, X^1} \le Ce^{-\kappa_1 t} |||u_0|||_{\ep, X^1}
\]
holds uniformly in $t\ge0$ and $0<\ep\le\ep_{1}$. 
Furthermore, if $\re\lambda > -\kappa_1$, then $(I - e^{-\frac{2\pi}{a}\lambda} \mathscr{V}^{\ep}_{0} (\frac{2\pi}{a}) ) Q^\ep$ has a bounded inverse on $Q^\ep X^1$ and its inverse $[(I - e^{-\frac{2\pi}{a}\lambda} \mathscr{V}^{\ep}_{0} (\frac{2\pi}{a}) ) Q^\ep]^{-1}$ satisfies
\[
\| [ (I - e^{-\frac{2\pi}{a}\lambda} \mathscr{V}^{\ep}_{0} (\ts \frac{2\pi}{a} ) ) Q^\ep ]^{-1} \|
\le \frac{C}{1-e^{-\frac{2\pi}{a}(\re\lambda + \kappa_1)}}
\]
uniformly in $\ep\in(0,\ep_{1}]$.
\end{lem}

We also have the following spectral property of the monodromy operator $\mathscr{V}^{\ep}_{\omega}({\ts \frac{2\pi}{a} })$. 

\begin{lem}\label{rem5.9} {\rm (\cite[Remark 5.6]{Hsia-Kagei-Nishida-Teramoto1})} 
For any $r > 0$, there exists a positive constant $\omega_{r} = O(r)$ as $r \to 0$ such that if $|\omega| \le \omega_{r}$ then $\mu -  \mathscr{V}^{\ep}_{\omega}({\ts \frac{2\pi}{a} } ) $ has a bounded inverse on $X^1$ for $\mu \in \mathbb{C}$ satisfying $|\mu - 1| \ge r$ and $|\mu| \ge e^{-\frac{3}{4}\kappa_{1}\frac{2\pi}{a} }$ and $(\mu -  \mathscr{V}^{\ep}_{\omega}({\ts \frac{2\pi}{a} }) )^{-1} $ satisfies the estimate 
\[
||| (\mu -  \mathscr{V}^{\ep}_{\omega}({\ts \frac{2\pi}{a} } ) )^{-1} |||_{\ep, X^{1} } \le C\left(\frac{1}{r} + \frac{1}{|\mu|} \right) ||| F |||_{\ep, X^{1}} 
\]
uniformly in $\ep\in(0,\ep_{1}]$. 
\end{lem}

As for the operator $B^{\ep}$, we have the following result on its spectrum (\cite{Hsia-Kagei-Nishida-Teramoto1}). 

\begin{thm}\label{thm5.6} {\rm (\cite[Theorem 4.2]{Hsia-Kagei-Nishida-Teramoto1})} 
{\rm (i)} If $\ep\in(0,\ep_{1}]$, then
\[
\rho(-B^{\ep})\supset \Sigma_1 \setminus (\cup_{k \in \mathbb{Z}} \{ika^{\ep} \}).
\]
Here $\Sigma_1=\{\lambda\in\mathbb{C};\re\lambda>-\frac{a^{\ep} }{a} \kappa_1 \}$ with $\kappa_1$ being the positive number given in Lemma {\rm \ref{lem5.5}}; and each $ika^{\ep} $ is a semisimple eigenvalue of $B^{\ep}$ and the corresponding eigenspace is spanned by $ e^{i (1 - k) at}u_{+}^{(\ep)}$ and $ e^{- i (k + 1) at}u_{-}^{(\ep)}$.

\vspace{1ex} 
{\rm (ii)} If $\lambda\in\Sigma_1 \setminus (\cup_{k\in\mathbb{Z}}\{ika^{\ep} \})$, then
\begin{align*}
& (\lambda + B^{\ep})^{-1}F \\
& \quad =  \frac{ 2\pi e^{-\frac{a}{a^{\ep} } \lambda t} }{a^{\ep} ( 1 - e^{-\frac{2\pi}{a^{\ep} }\lambda}) } \left( \left[ e^{-\frac{a}{ a^{\ep} } \lambda (\frac{2\pi}{a}-s) } F \right]_{+, \ep} z_+^\ep + \left[ e^{-\frac{a}{ a^{\ep} } \lambda (\frac{2\pi}{a}-s) } F \right]_{-, \ep} z_-^\ep \right)  
\\
& \quad \quad + \frac{a}{ a^{\ep} } e^{-\frac{a}{ a^{\ep} }\lambda t} \mathscr{V}^{\ep}_{0}(t) \left[ (I- e^{-\frac{2\pi}{ a^{\ep} } \lambda } \mathscr{V}^{\ep}_{0}({\ts\frac{2\pi}{a} }) ) Q^\ep  \right]^{-1} 
\\
& \quad \quad \quad \quad \cdot Q^\ep
 \int_0^{\frac{2\pi}{a}} e^{-\frac{a}{ a^{\ep} } \lambda (\frac{2\pi}{a} - s ) } \mathscr{V}^{\ep}_{0}( {\ts\frac{2\pi}{a} } - s ) F(s) \,ds 
\\
&\quad \quad + \frac{a}{ a^{\ep} } \int_0^t  e^{-\frac{a}{a^{\ep}} \lambda (t-s)} \mathscr{V}^{\ep}_{0}( t - s ) F(s)\,ds.
\end{align*}
Furthermore, 
\[ ||| (\lambda + B^{\ep})^{-1}F |||_{\mathcal{Y}_a} \le \frac{ C}{| 1 - e^{-\frac{2\pi}{ a^{\ep} }\lambda } | } \sum_{j = +, -} \left|\left[ e^{-\frac{a}{a^{\ep} } \lambda (\frac{2\pi}{a}-s) } F \right]_{j,\ep} \right| + C|||F |||_{\mathcal{X}_a}. 
\]

\vspace{1ex}
{\rm (iii)} $0$ is a semisimple eigenvalue of $B^{\ep}$ and 
\[
\mathcal{X}_a = N(B^{\ep}) \oplus R(B^{\ep}), \ \ \ N(B^{\ep}) = \mathrm{span}\,\{z_+^\ep, z_-^\ep\}.
\]
The projection $\mathcal{P}^\ep$ is an eigenprojection for the eigenvalue $0$ and $\mathcal{Q}^\ep$ is a projection on $R(B^{\ep})$ along $N(B^{\ep})$. 

\vspace{1ex}
{\rm (iv)} Let $F = \trans(f,\b F)\in \mathcal{X}_a$. Then the equation 
\[ 
B^{\ep}u = F, \ \  u \in \mathcal{Y}_a
\]
is solvable if and only if $F \in \mathcal{Q}^\ep \mathcal{X}_a$, i.e., $F$ satisfies $[F]_{+,\ep} = [F]_{-,\ep} = 0$. If this condition for $F$ is satisfied, then the problem $B^{\ep}u = F$ with $u \in \mathcal{Q}^\ep \mathcal{Y}_a$ is uniquely solvable and the solution $u \in \mathcal{Q}^\ep \mathcal{Y}_a$ is given by 
\begin{align*}
u(t)  & = \frac{a}{ a^{\ep} } \mathcal{P}^\ep (sF(s))  
\\ 
& \quad + \frac{a}{ a^{\ep} } \mathscr{V}^{\ep}_{0}(t) \left[ (I- \mathscr{V}^{\ep}_{0}({\ts \frac{a}{ a^{\ep} } } ) ) Q^\ep \right] ^{-1}\int_0^{\frac{2\pi}{a}} \mathscr{V}^{\ep}_{0}({\ts \frac{2\pi}{a} }-s ) F(s)\,ds 
\\
& \quad  +\frac{a}{ a^{\ep} } \int_0^t \mathscr{V}^{\ep}_{0}( t - s ) F(s)\,ds. 
\end{align*}
Furthermore, $u$ satisfies
\[
||| u |||_{\ep, \mathcal{Y}_a} \le C \left\{ \int_0^\frac{2\pi}{a} \left( \ep^2\| f \|_2^2 + |\re(\b F, \bu)| +\ep^2 \|\b F\|_2^2 + \ep^6 \|\delx f\|_2^2 \right) \, dt \right\}^{\frac{1}{2}}
\]
and 
\[
||| u |||_{\ep, \mathcal{Y}_a} \le C ||| F |||_{\ep, \mathcal{X}_a} 
\]
uniformly in $\ep\in(0,\ep_{1}]$.
\end{thm}

\vspace{1ex}  
To prove the Hopf bifurcation in \eqref{4.6}, Theorem \ref{thm5.6} is insufficient; the term $\omega \delt u$ on the left-hand side of the $\mathcal{Q}^\ep$ part of the equation \eqref{4.6}  cannot be regarded as a perturbation of $B^{\ep} u$. 
To overcome this we shall put the term $\omega \delt$ into the principal part, so we define the operator $B^{\ep}(\omega)$ on $\mathcal{X}_a$ by 
\[
D(B^{\ep}( \omega))=\mathcal{Y}_a, \ \ B^{\ep}( \omega) u = \frac{a^{\ep}}{a} (1+\omega)\delt u + L^{\ep}_{\mathcal{R}_{1,c}^\ep} u \ \ (u \in D(B^{\ep}( \omega))). 
\]
As for $B^{\ep}(\omega)$, we have the following estimate. 
  
\begin{lem}\label{lem5.8}{\rm (\cite[Lemma 5.6]{Hsia-Kagei-Nishida-Teramoto1})}  
There exist a positive constant $\ep_{1}$ such that the following assertion holds.  
For a given $F \in \mathcal{Q}^\ep \mathcal{X}_a$, there exists a unique solution $u \in \mathcal{Q}^\ep \mathcal{Y}_a$ of 
$B^{\ep}( \omega) u = F $, and the solution $u$ satisfies the estimate 
\[
||| u |||_{\ep, \mathcal{Y}_a} \le C \left\{ \int_0^\frac{2\pi}{a} \left( \ep^2\| f \|_2^2 + |\re(\b F, \bu)| +\ep^2 \|\b F\|_2^2 + \ep^6 \|\delx f\|_2^2 \right) \, dt \right\}^{\frac{1}{2}}
\]
and  
\[
||| u |||_{\ep, \mathcal{Y}_a} \le C||| F |||_{\ep, \mathcal{X}_a}
\]
uniformly for $0 <\ep \le \ep_{1}$ and $|\omega| \le \frac{1}{4}$.  
\end{lem}

%Similarly to above, one can also obtain the following solution formula for $( \lambda + B^{\ep}(\omega) ) u = F$ for $\lambda$ satisfying $\re\lambda \ge - \frac{a^{\ep} }{a} \kappa_{1}$ and $\lambda \neq i a^{\ep} ((1 + \omega ) k \pm 1)$ $(k \in \mathbb{Z})$: 
%\begin{equation}\label{5.39}
%u (t)  =  u_1(t) + u_2(t) + u_3(t), 
%\end{equation} 
%where 
%\begin{align*}
%u_1(t) & =  \frac{2\pi}{a^{\ep} ( 1 + \omega )} \sum_{ j= +, -} \frac{1}{ e^{ \frac{2\pi (\lambda - j i a^{\ep} ) }{ a^{\ep} ( 1 + \omega )} } - 1 } \left[ e^{ \frac{ a( \lambda + j i a^{\ep} \omega ) }{ a^{\ep} (1 + \omega) } s } F \right]_{j, \ep} e^{ - \frac{ a ( \lambda  - j i a^{\ep} \omega ) }{ a^{\ep} (1 + \omega) } t } u_j^\ep, \\
%u_2(t) & = \frac{a}{a^{\ep} (1 + \omega)} \mathscr{V}^{\ep}_{\omega} (t) \left[ \left( I -  e^{ -\frac{2 \pi \lambda }{ a^{\ep} ( 1 + \omega ) } } \mathscr{V}^{\ep}_{\omega}({\ts \frac{2\pi}{a} } ) \right) Q^\ep \right]^{-1} \\[2ex]
%& \quad \quad \cdot Q^\ep \int_0^{\frac{2\pi}{a}} e^{ -\frac{ a \lambda }{ a^{\ep} ( 1 + \omega ) } \left( \frac{ 2 \pi }{a} - s \right) } \mathscr{V}^{\ep}_{\omega} ({\ts \frac{2\pi}{a} -s } ) F(s)\,ds, \\
%u_3(t) & = \frac{a}{a^{\ep}(1+\omega)} \int_0^t  e^{ -\frac{ a \lambda }{ a^{\ep} ( 1 + \omega ) } ( t - s ) } \mathscr{V}^{\ep}_{\omega}( t - s ) F(s)\,ds.
%\end{align*}
%This representation will be employed to study the stability of the bifurcating time periodic solutions. 
%\end{rem} 

\section{Proof of Theorem \ref{thm4.2}}\label{Proof of thm4.2}

In this section we prove Theorem \ref{thm4.2}, i.e., we shall prove that a nontrivial solution branch of \eqref{4.6} bifurcates from $\{\eta, u\}=\{0,0\}$. 
In what follows, we will use the same letters $\ep_{0}$ and $\delta_{0}$ for bounds of the ranges of $\ep$ and $\delta$, even when they should be taken suitably smaller than those of the previous ones if no confusion will occur from the context. 

\subsection{Proof of Theorem \ref{thm4.2} (i) } 

We look for a solution $u$ of \eqref{4.6} in the form
\begin{equation}\label{6.1}
u = \delta(z_0^\ep + \delta U^\ep),
\end{equation}
where $\delta$ is a small parameter; $z_0^\ep = 2\re z_+^\ep$; and $U^\ep$ is a real valued function in $\mathcal{Y}_a$ satisfying $[U^\ep ]_{+,\ep}$ = 0. 

We note that if $U$ is real valued, then $[U]_{-,\ep} = \overline{[U]_{+,\ep}}$. Therefore, the condition $[U]_{+,\ep} = 0$ automatically implies the condition $[U]_{-,\ep} = 0$, which in turn implies that $U = \mathcal{Q}^\ep U \in R(B^{\ep}) = \mathcal{Q}^\ep \mathcal{X}_a$ by Theorem \ref{thm5.6}.

We first collect basic facts concerning $\mathcal{P}^{\ep}$. 

\begin{lem}\label{lem6.1} 
{\rm (i)} It holds that $[z_0^\ep]_{+,\ep} = 1$, $[\delt z_0^\ep]_{+,\ep} = ia$ and $\mathcal{P}^\ep z_0^\ep = z_0^\ep$. 

\vspace{1ex}
{\rm (ii)} It holds that $\mathcal{P}^{\ep}_{\pm} \delt \subset \delt \mathcal{P}^{\ep}_{\pm} = \pm ia \mathcal{P}^{\ep}_{\pm}$ and $\mathcal{P}^{\ep} \delt \subset \delt\mathcal{P}^{\ep}$. 

\vspace{1ex}
{\rm (iii)} If $u$ is real valued, then $\mathcal{P}^\ep u = 0$ is equivalent to $[u]_{+, \ep} = 0$ that is equivalent to $\re[u]_{+,\ep} = \im[u]_{+,\ep} = 0$. 
\end{lem}

Lemma \ref{lem6.1} can be proved by straightforward computations. We omit the proof. 

The following estimates will be used to estimate nonlinear terms. 
We define the bilinear form $M(\,\cdot \, , \, \cdot \, )$ by 
\[
M(u_{1}, u_{2}) = 
\begin{pmatrix}
0 \\ 
\b{M}(\bu_{1}, \bu_{2}) 
\end{pmatrix}, 
\quad \b{M}(\bu_{1}, \bu_{2}) = \b{N}(\bu_{1}, \bu_{2}) + \b{N}(\bu_{2}, \bu_{1}) 
\]
for $u_{j} = \trans(\phi_{j}, \bu_{j})$ $(j= 1,2)$. 

\begin{lem}\label{lem6.2}
The following inequalities hold for $u_j=\trans(\phi_j, \bu_j) \in \mathcal{Y}_a$ $(j = 1,2,3)$ uniformly in $0 < \ep \le \ep_{0}$:

\vspace{1ex} 
{\rm (i)} \quad $|[u]_{+,\ep}|\le C\opnorm{u}_{\ep, \mathcal{X}_a}$, 

\vspace{1ex} 
{\rm (ii)} \quad $|[N(u_{1}, u_{2})]_{+,\ep}|\le C\opnorm{u_{1}}_{\ep, \mathcal{Y}_a} \opnorm{u_{2}}_{\ep, \mathcal{Y}_a}$, 

\vspace{1ex} 
{\rm (iii)} \quad $||| M(z^{\ep}_{0}, u) |||_{\ep, \mathcal{X}_{a} } \le C ||| u |||_{\ep, \mathcal{Y}_{a} }$, 

\vspace{1ex} 
{\rm (iv)} \quad $\ds \int_0^{\frac{2\pi}{a}} \left|(\b N(\bu_1,\bu_2),\bu_3)\right| \, dt \le C ||| u_1 |||_{\ep, \mathcal{Y}_a} ||| u_2 |||_{\ep, \mathcal{Y}_a} ||| u_3 |||_{\ep, \mathcal{Y}_a}$,  

\vspace{1ex}
{\rm (v)} \quad $\ds \int_0^{\frac{2\pi}{a}} \ep^2 \|\b N(\bu_1,\bu_2)\|_2^2 \, dt \le C ||| u_1 |||_{\ep, \mathcal{Y}_a}^2 ||| u_2 |||_{\ep, \mathcal{Y}_a}^2$. 
\end{lem}

We look for a solution $u$ in the form \eqref{6.1} of the equation \eqref{4.6} for $\eta = \delta^{2} \tilde{\eta}^\ep$ and $\omega = \delta^{2}\tilde{\omega}^\ep$. This scaling of $\eta$ and $\omega$ is due to the fact $[N(z^{\ep}_{0}) ]_{+, \ep} = 0$. We substitute $\{\eta,\omega,u\} = \{\delta^{2}\tilde{\eta}^\ep, \delta^{2} \tilde{\omega}^\ep,\delta(z_0^\ep+ \delta U^\ep)\}$ into \eqref{4.6} and divide the resulting equation by $\delta^2$. We then obtain
\begin{equation}\label{6.2}
B^{\ep} U^\ep + \frac{a^{\ep}}{a} \delta \tilde{\omega}^\ep (\delt z_0^\ep + \delta \delt U^\ep ) + \delta \tilde{\eta}^\ep K (z_0^\ep + \delta U^\ep) + N (z_0^\ep + \delta U^\ep) = 0.
\end{equation} 

We decompose \eqref{6.2} into its $\mathcal{P}^\ep$ and $\mathcal{Q}^\ep$ parts. 
Applying $\mathcal{P}^\ep$ and $\mathcal{Q}^\ep$ to \eqref{6.2} and noting that $[N(z^{\ep}_{0}) ]_{+, \ep} = 0$, we see from Lemma \ref{lem6.1} that 
\begin{equation}\label{6.5}
\tilde{\eta}^\ep\re[K z_0^\ep]_{+,\ep} + \delta \tilde{\eta}^\ep\re[K U^\ep]_{+,\ep} + \re[M(z_0^\ep, U^{\ep}) + \delta N(U^\ep)]_{+,\ep} = 0,
\end{equation}
\begin{equation}\label{6.6}
a^\ep \tilde{\omega}^\ep +  \tilde{\eta}^\ep\im[K z_0^\ep]_{+,\ep} + \delta \tilde{\eta}^\ep \im[K U^\ep]_{+,\ep} + \im[M(z_0^\ep, U^{\ep}) + \delta N(U^\ep)]_{+,\ep} = 0,
\end{equation}
\begin{equation}\label{6.7}
B^{\ep}(\delta^{2} \tilde{\omega}^\ep) U^\ep + \mathcal{Q}^\ep \left(\delta \tilde{\eta}^\ep K (z_0^\ep + \delta U^\ep) +N (z_0^\ep + \delta U^\ep) \right) = 0.
\end{equation} 
Here, as explained in section \ref{Spectrum of B}, we regard the term $\frac{a^{\ep} }{a} \delta^{2} \tilde{\omega}^{\ep} \delt U^{\ep}$ in \eqref{6.7} as a part of the principal part of the equation \eqref{6.7}. 

We define the map $\mathcal{B}^{\ep}(\omega):\mathbb{R}\times \mathbb{R}\times\mathcal{Q}^\ep \mathcal{Y}_a\rightarrow\mathbb{R}\times\mathbb{R}\times\mathcal{Q}^\ep \mathcal{X}_a$ by
\[
\mathcal{B}^{\ep}(\omega) =
\begin{pmatrix}
\re[K z_0^\ep]_{+,\ep} & 0 & \re [M(z_0^\ep, \, \cdot \, )]_{+, \ep} \\
\im[K z_0^\ep]_{+,\ep} & a^\ep & \im [M(z_0^\ep, \, \cdot \, )]_{+, \ep} \\
0 & 0 & B^{\ep}(\omega)\mathcal{Q}^\ep
\end{pmatrix}.
\]
The problem \eqref{6.5}--\eqref{6.7} is then written as
\begin{equation}\label{6.8}
\mathcal{B}^{\ep} (\delta^{2} \tilde{\omega}^\ep) \mathcal{U}^\ep = - \mathcal{N}^{\ep}(\delta;\mathcal{U}^\ep),
\end{equation}
where $\mathcal{U}^\ep = \trans ( \tilde{\eta}^\ep, \tilde{\omega}^\ep, U^\ep) \in\mathbb{R}\times\mathbb{R}\times \mathcal{Q}^\ep \mathcal{Y}_a$ and 
\[
\mathcal{N}^{\ep}(\delta;\mathcal{U}^\ep) = 
\begin{pmatrix}
\delta \tilde{\eta}^\ep \re[K U^\ep]_{+,\ep} + \delta \re[N (U^\ep)]_{+,\ep} \\
\delta \tilde{\eta}^\ep \im[K U^\ep]_{+,\ep} + \delta \im[N (U^\ep)]_{+,\ep} \\
\delta \tilde{\eta}^\ep\mathcal{Q}^\ep K (z^{\ep}_{0} + \delta U^\ep) + \mathcal{Q}^\ep N (z_0^\ep + \delta U^\ep)
\end{pmatrix}
\]

Let us suppose that $|\tilde{\eta}^\ep| \le 1$ and $|\tilde{\omega}^\ep| \le 1$. By Theorem \ref{thm4.1} we see that $\re[K z_0^\ep]_{+,\ep} = \re(\b K\bu_+,\boldsymbol{u}_+^*) + O(\ep^2) < 0$ and $a^\ep = a + O(\ep^2) > 0$ if $\ep>0$ is sufficiently small. Furthermore, Lemma \ref{lem5.8} implies that $B^{\ep}(\delta^{2} \tilde{\omega}^\ep): \mathcal{Q}^\ep \mathcal{Y}_a \rightarrow\mathcal{Q}^\ep \mathcal{X}_a$ has a bounded inverse if $\ep>0$ and $\delta$ are sufficiently small. As a result, $\mathcal{B}^{\ep}(\delta^2 \hat \omega^\ep)$ has a bounded inverse if $0 < \ep \le \ep_{0}$ and $\delta^{2} \le \frac{1}{4}$ for some small $\ep_{0} >0$. Therefore, if $0 < \ep \le \ep_{0}$ and $|\delta| \le \delta_{0}$ for some small $\ep_{0}, \delta_{0}>0$, then \eqref{6.8} is rewritten as 
\begin{equation}\label{6.9}
\mathcal{U}^\ep = - \mathcal{B}^{\ep}(\delta^{2} \tilde{\omega}^\ep)^{-1}\mathcal{N}^{\ep}(\delta;\mathcal{U}^\ep). 
\end{equation}
We note that if $\tilde{\mathcal{U} } = \trans(\tilde{\eta}, \tilde{\omega}, U)$ with $U = \mathcal{Q}^{\ep} U$, then   
\[
\mathcal{B}^{\ep}(\omega)^{-1} \tilde{\mathcal{U}} 
= 
\begin{pmatrix} 
(\mathcal{B}^{\ep}(\omega)^{-1} \tilde{\mathcal{U}})^{1} \\
(\mathcal{B}^{\ep}(\omega)^{-1} \tilde{\mathcal{U}})^{2} \\
(\mathcal{B}^{\ep}(\omega)^{-1} \tilde{\mathcal{U}})^{3} \\
\end{pmatrix},  
\]
where 
\begin{align*}
(\mathcal{B}^{\ep}(\omega)^{-1} \tilde{\mathcal{U}})^{1}  & = \frac{1}{\re[ K z_0^\ep]_{+,\ep}}\left\{ \tilde {\eta} - \re[M(z^{\ep}_0, B^{\ep}(\omega)^{-1} U) ]_{+, \ep} \right\}, \\[1ex]
(\mathcal{B}^{\ep}(\omega)^{-1} \tilde{\mathcal{U}})^{2}  & = \frac{1}{a^{\ep} } \tilde{\omega}  - \frac{\im[ K z_0^\ep ]_{+,\ep}}{a^{\ep} } (\mathcal{B}^{\ep}(\omega)^{-1} \tilde{\mathcal{U}})^{1} \\[1ex]
& \quad - \frac{1}{a^{\ep} } \im [M(z^{\ep}_0, B^{\ep}(\omega)^{-1} U) ]_{+, \ep}, \\[1ex]
(\mathcal{B}^{\ep}(\omega)^{-1} \tilde{\mathcal{U}})^{3} & = \left( B^{\ep}(\omega)\mathcal{Q}^\ep \right)^{-1} U. 
\end{align*}

To solve \eqref{6.8}, we introduce a norm of $\mathbb{R} \times \mathbb{R} \times \mathcal{Y}_a$. 
We define the norm $||| \mathcal{U} |||_{\ep, \mathbb{R} \times \mathbb{R} \times \mathcal{Y}_a}$ of $\mathcal{U} = \trans( \hat \eta, \hat \omega, U) \in \mathbb{R} \times \mathbb{R} \times \mathcal{Y}_a$ by 
\[
||| \mathcal{U} |||_{\ep, \mathbb{R} \times \mathbb{R} \times \mathcal{Y}_a} = | \hat \eta| + | \hat \omega| + |||U |||_{\ep, \mathcal{Y}_a}.
\]
   
We construct approximate solutions $\{\mathcal{U}_n^\ep\}_{n=0}^\infty$ in the following way. Consider first the problem   
\begin{equation}\label{6.22}
B^{\ep} U_0^\ep + N(z_0^\ep)=0. 
\end{equation}
Since $[N(z^{\ep}_{0} ) ]_{+, \ep} = 0$, i.e., $N(z^{\ep}_{0}) = \mathcal{Q}^{\ep} N(z^{\ep}_{0})$, it follows from Theorems \ref{thm4.1} and \ref{thm5.6} and Lemma \ref{lem6.2} that \eqref{6.22} has a unique solution $U^{\ep}_{0} \in \mathcal{Q}^{\ep}\mathcal{Y}_{a}$ with estimate $|||U_0^\ep|||_{\ep,\mathcal{Y}_a} \le C |||N(z_0^\ep)|||_{\ep, \mathcal{X}_a} \le C$ uniformly for $0< \ep \le \ep_{0}$. We also note that $|||\delt U_0^\ep|||_{\ep,\mathcal{Y}_a} \le C |||M(\delt z_0^\ep, z_0^\ep)|||_{\ep, \mathcal{X}_a} \le C$ uniformly for $0< \ep \le \ep_{0}$. 

We now define $\trans(\tilde \eta_0^\ep, \tilde \omega_0^\ep)$ by the solution of 
\begin{equation}\label{6.23}
\tilde \eta_0^\ep \re[K z_0^\ep]_{+,\ep} + \re[M(z_0^\ep, U_0^\ep) ]_{+,\ep} = 0,
\end{equation}
\begin{equation}\label{6.24}
a^\ep  \tilde \omega_0^\ep +  \tilde \eta_0^\ep\im[K z_0^\ep]_{+,\ep} + \im[M(z_0^\ep, U_0^\ep) ]_{+,\ep} = 0.
\end{equation}
Using the estimate of $U^{\ep}_{0}$, we have $|\tilde \eta_0^\ep| + |\tilde \omega_0^\ep| \le C$ uniformly for $0< \ep \le \ep_{0}$. 
We thus obtain $\mathcal{U}^{\ep}_{0} = \trans(\tilde{\eta}_{0}, \tilde{\omega}_{0}, U^{\ep}_{0}) $ which is the solution of  
\[
\mathcal{B}^{\ep}(0)\mathcal{U}^\ep_{0} = - \mathcal{N}^{\ep}(0; 0) 
= 
-
\begin{pmatrix}
0 \\
0 \\
N (z_0^\ep)
\end{pmatrix}. 
\]
We set 
\begin{equation}\label{6.10} 
D_{1} = \sup \left\{\ts{ |||\mathcal{B}^{\ep}(\omega)^{ -1} \mathcal{N}^{\ep}(0; 0) |||_{\ep, \mathbb{R} \times \mathbb{R} \times \mathcal{Y}_{a} } : |\omega| \le \frac{1}{4}, \, 0< \ep \le \ep_{0}}  \right\}.  
\end{equation}
By Theorem \ref{thm4.1}, Lemmas \ref{lem5.8} and \ref{lem6.2}, we see that $D_{1} < \infty$ and $\mathcal{U}^{\ep}_{0}$ satisfies 
\[
||| \mathcal{U}^{\ep}_{0} |||_{\ep, \mathbb{R} \times \mathbb{R} \times \mathcal{Y}_{a} } \le D_{1} 
\]
uniformly in $0 < \ep \le \ep_{0}$.   

We now define $\mathcal{U}_n^\ep = \trans(\tilde{\eta}_n^\ep, \tilde{\omega}_n^\ep,U_n^\ep)$ $(n = 1, 2, \cdots )$ by 
\begin{align*}
\mathcal{U}_n^\ep & = - \mathcal{B}^{\ep}(\delta^{2} \tilde{\omega}^{\ep}_{n-1})^{-1} \mathcal{N}^{\ep}(\delta;\mathcal{U}_{n-1}^\ep)\,\,\,(n = 1,2,\cdots).
\end{align*}
We shall prove that the sequence $\{\mathcal{U}_n^\ep \}_{n=0}^\infty$ converges to a solution $\mathcal{U}^\ep = \trans(\tilde{\eta}^{\ep}, \tilde{\omega}^{\ep}, U^{\ep} )$ of \eqref{6.8} satisfying $\delt U^{\ep} \in \mathcal{Q}^{\ep}\mathcal{Y}_{a}$ with estimates for $\mathcal{U}^\ep$ and $\delt U^{\ep}$ uniformly in $0 < \ep \le \ep_{0}$ and $|\delta| \le \delta_{0}$. 

We claim that there exists a positive constant $\delta_{0}$ such that if $|\delta| \le \delta_{0}$, then
\begin{equation}\label{6.11}
||| \mathcal{U}_n^\ep |||_{\ep, \mathbb{R}\times\mathbb{R}\times \mathcal{Y}_a}\le 2D_{1}
\end{equation}
for all $n = 0, 1, 2, \cdots $. 

Indeed, \eqref{6.11} can be proved by induction. We have already seen that \eqref{6.11} holds for $n = 0$. Suppose that \eqref{6.11} holds with $n$ replaced by $n-1$ $(n \ge 1)$. We shall show that \eqref{6.11} also holds for $n$. 

We take $\delta$ so that $|\delta| \le 1$ and $D_{1} \delta^{2} \le \frac{1}{4}$. Writing 
\[
\mathcal{U}_n^\ep = -\mathcal{B}^{\ep}(\delta^2 \tilde{\omega}_{n-1}^\ep)^{-1} \left\{ \mathcal{N}^{\ep}(0;0) + (\mathcal{N}^{\ep}(\delta;\mathcal{U}_{n-1}^\ep) - \mathcal{N}^{\ep}(0;0))\right\},
\]
we deduce from Theorem \ref{thm5.6}, Lemmas \ref{lem5.8} and \ref{lem6.2} that 
\[
||| \mathcal{U}_n^\ep |||_{\ep, \mathbb{R} \times \mathbb{R} \times \mathcal{Y}_a} 
\le D_{1} + C |\delta| \left\{ D_{1} + (1+ |\delta|) D_{1}^{2} +  |\delta| D_{1} ||| U^{\ep}_{n} |||_{\ep, \mathcal{Y}_{a} }^{\frac{1}{2} } \right\}, 
\]
and hence,  
\[
||| \mathcal{U}_n^\ep |||_{\ep, \mathbb{R} \times \mathbb{R} \times \mathcal{Y}_a } \le \frac{4}{3} D_{1} + C |\delta| D_{1} \left\{ 1 + (1+ |\delta| + |\delta|^{3}) D_{1} \right\}.  
\]
Taking $\delta$ in such a way that $|\delta| \le \min\left\{ 1, \frac{1}{4D_{1}}, \frac{2}{3C( 1+ 3 D_{1} )}\right\}$, we have 
\[
||| \mathcal{U}_n^\ep |||_{\ep, \mathbb{R} \times \mathbb{R} \times \mathcal{Y}_a}  \le 2D_{1}. 
\]
This proves that \eqref{6.11} holds for $n$. We thus conclude by induction that \eqref{6.11} holds for all $n = 0, 1, 2, \cdots$. 

We next estimate $||| \delt U^{\ep}_{n} |||_{\ep, \mathcal{Y}_{a} }$. One can verify that $\delt U^{\ep}_{n}$ satisfies the linear problem 
\begin{align*}
B^{\ep}( \delta^{2} \tilde{\omega}^{\ep}_{n - 1}) U_n^\ep & = - M(z^{\ep}_{0}, \delt z^{\ep}_{0}) -  \delta \tilde{\eta}^{\ep}_{n -1} \mathcal{Q}^\ep  K (\delt z^{\ep}_{0} + \delta \delt U^{\ep}_{n - 1}) \\ 
& \quad - \delta \mathcal{Q}^\ep  (M (z_0^\ep, \delt U^{\ep}_{n - 1})  +  M(U^\ep_{n - 1}, \delt z^{\ep}_{0}) + \delta M(U^{\ep}_{n -1}, \delt U^{\ep}_{n - 1} ) ) 
\end{align*}    
and this linear problem has a unique solution in $\mathcal{Q}^{\ep} \mathcal{Y}_{a}$. Note that $\mathcal{Q}^\ep  M(z^{\ep}_{0}, \delt z^{\ep}_{0}) =  M(z^{\ep}_{0}, \delt z^{\ep}_{0})$ since $[ M(z^{\ep}_{0}, \delt z^{\ep}_{0})]_{+, \ep} = 0$.  

We set 
\begin{equation}\label{6.10'}
D_{2} = \max\{ D_{1}, \tilde{D}_{2} \}, 
\end{equation} 
where $\tilde{D}_{2} = \sup \left\{\ts{ |||B^{\ep}(\omega)^{ -1}  M(z^{\ep}_{0}, \delt z^{\ep}_{0}) |||_{\ep, \mathcal{Y}_{a} } : |\omega| \le \frac{1}{4}, \, 0< \ep \le \ep_{0}}  \right\}$.   

As above, one can prove
\begin{equation}\label{6.11'}
||| \delt U^{\ep}_{n} |||_{\ep, \mathcal{Y}_{a} } \le 2D_{2}
\end{equation} 
for all $n = 0, 1,2, \cdots$. Indeed, this is true for $n = 0$. Suppose that \eqref{6.11'} holds for $n$ replaced by $n - 1$ $(n \ge 1)$. We then see from \eqref{6.11}, Lemmas \ref{lem5.8}, \ref{lem6.1} and \ref{lem6.2} that 
\[
||| \delt U^{\ep}_{n} |||_{\ep, \mathcal{Y}_{a} } \le D_{2} + C|\delta| \left\{ D_{2} + |\delta| D_{2}^{2} +D_{2}^{\frac{1}{2} } \left( 1+ D_{2}^{\frac{1}{2} } \right) ||| \delt U^{\ep}_{n} |||_{\ep, \mathcal{Y}_{a} } \right\}, 
\]  
and hence, by the same argument as above, we see that \eqref{6.11'} also holds for $n$ if $|\delta| \le \delta_{0}$ for some small $\delta_{0} > 0$. We thus conclude that \eqref{6.11'} holds for all $n = 0, 1, 2, \cdots$. 

We next show that $\{\mathcal{U}^{\ep}_{n} \}_{n \ge 0}$ is a Cauchy sequence in $\mathbb{R} \times \mathbb{R} \times \mathcal{Q}^{\ep} \mathcal{Y}_{a}$. Since $\mathcal{U}^{\ep}_{n + 1} - \mathcal{U}^{\ep}_{n}$ satisfies 
\begin{align*}
\mathcal{B}^{\ep}(\delta^{2}\tilde{\omega}^{\ep}_{n}) (\mathcal{U}^{\ep}_{n + 1} - \mathcal{U}^{\ep}_{n}) & = -\frac{a^{\ep} }{a} \delta^{2} ( \tilde{\omega}^{\ep}_{n} - \tilde{\omega}^{\ep}_{n - 1})
\begin{pmatrix} 
0 \\ 
0 \\ 
\delt U^{\ep}_{n} 
\end{pmatrix} 
\\[1ex]
& \quad - (\mathcal{N}^{\ep} (\delta; \mathcal{U}^{\ep}_{n}) - \mathcal{N}^{\ep} (\delta; \mathcal{U}^{\ep}_{n - 1}) ), 
\end{align*}
we see from \eqref{6.11} and \eqref{6.11'} that 
\begin{align*} 
& ||| \mathcal{U}^{\ep}_{n + 1} - \mathcal{U}^{\ep}_{n} |||_{\ep, \mathbb{R} \times \mathbb{R} \times \mathcal{Y}_{a} } \\
& \quad \le C|\delta| (1 +|\delta| D_{1} +  |\delta| D_{2}) ||| \mathcal{U}^{\ep}_{n} - \mathcal{U}^{\ep}_{n - 1} |||_{\ep, \mathbb{R} \times \mathbb{R} \times \mathcal{Y}_{a} } \\ 
& \quad \quad  + C \delta^{2} D_{1}^{\frac{1}{2} } ||| \mathcal{U}^{\ep}_{n } - \mathcal{U}^{\ep}_{n - 1} |||_{\ep, \mathbb{R} \times \mathbb{R} \times \mathcal{Y}_{a} }^{\frac{1}{2} } ||| \mathcal{U}^{\ep}_{n + 1} - \mathcal{U}^{\ep}_{n} |||_{\ep, \mathbb{R} \times \mathbb{R} \times \mathcal{Y}_{a} }^{\frac{1}{2} }, 
\end{align*} 
from which we have 
\[
||| \mathcal{U}^{\ep}_{n + 1} - \mathcal{U}^{\ep}_{n} |||_{\ep, \mathbb{R} \times \mathbb{R} \times \mathcal{Y}_{a} } \le  C|\delta| ||| \mathcal{U}^{\ep}_{n } - \mathcal{U}^{\ep}_{n - 1} |||_{\ep, \mathbb{R} \times \mathbb{R} \times \mathcal{Y}_{a} }
\]
uniformly in $0 < \ep \le \ep_{0}$, $|\delta| \le \delta_{0}$ and $n = 1, 2, \cdots$. This implies that $\{\mathcal{U}^{\ep}_{n} \}_{n \ge 0}$ is a Cauchy sequence in $\mathbb{R} \times \mathbb{R} \times \mathcal{Q}^{\ep} \mathcal{Y}_{a}$ if $\delta_{0}$ is taken so that $0 < \delta_{0} \le 1/(2C)$. 

Consequently, we deduce that there exists $\mathcal{U}^\ep = \trans(\tilde{\eta}^{\ep}, \tilde{\omega}^{\ep}, U^{\ep}) \in \mathbb{R} \times \mathbb{R} \times \mathcal{Q}^{\ep} \mathcal{Y}_a$ with $\delt U^{\ep} \in \mathcal{Q}^{\ep} \mathcal{Y}_{a}$ such that  
\begin{align*} 
& \mbox{\rm  $\mathcal{U}_{n}^\ep$ $\to$ $ \mathcal{U}^\ep$ strongly in $\mathbb{R} \times \mathbb{R} \times \mathcal{Q}^{\ep} \mathcal{Y}_a$ }, 
\\
&  \mbox{\rm  $\delt U_{n}^\ep$ $\to$ $\delt {U}^\ep$ weakly in $\mathcal{Q}^{\ep} \mathcal{Y}_a$}. 
\end{align*}
Furthermore, this $\mathcal{U}^{\ep}$ is a solution of \eqref{6.8} and  $\mathcal{U}^{\ep}  = \trans(\tilde{\eta}^{\ep}, \tilde{\omega}^{\ep}, U^{\ep})$ satisfies the estimates 
\begin{equation}\label{6.11''}
||| \mathcal{U}^\ep |||_{\ep, \mathbb{R}\times\mathbb{R}\times \mathcal{Y}_a} \le 2D_{1}, \quad  ||| \delt U^{\ep} |||_{\ep, \mathcal{Y}_a} \le 2D_{2}
\end{equation}
uniformly for $0 < \ep \le \ep_{0}$ and $|\delta| \le \delta_{0}$.  
 
Similarly to above, one can show  
\begin{equation}\label{6.11'''}
||| \mathcal{U}^{\ep} - \mathcal{U}^{\ep}_{0} |||_{\ep, \mathbb{R} \times \mathbb{R} \times \mathcal{Y}_{a} } \le  C|\delta| 
\end{equation}
uniformly in $0 < \ep \le \ep_{0}$ and $|\delta| \le \delta_{0}$. This implies that $\tilde{\eta}^{\ep}$ and $\tilde{\omega}^{\ep}$ take the forms 
\begin{equation}\label{6.21'}
\tilde{\eta}^{\ep} = \tilde{\eta}^{\ep}_{0} + \tilde{\eta}^{\ep}_{1}(\delta) \delta, 
\quad  
\tilde{\omega}^{\ep} = \tilde{\omega}^{\ep}_{0} + \tilde{\omega}^{\ep}_{1}(\delta) \delta  
\end{equation}
for some $\tilde{\eta}^{\ep}_{1}(\delta)$ and $\tilde{\omega}^{\ep}_{1}(\delta)$ satisfying $|\tilde{\eta}^{\ep}_{1}(\delta)| + |\tilde{\omega}^{\ep}_{1}(\delta)| \le C$ uniformly for $0 < \ep \le \ep_{0}$ and $|\delta| \le \delta_{0}$. This completes the proof of Theorem \ref{thm4.2} (i). 
\hfill$\square$

\subsection{ Proof of Theorem \ref{thm4.2} (ii)}

We next prove Theorem \ref{thm4.2} (ii). 
In this subsection we denote the solution branch obtained in the previous subsection by $\{\eta^{\ep}_{\delta}, \omega^{\ep}_{\delta}, u^{\ep}_{\delta} \}$ with 
\[
\eta^{\ep}_{\delta} = \delta^{2} \tilde{\eta}^{\ep}_{\delta}, \quad \omega^{\ep}_{\delta} = \delta^{2} \tilde{\omega}^{\ep}_{\delta}, \quad u^{\ep}_{\delta} = \delta(z^{\ep}_{0} + \delta U^{\ep}_{\delta}), 
\]
where $\{\tilde{\eta}^{\ep}_{\delta}, \tilde{\omega}^{\ep}_{\delta}, U^{\ep}_{\delta} \}$ is the solution of \eqref{6.5}--\eqref{6.7} obtained in the previous subsection. Therefore, $ \tilde{\eta}^{\ep}_{\delta}$ and $\tilde{\omega}^{\ep}_{\delta} $ are  written in the forms in \eqref{6.21'}. 

We first consider the behavior of the time periodic solution branch for the incompressible problem \eqref{3.8} as $\delta \to 0$. 

Let $\{\eta_{\delta}, \omega_{\delta}. u_{\delta} \}$ be the time periodic solution branch for \eqref{3.8} given in Proposition \ref{prop3.4} which takes the form  
\[
\eta_{\delta} = \tilde{\eta}_{0} \delta^{2} + \tilde{\eta}_{1}(\delta)\delta^{3}, \quad \omega_{\delta} = \tilde{\omega}_{0} \delta^{2} + \tilde{\omega}_{1}(\delta)\delta^{3}, \quad u_{\delta} =\delta( z_{0} + \delta U_{\delta}). 
\]
Since $[\b N(\b z_0 )]_{+} = 0$, there exists a unique solution $\b{U}_{0} \in L^2(\mathbb{T}_{\frac{2\pi}{a}}; \mathbb{Y}_{\sigma} ) \cap H^1(\mathbb{T}_{\frac{2\pi}{a}};\mathbb{X}_\sigma)$ of 
\begin{equation}\label{6.25'}
\delt \b U_0 + \mathbb{L}_{\mathcal{R}_{1,c}} \b U_0 + \mathbb{P} \b N(\b z_0 ) = \b 0
\end{equation}
with $[\b{U}_{0 }]_{+} = 0$. Let $U_{0} = \trans(\Phi_{0}, \b{U}_{0}) \in \mathcal{Y}_{a}$ with $\Phi_{0}$ being the associated pressure of $\b{U}_{0}$. Then $U_0=\trans(\Phi_0, \b U_0)$ is the unique solution of   
\begin{equation}\label{6.25}
\delt 
\begin{pmatrix}
0 \\
\b U_0
\end{pmatrix}
+ L^{\ep}_{\mathcal{R}_{1,c}} U_0 + N(z_0)  = 0
\end{equation}
with $[\b{U}_{0 }]_{+} = 0$. 

Similarly to the proof of Theorem \ref{thm4.2} (i), one can show the following proposition. 

\begin{prop}\label{prop6.3'}
Let $\{\eta_{\delta}, \omega_{\delta}, u_{\delta}\}$ be the time periodic solution branch for \eqref{3.8} given in {\rm Proposition \ref{prop3.4}} and let $U_{0}= \trans(\Phi_{0}, \b{U}_{0}) \in \mathcal{Y}_{a}$ be the unique solution of \eqref{6.25} with $[\b{U}_{0 }]_{+} = 0$. Then  $\{\eta_{0}, \omega_{0}, U_{0}\}$ satisfies 
\begin{equation}\label{6.26} 
\tilde \eta_0 \re[\b K \b z_0]_{+} + \re[\b M (\b z_0, \b U_0) ]_{+} = 0,
\end{equation}
\begin{equation}\label{6.27}
a \tilde\omega_0 + \tilde \eta_0 \im[\b K \b z_0]_{+} + \im[\b M (\b z_0, \b U_0) ]_{+} 
= 0
\end{equation}
and  
\[
|\tilde{\eta}_{\delta} - \tilde{\eta}_{0}| + |\tilde{\omega}_{\delta} - \tilde{\omega}_{0}| + \| U_{\delta} - U_{0} \|_{\mathcal{Y}_{a} } \le C|\delta|
\]
uniformly for $|\delta| \le \delta_{0}$. 
\end{prop}

We now estimate the difference between the solution branches for the incompressible system and the artificial compressible system. 

Let $\{\eta_{\delta}, \omega_{\delta}, u_\delta \}$ be the time periodic solution branch of the incompressible problem \eqref{3.8} given in Proposition \ref{prop3.4}. Then $\{\eta_{\delta}, \omega_{\delta}, u_\delta \}$ satisfies 
\begin{equation}\label{6.13}
(1+\omega_{\delta}) \delt 
\begin{pmatrix}
0 \\
\bu_\delta
\end{pmatrix}
+L^{\ep}_{\mathcal{R}_{1,c}} u_\delta + \eta_{\delta} K u_\delta + N(u_\delta) =0. 
\end{equation}

We set $\eta_{\delta}=\delta^{2} \tilde \eta_{\delta}$ and $\omega_{\delta} = \delta^{2} \tilde \omega_{\delta}$. Then $\{\tilde \eta_{\delta}, \tilde \omega_{\delta}, U_{\delta} \}$ satisfies the equation    
\begin{equation}\label{6.14} 
\begin{split}
(1+\delta^{2} \tilde \omega_{\delta})\delt 
\begin{pmatrix}
0 \\
\b U_{\delta}
\end{pmatrix}
+ L^{\ep}_{\mathcal{R}_{1,c}} U_{\delta} 
+ \delta \tilde \omega_{\delta} \delt 
\begin{pmatrix}
0 \\
\b z_0
\end{pmatrix}
+ \delta \tilde \eta_{\delta} K(z_0 +\delta U_{\delta})  & \\  
+ N(z_0 + \delta U_{\delta}) &  = 0. 
\end{split} 
\end{equation}
 
As in the proof of Theorem \ref{thm4.2} (i), we see from \eqref{6.14} that  
\begin{equation}\label{6.15}
\tilde \eta_{\delta} \re[\b K \b z_0]_{+} + \delta \tilde \eta_{\delta} \re[\b K \b U_{\delta}]_{+} + \re[\b{M} (\b z_0, \b U_{\delta}) + \delta \b{N}(U_{\delta}) ]_{+} = 0,
\end{equation}
\begin{equation}\label{6.16}
a \tilde \omega_{\delta} + \tilde \eta_{\delta} \im[\b K \b z_0]_{+} + \delta \tilde \eta_{\delta} \im [\b K \b U_{\delta}]_{+} + \im[\b{M} (\b z_0, \b U_{\delta}) + \delta \b{N}(U_{\delta})]_{+} 
= 0.
\end{equation}
Applying $\mathcal{Q}^{\ep}$ to \eqref{6.14}, we have 
\begin{equation}\label{6.17}
\begin{split}
& B^{\ep}(\delta^{2} \tilde{\omega}^{\ep}_{\delta})\mathcal{Q}^{\ep} U_{\delta} + \mathcal{Q}^{\ep} H^{\ep} (\delta;\tilde{\omega}^{\ep}_{\delta},\tilde{\omega}_{\delta}, U_{\delta})  \\[1ex]  
& \quad\quad  + \mathcal{Q}^{\ep} \left( \delta \tilde{\eta}_{\delta} K(z_{0} + \delta U_{\delta}) + N(z_{0} + \delta U_{\delta}) \right) = 0.  
\end{split} 
\end{equation}
Here $H^{\ep} (\delta; \tilde{\omega}^{\ep}_{\delta}, \tilde{\omega}_{\delta}, U_{\delta})$ is a function of $\tilde{\omega}^{\ep}_{\delta}$, $\tilde{\omega}_{\delta}$ and $U_{\delta} = \trans(\Phi_{\delta}, \b{U}_{\delta})$ with $\hat{\mathcal{Q} }_{0} \b{U}_{\delta} = \b{0} $ given by  
\begin{align*}
H^{\ep} (\delta; \tilde{\omega}^{\ep}_{\delta}, \tilde{\omega}_{\delta}, U_{\delta})  & =  - \frac{a^\ep -a}{a} ( 1 + \delta^{2}\tilde{\omega}^{\ep}_{\delta} ) 
\begin{pmatrix} 
0 \\
\delt \b U_{\delta}
\end{pmatrix}
- \delta ^{2}( \tilde{\omega}^{\ep}_{\delta} - \tilde{\omega}_{\delta} ) 
\begin{pmatrix} 
0 \\
\delt \b U_{\delta}
\end{pmatrix}
\\[1ex] 
& \quad  - \frac{a^{\ep} }{a} (1 + \delta^{2} \tilde{\omega}^{\ep}_{\delta} ) 
\begin{pmatrix}
\delt \Phi_{\delta} \\ 
\b 0
\end{pmatrix} 
+\delta \tilde{\omega}_{\delta} 
\begin{pmatrix} 
0 \\
\delt \b{z}_{0}
\end{pmatrix}
 - (\mathcal{R}^\ep_{1,c} - \mathcal{R}_{1,c}) K U_{\delta}.
\end{align*}

We set $\mathcal{U}_{\delta} = \trans( \tilde \eta_{\delta}, \tilde \omega_{\delta}, U_{\delta})$. It then follows from \eqref{6.5}--\eqref{6.7} and \eqref{6.15}--\eqref{6.17} that 
\[
\mathcal{B}^{\ep}(\delta^{2} \tilde{\omega}^{\ep}_{\delta})(\mathcal{U}^{\ep}_{\delta}-\mathcal{U}_{\delta}) = - \mathcal{M}^{\ep}(\delta;\mathcal{U}^{\ep}_{\delta}, \mathcal{U}_{\delta}),
\]
Here
\[
\mathcal{M}^{\ep} (\delta;\mathcal{U}^{\ep},\mathcal{U}_{\delta}) = 
\begin{pmatrix}
\mathcal{M}^{\ep, 1}(\delta;\mathcal{U}^{\ep}_{\delta}, \mathcal{U}_{\delta}) \\ 
\mathcal{M}^{\ep, 2}(\delta;\mathcal{U}^{\ep}_{\delta}, \mathcal{U}_{\delta}) \\ 
\mathcal{Q}^{\ep} H^{\ep} (\delta; \tilde{\omega}^{\ep}_{\delta}, \tilde{\omega}_{\delta}, U_{\delta}) + \mathcal{M}^{\ep, 3}(\delta;\mathcal{U}^{\ep}_{\delta}, \mathcal{U}_{\delta}) 
\end{pmatrix},
\]
where
\begin{align*}
\mathcal{M}^{\ep, 1} (\delta;\mathcal{U}^{\ep}_{\delta},\mathcal{U}_{\delta}) & = \tilde \eta_{\delta} \re\{[Kz_0^\ep]_{+,\ep} - [\b K\b z_0]_+\} + \delta \re \left\{ \tilde{\eta}^{\ep}_{\delta} [K U^{\ep}_{\delta}]_{+,\ep} -\tilde{\eta}_{\delta} [\b K \b U_{\delta}]_ +\right\}  \\
& \quad + \re \left\{ [M(z^{\ep}_{0}, U_{\delta}) ]_{+, \ep} - [\b{M}(\b{z}_{0}, \b{U}_{\delta}) ]_{+} \right\} \\
& \quad +\delta \re\left\{ [N(U^\ep_{\delta})]_{+,\ep} -[\b N (\b U_{\delta})]_+ \right\} \\
\mathcal{M}^{\ep, 2} (\delta;\mathcal{U}^\ep,\mathcal{U})&= (a^{\ep} - a ) \tilde{\omega}_{\delta} 
\\ & \quad + \tilde \eta_{\delta} \im\{[Kz_0^\ep]_{+,\ep} - [\b K\b z_0]_+\} + \delta \im \left\{ \tilde{\eta}^{\ep}_{\delta} [K U^{\ep}_{\delta}]_{+,\ep} -\tilde{\eta}_{\delta} [\b K \b U_{\delta}]_ +\right\}  \\
& \quad + \im \left\{ [M(z^{\ep}_{0}, U_{\delta}) ]_{+, \ep} - [\b{M}(\b{z}_{0}, \b{U}_{\delta}) ]_{+} \right\} \\
& \quad +\delta \im \left\{ [N(U^\ep_{\delta})]_{+,\ep} -[\b N (\b U_{\delta})]_+ \right\} \\
\mathcal{M}^{\ep, 3} (\delta;\mathcal{U}^\ep,\mathcal{U})& = \mathcal{Q}^{\ep} \left\{ \delta (\tilde{\eta}^{\ep}_{\delta} K z_0^\ep - \tilde{\eta}_{\delta} K z_0) + \delta^{2} ( \tilde{\eta}^{\ep}_{\delta}  K U^{\ep}_{\delta} - \tilde{\eta}_{\delta} K U_{\delta} )  \right. \\ 
& \quad \quad \quad \left.  + N(z_0^{\ep} + \delta U^{\ep}_{\delta}) -  N(z_0 + \delta U_{\delta}) \right\}. \\
\end{align*}

To estimate $\mathcal{U}^{\ep}_{\delta} - \mathcal{U}_{\delta}$, we use the following lemma. 

\begin{lem}\label{lem6.3}
Let $u_{j} = \trans(\phi_{j}, \bu_{j})$, $\tilde{u}_{j} = \trans(\tilde{\phi}_{j}, \tilde{\bu}_{j})$, $j =1, 2$. Then the following estimates hold uniformly for $0 < \ep \le \ep_{0}${\rm :} 
\begin{itemize}
\item[\rm (i)] 
$\begin{aligned}[t]
\left| [u_1]_{+,\ep}-[\bu_2]_{+} \right|  & \le C \left\{ \ep^{2} \|\phi_{1} \|_{L^2(\mathbb{T}_{\frac{2\pi}{a}};L^{2}(\Omega) ) } + \ep^2\|\bu_2\|_{L^2(\mathbb{T}_{\frac{2\pi}{a}};(\b X^1)^*)}  \right. \\ 
& \quad \quad \quad \left. + \|\bu_1-\bu_2\|_{L^2(\mathbb{T}_{\frac{2\pi}{a}};(\b X^1)^*)} \right\},
\end{aligned}$
\item[\rm (ii)] 
$\begin{aligned}[t]
& \left| [N(u_{1}, u_{2})]_{+,\ep}-[\b{N}(\tilde{\bu}_{1}, \tilde{\bu}_{2})]_{+} \right| \\[1ex] 
& \quad \le C \left\{ \ep^2|||u_{1} |||_{\ep, \mathcal{Y}_{a} } |||u_{2} |||_{\ep, \mathcal{Y}_{a} } 
%\right. \\ 
%& \quad \quad 
+ ||| u_{1} |||_{\ep, \mathcal{Y}_{a} } ||| u_{2} - \tilde{u}_{2} |||_{\ep, \mathcal{Y}_{a} } \right. \\[1ex] 
& \quad \quad \quad \quad \left. + ||| u_{1} - \tilde{u}_{1} |||_{\ep, \mathcal{Y}_{a} } ||| \tilde{u}_{2} |||_{\ep, \mathcal{Y}_{a} } \right\}, 
\end{aligned}$
%\item[\rm (iii)] 
%$\begin{aligned}[t]
%\opnorm{\mathcal{Q}^\ep u_1 - \tilde{\mathcal{Q}}_0 u_2}_{\ep,\mathcal{X}_a} & \le C\{\ep\opnorm{u_2}_{\ep,\mathcal{X}_a} + \opnorm{u_1-u_2}_{\ep,\mathcal{X}_a}\},   
%\end{aligned}$ 
%
%where $\tilde{\mathcal{Q} }_{0} = I - \tilde{\mathcal{P} }_{0} $ with $\tilde{\mathcal{P} }_{0} u = \trans(0, \hat{ \mathcal{P} }_{0}\bu)$ for $u = \trans(\phi, \bu)$.   
\item[\rm (iii)] 
$\begin{aligned}[t]
|||\mathcal{Q}^{\ep} H^{\ep} (\delta; \tilde{\omega}^{\ep}_{\delta}, \tilde{\omega}_{\delta}, U_{\delta})|||_{\ep, \mathcal{X}_{a} }
 \le C\left\{ \ep (1 + |\delta| |\tilde{\omega}_{\delta} | + \delta^{2} |\tilde{\omega}^{\ep}_{\delta} |) + \delta^{2} | \tilde{\omega}^{\ep}_{\delta} - \tilde{\omega}_{\delta} |\right\}. 
\end{aligned}$
\end{itemize}
\end{lem}

We will give a proof of Lemma \ref{lem6.3} in the end of this section. 

\vspace{2ex}
Noting that $Ku = \trans(0, \b K \bu)$ for $u=\trans(\phi,\bu)$, we apply Theorems \ref{thm4.1}, \ref{thm4.2} (i), Lemmas \ref{lem5.8}, \ref{lem6.1}, \ref{lem6.2} and \ref{lem6.3} to obtain 
\begin{equation}\label{6.18} 
\begin{split}
&| \tilde \eta^\ep_{\delta} -  \tilde \eta_{\delta}| + | \tilde \omega^\ep_{\delta} -  \tilde \omega_{\delta}|  \\
& \quad \le C \{\ep^{2} + |\delta|(| \tilde \eta^\ep_{\delta} -  \tilde \eta_{\delta}| + | \tilde \omega^\ep_{\delta} -  \tilde \omega_{\delta}|) + \opnorm{U^\ep_{\delta} - U_{\delta}}_{\ep,\mathcal{Y}_a} \},   
\end{split}
\end{equation}
and 
\begin{equation}\label{6.19} 
\opnorm{\mathcal{Q}^{\ep} (U^\ep_{\delta} - U_{\delta}) }_{\ep,\mathcal{Y}_a} \\
\le C\{ \ep + |\delta|(| \tilde \eta^\ep_{\delta} -  \tilde \eta_{\delta}| + | \tilde \omega^\ep_{\delta} -  \tilde \omega_{\delta}|) + |\delta|\opnorm{ U^\ep_{\delta} - U_{\delta} }_{\ep,\mathcal{Y}_a} \}   
\end{equation}
uniformly in $0 < \ep \le \ep_{0}$ and $|\delta| \le \delta_{0}$. 

Since $[\b{U}_{\delta} ]_{+} = (\b{U}_{\delta}, \b{z}^{*}_{+} ) = 0$, $[U_{\delta} ]_{+, \ep} = \ep^{2}(\Phi_{\delta}, p^{\ep*}_{+}) + (\b{U}_{\delta}, \b{z}^{\ep*}_{+})$ and  $\b{z}_{+}^{\ep*} = \b{z}^{*}_{+} + O(\ep^{2})$, we have 
\[
|[U_{\delta} ]_{+, \ep}| 
\le C \{ \ep^{2} \|\Phi_{\delta} \|_{L^{2}(\mathbb{T}_{\frac{2\pi}{a} }; L^2) } + \ep^{2} \| \b{U}_{\delta} \|_{L^{2}(\mathbb{T}_{\frac{2\pi}{a} }; (\b{X}^{1})^{*}) } \} 
\le C \ep^{2} \| U_{\delta} \|_{\mathcal{X}_{a} }
\]
uniformly for $0 < \ep \le \ep_{0}$ and $|\delta| \le \delta_{0}$, and hence,  
\[ 
|\mathcal{P}^{\ep} U_{\delta} | = | 2\re ( [U_{\delta} ]_{+, \ep} z^{\ep}_{+} ) | \le C\ep^{2} \| U_{\delta} \|_{\mathcal{X}_{a} } |z^{\ep}_{+}| .
\]
Combining this with $\mathcal{P}^\ep U_\delta^\ep = 0$, we see that 
\begin{equation*}
\begin{split}
|||U^\ep_{\delta} - U_{\delta}|||_{\ep, \mathcal{Y}_a} & \le |||\mathcal{Q}^\ep (U^\ep_{\delta} - U_{\delta} )|||_{\ep, \mathcal{Y}_a} + |||\mathcal{P}^\ep (U^\ep_{\delta} - U_{\delta} )|||_{\ep, \mathcal{Y}_a} \\
& = |||\mathcal{Q}^\ep (U^\ep_{\delta} - U_{\delta} )|||_{\ep, \mathcal{Y}_a} + |||\mathcal{P}^\ep U_{\delta} |||_{\ep, \mathcal{Y}_a} \\
& \le |||\mathcal{Q}^\ep (U^\ep_{\delta} - U_{\delta})|||_{\ep, \mathcal{Y}_a} +C \ep^2 \| U_{\delta} \|_{\mathcal{X}_{a} } \\
& \le |||\mathcal{Q}^\ep (U^\ep_{\delta} - U_{\delta})|||_{\ep, \mathcal{Y}_a} +C \ep^2.  
\end{split}
\end{equation*}
This, together with \eqref{6.18} and \eqref{6.19}, implies that the estimate 
\begin{equation}\label{6.19'} 
\opnorm{\mathcal{U}^\ep_{\delta} -\mathcal{U}_{\delta} }_{\ep,\mathbb{R}\times\mathbb{R}\times\mathcal{Y}_a} \le C\ep
\end{equation}
holds uniformly for $0<\ep\le \ep_{0}$ and $|\delta| \le \delta_{0}$. 

Before proceeding further, we make an observation. It is not difficult to see that $\|\delt^{2} U_{0}\|_{\mathcal{X}_{a} } \le C$, from which one can also prove, in a similar manner to above, that 
\begin{equation}\label{6.19''} 
||| \delt (U^{\ep}_{0} - U_{0}) |||_{\ep, \mathcal{Y}_{a} } \le C\ep. 
\end{equation}

It remains to show that 
\begin{equation}\label{6.29}
|\tilde \eta_1^\ep(\delta) - \tilde{\eta}_1(\delta)| + |\tilde \omega_1^\ep(\delta) - \tilde{\omega}_1(\delta)| \le C \ep. 
\end{equation}

To this end, we shall estimate $\mathcal{V}^\ep_{\delta} = (\mathcal{U}^\ep_{\delta} - \mathcal{U}_0^\ep) - (\mathcal{U}_{\delta} - \mathcal{U}_0)$. Since $\mathcal{V}^\ep_{\delta} = (\mathcal{U}^\ep_{\delta} - \mathcal{U}_{\delta}) - (\mathcal{U}^{\ep}_{0} - \mathcal{U}_0)$, we have 
\begin{align*}
\mathcal{B}^\ep (\delta^2 \tilde{\omega}^\ep_{\delta}) \mathcal{Q}^{\ep} \mathcal{V}^\ep_{\delta} 
& 
= - \frac{a^{\ep} }{a} \delta^{2} \tilde{\omega}^{\ep}_{\delta} 
\begin{pmatrix}
0 \\ 
0 \\
\mathcal{Q}^{\ep}\delt(U^{\ep}_{0} - U_{0}) 
\end{pmatrix} 
\\[2ex]
& \quad \quad  
- \left( \mathcal{M}^\ep (\delta; \mathcal{U}^\ep_{\delta}, \mathcal{U}_{\delta}) - \mathcal{M}^{\ep} (0; \mathcal{U}^{\ep}_{0}, \mathcal{U}_{0} ) \right).  
\end{align*}
By using \eqref{6.19''}, we have $||| \frac{a^{\ep} }{a} \delta^{2} \tilde{\omega}^{\ep}_{\delta} \delt(U^{\ep}_{0} - U_{0}) |||_{\ep, \mathcal{Y}_{a} } \le C \ep |\delta|$. 
%+ (H^{\ep}(\delta; \tilde{\omega}^{\ep}_{\delta}, \tilde{\omega}_{\delta}, U_{\delta}) - H^{\ep}(0; \tilde{\omega}^{\ep}_{0}, \tilde{\omega}_{0}, U_{0}) ) |||_{\ep, \mathcal{X}_{a} } \le C \ep |\delta|.  \]
Based on this, together with \eqref{6.27} and \eqref{6.19'}, we have, in a similar manner to above, 
\begin{align*}
& |(\tilde{\eta}^{\ep}_{\delta} - \tilde{\eta}^{\ep}_{0}) - (\tilde{\eta}_{\delta} - \tilde{\eta}_{0})|
+ |(\tilde{\omega}^{\ep}_{\delta} - \tilde{\omega}^{\ep}_{0}) - (\tilde{\omega}_{\delta} - \tilde{\omega}_{0})| \\
& \quad \le C\{ \ep |\delta| + ||| (U^{\ep}_{\delta} - U^{\ep}_{0}) - (U_{\delta} - U_{0}) |||_{\ep, \mathcal{Y}_{a} } \},     
\end{align*} 
and 
\[
|||\mathcal{Q}^{\ep}[ (U^{\ep}_{\delta} - U^{\ep}_{0}) - (U_{\delta} - U_{0}) ]|||_{\ep, \mathcal{Y}_{a} } \le C\ep |\delta|. 
\]
Since $[\b{U}_{\delta} - \b{U}_{0}]_{+} = 0$, we see, as above, that 
\begin{align*}
& ||| (U^{\ep}_{\delta} - U^{\ep}_{0}) - (U_{\delta} - U_{0}) |||_{\ep, \mathcal{Y}_{a} } \\
& \quad \le C\{|||\mathcal{Q}^{\ep}[ (U^{\ep}_{\delta} - U^{\ep}_{0}) - (U_{\delta} - U_{0}) ]|||_{\ep, \mathcal{Y}_{a} } + \ep^{2} ||U_{\delta} - U_{0} ||_{\mathcal{X}_{a} } \} \\
& \quad \le C\{|||\mathcal{Q}^{\ep}[ (U^{\ep}_{\delta} - U^{\ep}_{0}) - (U_{\delta} - U_{0}) ]|||_{\ep, \mathcal{Y}_{a} } + \ep^{2}|\delta| \}. 
\end{align*} 
We thus conclude that $|||\mathcal{V}^\ep|||_{\ep, \mathbb{R} \times \mathbb{R} \times \mathcal{Y}_a} \le C \ep |\delta|$ uniformly for $0< \ep \le \ep_{0}$ and $|\delta | \le \delta_{0}$. 
In particular, we obtain \eqref{6.29}. This completes the proof. 
\hfill$\square$

\subsection{Proof of Theorem \ref{thm4.2} (iii)} 
We next prove Theorem \ref{thm4.2} (iii), i.e., the uniqueness of time periodic solutions up to translations in the time variable $t$. 

We introduce notation. We set ${z}_{0}^{\ep, \tau} = 2\re(e^{ia\tau}z_{+}^{\ep})$ and $({z}_{0}^{\ep*})^{\tau} = 2\re(e^{ia\tau}z_{+}^{\ep*})$ for a constant $\tau \in [- \frac{\pi}{a}, \frac{\pi}{a})$ and define $[u]_{+,\ep}^{(\tau)}$ by $[u]_{+,\ep}^{(\tau)} =\langle u, e^{ia\tau}{z}_{+}^{\ep*} \rangle_{\ep} = e^{- ia\tau} [u]_{+,\ep}$. We note that $\mathcal{P}_{\pm}^\ep$ can be expressed as 
\[
\mathcal{P}_{\pm}^\ep u = \langle u, e^{\pm i a\tau} z_{\pm}^{\ep*} \rangle_{\ep} \,e^{\pm i a\tau} z_{\pm}^{\ep}. 
\]

\vspace{2ex}
\noindent
{\bf Proof of  Theorem \ref{thm4.2} (iii).} Let $u^\ep \in\mathcal{Y}_a$ be a solution of \eqref{4.6} with $\eta = \eta^{\ep}$ and $\omega = \omega^{\ep}$. Then, by Theorem \ref{thm5.6}, $u^\ep $ is written as 
\begin{equation}\label{6.30}
u^\ep  = u_0^\ep  + U^\ep,  \ \ u_0^\ep =\mathcal{P}^{\ep} u^{\ep}  \in N(B^{\ep}),\,\, U^\ep \in R(B^{\ep})
\end{equation}

Assume first that $u_0^\ep  = 0$ which is equivalent to $[u^{\ep}]_{+, \ep} = 0$. We then find that 
\[
B^{\ep} U^\ep  + \frac{a^{\ep} }{a} \omega^\ep \delt U^\ep  + \eta^\ep  K U^\ep  + N(U^\ep ) = 0,
\]
and hence, 
\[
B^{\ep}(\omega^\ep ) U^\ep  + \eta^\ep  \mathcal{Q}^\ep K U^\ep  + \mathcal{Q}^\ep N (U^\ep ) = 0.
\]
By using Lemmas \ref{lem5.8} and \ref{lem6.1}, we see that if $|\omega^{\ep} |\le\frac{1}{4}$, then
\[
\opnorm{U^\ep }_{\ep,\mathcal{Y}_a} \le C\{|\eta^\ep |\opnorm{U^\ep }_{\ep,\mathcal{Y}_a} + \opnorm{U^\ep }_{\ep,\mathcal{Y}_a}^{\frac{3}{2}} + \opnorm{U^\ep }_{\ep,\mathcal{Y}_a}^2\}.
\]
This implies that there exists a positive constant $r_0$ such that if $|\eta^\ep | + \opnorm{U^\ep }_{\ep,\mathcal{Y}_a} \le r_0$, then 
$\opnorm{U^\ep }_{\ep,\mathcal{Y}_a} \le \frac{1}{2}\opnorm{U^\ep }_{\ep,\mathcal{Y}_a}$, and therefore, $U^\ep  = 0$. We thus conclude that if $|\eta^\ep | + \opnorm{u^\ep }_{\ep,\mathcal{Y}_a} \le r_0$, $|\omega^\ep |\le\frac{1}{4}$ and $u^\ep  = U^\ep \in R(B^{\ep})$, then $u^\ep  = 0$.

Assume next that $u_0^\ep \neq 0$ in \eqref{6.30}. Then there exists a constant $\tau \in [- \frac{\pi}{a}, \frac{\pi}{a} )$ such that $u_0^{\ep} = \mathcal{P}^{\ep} u^{\ep}$ is written as $u_0^\ep  = \delta {z}_{0}^{\ep,\tau} $, where $\delta =| [u^{\ep}]_{+, \ep}|$ and ${z}_{0}^{\ep, \tau} = 2 \re (e^{ia\tau} z_{+}^{\ep})$. It then follows that $u^\ep $ is written as 
\[
u^\ep  = \delta {z}_{0}^{\ep, \tau}  + U^\ep, \ \ U^\ep \in R(B^{\ep}).
\]
We substitute this into \eqref{4.6} with $\eta = \eta^{\ep}$ and $\omega = \omega^{\ep}$. Since $[K {z}_{0}^{\ep,\tau} ]_{+, \ep}^{(\tau)} = [K z_{0}^{\ep} ]_{+, \ep}$, as in the proof of Theorem \ref{thm4.2} (i), we obtain
\[
\delta\eta^\ep \re[K z_0^\ep ]_{+,\ep} + \eta^\ep \re[K U^\ep ]_{+,\ep}^{(\tau)} + \re[N(\delta {z}_{0}^{\ep, \tau}  + U^\ep )]_{+,\ep}^{(\tau)} = 0,
\]
\[
\delta a^\ep  \omega^\ep  + \delta \eta^\ep \im[K z_0^\ep ]_{+,\ep} + \eta^\ep \im[K U^\ep ]_{+,\ep}^{(\tau)} + \im[N(\delta {z}_{0}^{\ep, \tau}  + U^\ep )]_{+,\ep}^{(\tau)} = 0,
\]
\[
B^{\ep}(\omega^\ep ) U^\ep  + \delta \eta^\ep \mathcal{Q}^\ep K {z}_{0}^{\ep, \tau}  + \eta^\ep \mathcal{Q}^\ep K U^\ep  + \mathcal{Q}^\ep N (\delta {z}_{0}^{\ep, \tau}  + U^\ep ) = 0.
\]
Setting $\hat{\mathcal{V}}_\delta^\ep  = \trans(\delta \eta^\ep ,\delta  \omega^\ep ,U^\ep )$, we have
\[
\hat{ \mathcal{B} }^{\ep} (\omega^\ep )\hat{\mathcal{V}}_\delta^\ep 
= -
\begin{pmatrix}
\eta^\ep \re[K U^\ep ]_{+,\ep}^{(\tau)} + \re[N (\delta {z}_{0}^{\ep, \tau} + U^\ep )]_{+,\ep}^{(\tau)} \\[1ex]
\eta^\ep \im[K U^\ep ]_{+,\ep}^{(\tau)} + \im[N (\delta {z}_{0}^{\ep, \tau} + U^\ep )]_{+,\ep}^{(\tau)} \\[1ex]
\eta^\ep \mathcal{Q}^\ep K U^\ep  + \mathcal{Q}^\ep N (\delta {z}_{0}^{\ep, \tau} + U^\ep )
\end{pmatrix}, 
\]
where $\hat{\mathcal{B}}^{\ep}(\omega):\mathbb{R}\times \mathbb{R}\times\mathcal{Q}^\ep \mathcal{Y}_a\rightarrow\mathbb{R}\times\mathbb{R}\times\mathcal{Q}^\ep \mathcal{X}_a$ is defined by
\[
\hat{ \mathcal{B} }^{\ep} (\omega^\ep ) 
= 
\begin{pmatrix}
\re[K z_0^\ep]_{+,\ep} & 0 & 0\\
\im[K z_0^\ep]_{+,\ep} & a^\ep & 0\\
\mathcal{Q}^\ep K {z}_{0}^{\ep, \tau} & 0 & B^{\ep}(\omega^\ep)\mathcal{Q}^\ep
\end{pmatrix}.
\]
Applying Lemmas \ref{lem5.8} and \ref{lem6.2}, we see that
\[
\opnorm{\hat{\mathcal{V}}_\delta^\ep }_{\ep,\mathbb{R}\times\mathbb{R}\times\mathcal{Y}_a} \le C\{|\delta|^2 + |\delta| \opnorm{U^\ep }_{\ep,\mathcal{Y}_a} + |\eta^\ep |\opnorm{U^\ep }_{\ep,\mathcal{Y}_a} + \opnorm{U^\ep }_{\ep,\mathcal{Y}_a}^{\frac{3}{2}} + \opnorm{U^\ep }_{\ep,\mathcal{Y}_a}^2\}.
\]
This implies that there exists a positive number $r_1$ such that if $|\eta^\ep | +\opnorm{U^\ep }_{\ep,\mathcal{Y}_a} \le r_{1}$ and $|\omega^{\ep} | \le \frac{1}{4}$, then 

\begin{equation}\label{6.33}
\opnorm{\hat{\mathcal{V}}_\delta^\ep }_{\ep,\mathbb{R}\times\mathbb{R}\times\mathcal{Y}_a} \le C_1|\delta|^{2}. 
\end{equation}
This suggests us to write
\begin{equation}\label{6.34}
\eta^\ep = \delta \hat \eta^\ep, \ \ \ \omega^\ep = \delta \hat \omega^\ep 
\end{equation}
and 
\begin{equation}\label{6.35}
u = \delta( {z}_{0}^{\ep, \tau} + \delta \hat U^\ep ), \ \  \hat U^\ep \in R(B^{\ep}).
\end{equation}
We note that, by \eqref{6.33}, $\hat{\mathcal{U}}^\ep = \trans(\hat \eta^\ep, \hat \omega^\ep, \hat U^\ep)$ with $\hat \eta^\ep$, $\hat \omega^\ep$ and $\hat U^\ep$ given in \eqref{6.34} and \eqref{6.35} satisfy
\[
\opnorm{\hat{\mathcal{U}}^\ep}_{\ep,\mathbb{R}\times\mathbb{R}\times\mathcal{Y}_a} \le C_{1}.
\]
As in the proof of Theorem \ref{thm4.2} (i), we also find that $\delt \hat{U}^{\ep} \in \mathcal{Y}_{a}$ with estimate $||| \delt \hat{U}^{\ep} |||_{\mathcal{Y}_{a} } \le C_{2}$ uniformly in $0 < \ep \le \ep_{0}$ and $|\delta| \le \delta_{0}$. 

We thus find that if $u^\ep = \delta {z}_{0}^{\ep, \tau} + U^\ep $, $|\delta| \le \delta_{0}$, $|\eta^\ep | + \opnorm{U^\ep }_{\ep,\mathcal{Y}_a} \le r_1$ and $|\omega^\ep | \le \frac{1}{4}$, then $\trans(\eta^\ep , \omega^\ep , U^\ep )$ belongs to the set  
\[
\{ \trans( \delta \hat \eta^\ep , \delta \hat \omega^\ep ,\delta^{2} \hat U^\ep); \opnorm{(\hat \eta^\ep, \hat \omega^\ep, \hat U^\ep)}_{\ep,\mathbb{R}\times\mathbb{R}\times\mathcal{Y}_a} \le C_{1},  ||| \delt \hat{U}^{\ep} |||_{\mathcal{Y}_{a} } \le C_{2} \}. 
\]

We next claim that if $\hat{\mathcal{U}}_j^\ep = \trans (\hat \eta^\ep_j, \hat \omega_j^\ep,  \hat U_j^\ep)$, $j=1,2$, are solutions of \eqref{6.2} with $z_{0}^{\ep}$ replaced by $z_{0}^{\ep,\tau}$ for $0 < |\delta| \le \delta_{0}$ and $\hat{\mathcal{U}}_j^\ep  = \trans(\hat{\eta}_{j}^{\ep}, \hat{\omega}_{j}^{\ep}, \hat{U}_{j}^{\ep})$ satisfy $\opnorm{\hat{\mathcal{U}}_j^\ep }_{\ep,\mathcal{Y}_a} \le C_{1}$ and $ ||| \delt \hat{U}^{\ep}_{j} |||_{\mathcal{Y}_{a} } \le C_{2}$ $(j=1,2)$, then it holds $\hat{\mathcal{U}}_1^\ep = \hat{\mathcal{U}}_2^\ep$. 

Indeed, we see that $\hat{\mathcal{U}}_j^\ep$ $(j=1,2)$ satisfy
\[
\hat{\mathcal{B} }^{\ep} (\delta \hat \omega_j^\ep)\hat{\mathcal{U}}_j^\ep = - \hat{\mathcal{N} }^{\ep} (\delta; \hat{\eta}^{\ep}_{j}, \hat{U}_j^\ep),   
\]
where 
\[
\hat{\mathcal{N} }^{\ep}( \delta; \hat{\eta}^{\ep}, \hat{U}^{\ep}) 
=
\begin{pmatrix}
\delta \hat{\eta}^\ep \re[K \hat U^\ep ]_{+,\ep}^{(\tau)} + \re[N({z}_{0}^{\ep, \tau}  + \delta \hat{U}^\ep )]_{+,\ep}^{(\tau)} \\[1ex]
\delta \hat{\eta}^\ep \im[K \hat{U}^\ep ]_{+,\ep}^{(\tau)} + \im[N({z}_{0}^{\ep, \tau}  + \delta \hat{U}^\ep) ]_{+,\ep}^{(\tau)} \\[1ex]
\delta\hat{\eta}^\ep \mathcal{Q}^\ep K \hat{U}^\ep  + \mathcal{Q}^\ep N({z}_{0}^{\ep, \tau}  + \delta \hat{U}^\ep )
\end{pmatrix}
\]

%where $\mathcal{B}^{\ep}_{\tau}(\omega)$ is the map with $\mathcal{Q}^{\ep} K z_{0}^{\ep}$ in $\mathcal{B}^{\ep}(\omega)$ replaced by $\mathcal{Q}^{\ep} K \tilde{z}_{0}^{\ep}$. 
It follows that
\begin{align*}
\hat{ \mathcal{B} }^{\ep} (\delta\hat \omega_1^\ep)(\hat{\mathcal{U}}_1^\ep-\hat{\mathcal{U}}_2^\ep) 
& = \frac{a^\ep}{a}\delta(\hat \omega_2^\ep-\hat \omega_1^\ep) \begin{pmatrix}0 \\ 0 \\ \delt \hat{U}_2^\ep \end{pmatrix} 
%+ \delta(\hat{\eta}^{\ep}_{2} - \hat{\eta}^{\ep}_{1}) 
%\begin{pmatrix} 
%\re[K \hat{U}^{\ep}_{2} ]_{+, \ep} \\
%\im[K \hat{U}^{\ep}_{2} ]_{+, \ep} \\
%\mathcal{Q}^{\ep} K \hat{U}^{\ep}_{2} 
%\end{pmatrix}
%\\[1ex]
%& \quad 
+ \hat{\mathcal{N} }^{\ep} (\delta; \hat{\eta}_{1}^\ep, \hat{U}_1^\ep)-\hat{\mathcal{N}}^{\ep}(\delta; \hat{\eta}^{\ep}_{2}, \hat{U}_2^\ep).
\end{align*} 
By using Lemmas \ref{lem5.8} and \ref{lem6.1}--\ref{lem6.3}, we find that 
\begin{align*}
\opnorm{\hat{\mathcal{U}}_1^\ep - \hat{\mathcal{U}}_2^\ep }_{\ep,\mathbb{R}\times\mathbb{R}\times\mathcal{Y}_a} 
& \le C \{|\delta| |\hat \omega_1^\ep - \hat \omega_2^\ep | C_{2} + C_{1} |\delta|(|\hat \omega_1^\ep - \hat \omega_2^\ep | + |\hat \eta_1^\ep - \hat \eta_2^\ep|) 
\\ 
& \quad \quad + C_{1} |\delta|\opnorm{\hat U_1^\ep - \hat U_2^\ep }_{\ep,\mathcal{Y}_a}  + |\delta| \opnorm{\hat U_1^\ep - \hat U_2^\ep }_{\ep,\mathcal{Y}_a} \\
&\quad \quad  + (C_{1}^{\frac{1}{2}} + C_{1} )|\delta|^2\opnorm{\hat U_1^\ep - \hat U_2^\ep  }_{\ep,\mathcal{Y}_a} \} \\
& \le C|\delta| \opnorm{\hat{\mathcal{U}}_1^\ep - \hat{\mathcal{U}}_2^\ep }_{\ep,\mathbb{R}\times\mathbb{R}\times\mathcal{Y}_a}.
\end{align*} 
We thus conclude that there exists a positive number $\delta_{0}$ such that if $|\delta|\le\delta_{0}$, then $ \opnorm{\hat{\mathcal{U}}_1^\ep - \hat{\mathcal{U}}_2^\ep }_{\ep,\mathbb{R}\times\mathbb{R}\times\mathcal{Y}_a} =0$, i.e., $\hat{\mathcal{U}}_1^\ep = \hat{\mathcal{U}}_2^\ep$. 

Let $\trans(\eta^{\ep}_{\delta}, \omega^{\ep}_{\delta}, u^{\ep}_{\delta})$ be the solution branch obtained in Theorem \ref{thm4.2} (i). 
We denote $\eta^{\ep}_{\delta} = \delta \hat{\eta}^{\ep}_{\delta}$, $\omega^{\ep}_{\delta} = \delta \hat{\omega}^{\ep}_{\delta}$ and $u^{\ep}_{\delta} = \delta(z^{\ep}_{0} + \delta U^{\ep}_{\delta} )$. 
Set $\hat{\mathcal{U} }^{\ep}_{\delta} = \trans(\hat{\eta}^{\ep}_{\delta}, \hat{\omega}^{\ep}_{\delta}, \hat{U}^{\ep}_{\delta})$. 
We then see that  $|\eta^{\ep}_{\delta} | + |||U^{\ep}_{\delta} |||_{\ep, \mathcal{Y}_{a} } \le r_{1}$ and $|\omega^{\ep}_{\delta}| \le \frac{1}{4}$ by taking $\delta_{0}$ in the proof of Theorem \ref{thm4.2} (i) suitably smaller if necessary. 
The above argument is then applicable to $\hat{\mathcal{U}}^{\ep}_{\delta}$ to  obtain $||| \hat{\mathcal{U} }^{\ep}_{\delta} |||_{\ep, \mathbb{R} \times \mathbb{R} \times \mathcal{Y}_{a} } \le C_{1}$ and $||| \delt \hat{U}^{\ep}_{\delta} |||_{\ep, \mathcal{Y}_{a} } \le C_{2}$. 
We set $\hat{\mathcal{U}}_{\delta}^{\ep, \tau} = \trans(\hat{\eta}^{\ep}_{\delta}, \hat{\omega}^{\ep}_{\delta}, U^{\ep}_{\delta} (\cdot + \tau) )$. Since $u^{\ep}_{\delta}(\cdot + \tau)$ takes the form 
\[
u^{\ep}_{\delta}(\cdot + \tau) = \delta( {z}_{0}^{\ep, \tau}(\cdot) + \delta U^{\ep}_{\delta} (\cdot + \tau))
\] 
and $\hat{\mathcal{U}}_{\delta}^{\ep, \tau}$ satisfies $|||\hat{\mathcal{U} }_{\delta}^{\ep, \tau} |||_{\ep, \mathbb{R} \times \mathbb{R} \times \mathcal{Y}_{a} } \le C_{1}$ and $||| \delt \hat{U}^{\ep, \tau}_{\delta} |||_{\ep, \mathcal{Y}_{a} } \le C_{2}$, we see  that $\hat{\mathcal{U}}^{\ep} = \hat{\mathcal{U} }_{\delta}^{\ep, \tau}$, and hence, we have $\trans(\eta^{\ep}, \omega^{\ep}, u^{\ep}(\cdot) ) =\trans(\eta^{\ep}_{\delta}, \omega^{\ep}_{\delta}, u^{\ep}_{\delta}(\cdot + \tau) )$. This completes the proof. \hfill$\square$

\vspace{2ex}
We finally give proofs of Lemmas \ref{lem6.2} and \ref{lem6.3}.  

\vspace{1ex}
\noindent
{\bf Proof of Lemma \ref{lem6.2}.} 
By the definition of $[u]_{+,\ep}$ , we have 
\begin{align*}
|[u]_{+,\ep}| & \ \le \frac{a}{2\pi}\int_0^{\frac{2\pi}{a}} \ep^2|(\phi,p_+^{\ep*})_{L^2}| + |(\bu,\b z_+^{\ep*})|\,dt\\
&\ \le C\{\ep^2\|\phi\|_{L^2(\mathbb{T}_{\frac{2\pi}{a}};L^2(\Omega))} + \|\bu\|_{L^2(\mathbb{T}_{\frac{2\pi}{a}};(\b X^1)^*)}\}\\
&\ \le C\opnorm{u}_{\ep,\mathcal{X}_a}.
\end{align*}
This proves (i).

As for (ii), by Theorem \ref{thm4.1}, we have $\sup_{0\le t \le \frac{2\pi}{a}} \| z^{\ep}_{+}(t) \|_{L^{\infty}(\Omega) } \le C \| u^{\ep}_{+} \|_{ H^{2}(\Omega) } \le C$ uniformly in $0 < \ep \le \ep_{0}$. It follows that 
\begin{align*}
|[N(u_{1}, u_{2} ) ]_{+,\ep}| & \le \frac{a}{2\pi}\int_0^{\frac{2\pi}{a}} \|\b{N}(\bu_{1}, \bu_{2} )\|_{1} \| \b{z}^{\ep}_{+} \|_{\infty} \,dt \\
& \le C \int_0^{\frac{2\pi}{a}}  \|\bu_{1}\|_{2} \|\nabla \bu_{2} \|_{2} \, dt \\
& \le C\left(\sup_{0\le t \le \frac{2\pi}{a} } \| \bu_1(t)\|_2^2\right)^{\frac{1}{2} }\left(\int_0^{\frac{2\pi}{a}}  \| \nabla \bu_{2}\|_2^2 \, dt \right)^{\frac{1}{2} } \\
& \ \le C ||| u_1 |||_{\ep, \mathcal{Y}_a} ||| u_2 |||_{\ep, \mathcal{Y}_a}. 
\end{align*}
This shows (ii).  Similarly, by using $\sup_{0 \le t \le \frac{2\pi}{a} } \|\b{z}^{\ep}_{+}(t) \|_{\infty} + \|\nabla \b{z}^{\ep}_{+}(t) \|_{\infty} \le C$, we have 
\[
\int_{0}^{\frac{2\pi}{a} } \|\b{M}(\b{z}^{\ep}_{0}, \bu) \|_{(\b{X}_{1})^{*} }^{2} + \ep^{2} \|\b{M}(\b{z}^{\ep}_{0}, \bu) \|_{2}^{2} \,dt \le C \int_{0}^{\frac{2\pi}{a} } \| \nabla \bu \|_{2}^{2}\,dt,  
\]
which gives (iii). 

As for (iv), by the Gagliardo-Nirenberg inequality, we have 
\begin{align*}
\left|(\b N(\bu_1,\bu_2),\bu_3)\right| & \le \| \bu_1\|_4 \|\nabla \bu_2\|_2 \|\bu_3\|_4 \\
& \le \| \bu_1\|_2^{\frac{1}{2}} \| \nabla \bu_1\|_2^{\frac{1}{2}} \|\nabla \bu_2 \|_2 \|\bu_3\|_2^{\frac{1}{2}} \| \nabla \bu_3\|_2^{\frac{1}{2}}, 
\end{align*}
and hence, 
\begin{align*}
& \int_0^{\frac{2\pi}{a}} \left|(\b N(\bu_1,\bu_2),\bu_3)\right| \, dt 
\\
& \ \le C\left(\sup_{0\le t \le \frac{2\pi}{a} } \| \bu_1(t)\|_2^2 \int_0^{\frac{2\pi}{a}}  \| \nabla \bu_1\|_2^2 \, dt \right)^{\frac{1}{4} } \left( \int_0^{\frac{2\pi}{a}} \|\nabla \bu_2 \|_2^2 \, dt \right)^{\frac{1}{2}} \\
& \quad \quad  \cdot \left( \sup_{0\le t \le \frac{2\pi}{a} } \| \bu_3(t)\|_2^2\int_0^{\frac{2\pi}{a}} \| \nabla \bu_3\|_2^2\, dt \right)^{\frac{1}{4}} 
\\ 
& \ \le C ||| u_1 |||_{\ep, \mathcal{Y}_a} ||| u_2 |||_{\ep, \mathcal{Y}_a} ||| u_3 |||_{\ep, \mathcal{Y}_a}. 
\end{align*}
This proves (iv). Similarly, we have 
\[
\|\b N(\bu_1,\bu_2)\|_2^2 \, dt  
\le \|\bu_1\|_4^2 \|\nabla \bu_2\|_4^2 
\le C \|\bu_1\|_2 \|\nabla \bu_1\|_2 \|\nabla \bu_2\|_2 \|\delx^2 \bu_2\|_2,
\]
and hence, 
\begin{align*}
\int_0^{\frac{2\pi}{a}} \ep^2 \|\b N(\bu_1,\bu_2)\|_2^2 \, dt  
& \ \le C \int_0^{\frac{2\pi}{a}} \ep^2 \|\bu_1\|_2 \|\nabla \bu_1\|_2 \|\nabla \bu_2\|_2 \|\delx^2 \bu_2\|_2 \, dt
\\  
& \ \le C \left( \sup_{0\le t \le \frac{2\pi}{a} } \| \bu_1(t)\|_2^2 \int_0^{\frac{2\pi}{a}} \|\nabla \bu_1\|_2^2\, dt \right)^{\frac{1}{2}} 
\\ 
& \quad \quad \cdot \left( \sup_{0\le t \le \frac{2\pi}{a} } \ep^2 \| \nabla \bu_2(t)\|_2^2 \int_0^{\frac{2\pi}{a}}  \ep^2 \|\delx^2 \bu_2\|_2 \, dt \right)^{\frac{1}{2}}
\\  
& \ \le C ||| u_1 |||_{\ep, \mathcal{Y}_a}^2 ||| u_2 |||_{\ep, \mathcal{Y}_a}^2,
\end{align*}
and (v) is obtained. This completes the proof. 
\hfill$\square$ 

\vspace{1ex}
\noindent
{\bf Proof of Lemma \ref{lem6.3}.} 
As for (i), we have 
\begin{align*}
|[u_1]_{+,\ep}-[\bu_2]_{+}| &  \le \frac{a}{2\pi}\int_0^{\frac{2\pi}{a}} | \ep^{2}(\phi_{1}, p_{+}^{\ep *}) +  (\bu_1-\bu_2,\b z_+^{\ep*}) + (\bu_2,\b z_+^{\ep*}-\b z_+^*)|\,dt \\[1ex] 
& \le C \left\{ \ep^{2} \|\phi_{1} \|_{L^2(\mathbb{T}_{\frac{2\pi}{a}};L^{2}(\Omega) ) } + \ep^2\|\bu_2\|_{L^2(\mathbb{T}_{\frac{2\pi}{a}};(\b X^1)^*)}  \right. \\ 
& \quad \quad \quad \left. + \|\bu_1-\bu_2\|_{L^2(\mathbb{T}_{\frac{2\pi}{a}};(\b X^1)^*)} \right\}.
\end{align*}
We thus obtain (i). Similarly, we have 
\begin{align*}
& |[N(u_{1}, u_{2})]_{+,\ep}-[\b{N}(\tilde{\bu}_{1}, \tilde{\bu}_{2} )]_{+} | \\[1ex] 
& \quad = | \lang \b{N}(\bu_{1}, \bu_{2}), \b{z}^{\ep *}_{+} \rang - \lang \b{N}(\tilde{\bu}_{1}, \tilde{\bu}_{2}), \b{z}^{*}_{+} \rang | \\[1ex]
& \quad = | \lang \b{N}(\bu_{1}, \bu_{2}) , \b{z}^{\ep*}_{+} - \b{z}^{*}_{+} \rang + \lang \b{N}(\bu_{1}, \bu_{2}) - \b{N}(\tilde{\bu}_{1}, \tilde{\bu}_{2}), \b{z}^{*}_{+} \rang  | \\[1ex]
& \quad \le C \ep^{2} \int_{0}^{\frac{2\pi}{a} }\| \b{N}(\bu_{1}, \bu_{2}) \|_{1}\, dt 
+ C \int_{0}^{\frac{2\pi}{a} } \| \b{N}(\bu_{1}, \bu_{2}) - \b{N}(\tilde{\bu}_{1}, \tilde{\bu}_{2}) \|_{1} \, dt 
\\[1ex]
& \quad \le C \left\{ \ep^2|||u_{1} |||_{\ep, \mathcal{Y}_{a} } |||u_{2} |||_{\ep, \mathcal{Y}_{a} } 
+ ||| u_{1} |||_{\ep, \mathcal{Y}_{a} } ||| u_{2} - \tilde{u}_{2} |||_{\ep, \mathcal{Y}_{a} } \right. \\[1ex] 
& \quad \quad \quad \quad \left. + ||| u_{1} - \tilde{u}_{1} |||_{\ep, \mathcal{Y}_{a} } ||| \tilde{u}_{2} |||_{\ep, \mathcal{Y}_{a} } \right\}.  
\end{align*}
This proves (ii). 
 
The inequality (iii) is proved as follows. We first consider the term $\mathcal{Q}^{\ep} \trans(0, \delt\b{z}_{0})$. Since $\hat{ \mathcal{Q} }_{0} \delt \b{z}_{0} = \delt \b{z}_{0} - 2\re( [\delt \b{z}_{0}]_{+} \b{z}_{+} ) = \b{0}$, we have 
\begin{align*}
\mathcal{Q}^\ep\trans(0, \delt\b{z}_{0})  
%&  = - (\mathcal{Q}^\ep-\tilde{\mathcal{Q}}_0)u_2 \\ 
&  = -
\begin{pmatrix}
2\re( \lang \delt\b{z}_{0}, \b{z}^{\ep*}_{+} \rang p_+^\ep) 
\\[1ex] 
2\re( \lang \delt \b{z}_{0}, \b{z}^{\ep*}_{+} \rang \b z_+^\ep - \lang \delt \b{z}_{0}, \b{z}_+^{*} \rang \b z_+)
\end{pmatrix} 
\\
&  =  
-
\begin{pmatrix}
2\re( \lang \delt \b{z}_{0}, \b{z}^{\ep*}_{+} \rang p_+^\ep) 
\\[1ex] 
2\re( \lang \delt \b{z}_{0}, \b{z}^{\ep*}_{+} - \b{z}^{*}_{+} \rang) \b z_+^\ep + \lang \delt \b{z}_{0}, \b{z}^{*}_{+} \rang (\b z_+^\ep - \b z_+)) 
\end{pmatrix}
.
\end{align*}
We thus obtain $||| \mathcal{Q}^\ep\trans(0, \delt\b{z}_{0})  |||_{\ep, \mathcal{X}_{a} } \le C\ep$. We also note that $||| \trans(\delt \Phi_{\delta}, \b{0}) |||_{\ep, \mathcal{X}_{a} } \le C \ep || \trans(\delt \Phi_{\delta}, \b{0}) ||_{\mathcal{X}_{a} } \le C \ep$. It then follows that 
\[
||| \mathcal{Q}^{\ep} H^{\ep} (\delta; \tilde{\omega}^{\ep}_{\delta}, \tilde{\omega}_{\delta}, U_{\delta})|||_{\ep, \mathcal{X}_{a} } 
\le C\left\{ \ep (1 + |\delta| |\tilde{\omega}_{\delta} | + \delta^{2} |\tilde{\omega}^{\ep}_{\delta} |) + \delta^{2} | \tilde{\omega}^{\ep}_{\delta} - \tilde{\omega}_{\delta} | 
\right\}. 
\]
This completes the proof. 
\hfill$\square$

%\newpage
\section{Stability of bifurcating solutions}\label{Stability}

Throughout this section we denote $T_a = \frac{2\pi}{a}$. In what follows, we will use the same letters $\ep_{2}$ and $\delta_{2}$ for bounds of the ranges of $\ep$ and $\delta$, even when they should be taken suitably smaller than those of the previous ones if no confusion will occur from the context.

We investigate the linearized problem around $u_\delta^\ep (t)$. The linearized problem around $u_\delta^\ep (t)$ then takes the form 
\begin{align}\label{7.1}
\frac{a^{\ep} }{a} \delt u + \frac{1}{1 + \omega^{\ep}_{\delta} } L^{\ep}_{\mathcal{R}^{\ep}_{1,c} } u + \frac{\delta}{1 + \omega^{\ep}_{\delta} } M_{\delta}^{\ep}(t) u & = F, \\ 
\label{7.2} 
u|_{t = s} &=  v,  
\end{align}
where $v \in X^{1}$ and $F \in \mathcal{X}(s, s + T_{a})$ are given functions. 
Here and in what follows, we denote 
\begin{align*}
\hat {u}^{\ep}_\delta(t) & = \frac{1}{\delta}u_\delta^\ep (t) = z_0^\ep(t) +\delta U^\ep_\delta (t), 
\\
M_{\delta}^{\ep}(t) u & = \delta \tilde{\eta}_\delta^\ep K u  + M(\hat u_\delta^\ep(t), u).       
\end{align*}
Recall that $\eta^{\ep}_{\delta} = \tilde \eta_0^\ep \delta^2 + \tilde \eta_1^\ep(\delta) \delta^3$. 

By Theorem \ref{thm4.2}, we have $||| \hat{u}^{\ep}_{\delta} |||_{\ep, \mathcal{Y}_{a} } \le C$ uniformly in $0 < \ep \le \ep_{1}$ and $|\delta| \le \delta_{1}$. Based on this, we see, by a perturbation argument, that there are positve constants $\ep_{2}$ and $\delta_{2}$ such that if $0 < \ep \le \ep_{2}$ and $|\delta| \le \delta_{2}$, then the problem \eqref{7.1}--\eqref{7.2} has a unique solution $u \in \mathcal{Y}(s, s + T_{a})$ and it satisfies the estimate 
\[
||| u |||_{\ep, \mathcal{Y}(s, s + T_{a}) } \le C\{ ||| v |||_{X^{1} } + ||| F |||_{\ep, \mathcal{X}(s, s + T_{a}) } \}
\]
uniformly for $s \in \mathbb{R}$, $0 < \ep \le \ep_{2}$ and $|\delta| \le \delta_{2}$. 

Let $\mathscr{U}_{\delta}^\ep(t, s)$ be the solution operator for \eqref{7.1}--\eqref{7.2} with $F = 0$ and let $\wt{\mathscr{U} }_{\delta}^\ep(t, s)$ be the solution operator for \eqref{7.1}--\eqref{7.2} with $v = 0$. The solution $u$ of \eqref{7.1}--\eqref{7.2} is then written as 
\begin{equation}\label{7.2'} 
u(t) = \mathscr{U}_{\delta}^\ep(t, s) v + \wt{\mathscr{U} }_{\delta}^\ep(t, s) F
\end{equation} 
and $\mathscr{U}_{\delta}^\ep(t, s)$ and $\wt{\mathscr{U} }_{\delta}^\ep(t, s)$ satisfy the uniform estimates 
\begin{align}
||| \mathscr{U}_{\delta}^\ep(t, s) v |||_{\ep, \mathcal{Y}(s, s + T_{a}) } & \le C||| v |||_{X^{1} },  \label{7.2''} \\ 
||| \wt{\mathscr{U} }_{\delta}^\ep(t, s) F  |||_{\ep, \mathcal{Y}(s, s + T_{a}) } & \le C ||| F |||_{\ep, \mathcal{X}(s, s + T_{a}) } \label{7.2'''} 
\end{align}
uniformly for $s \in \mathbb{R}$, $0 < \ep \le \ep_{2}$ and $|\delta| \le \delta_{2}$. 

Since the bifurcating solution $u_\delta^\ep (t)=\delta(z_0^\ep(t) +\delta U^\ep(t))$ is a time periodic function of period $T_a$, its stability is determined by the spectrum of $\mathscr{U}_{\delta}^{\ep}(T_{a},0)$. We shall investigate the spectrum of $\mathscr{U}_{\delta}^{\ep}(T_{a},0)$ as a perturbation of that of $\mathscr{V}^{\ep}_{\omega^{\ep}_{\delta}} (T_{a})$. Here recall that $\mathscr{V}^{\ep}_{\omega^{\ep}_{\delta} } (t)
 = e^{- \frac{at}{a^{\ep}(1 + \omega^{\ep}_{\delta}) }  
L^{\ep}_{\mathcal{R}^{\ep}_{1,c} } }|_{X^{1}}$. 
%The spectrum of $\mathscr{V}_{\omega}^{\ep}(T_a)$ for $|\omega| \le \frac{1}{4}$ has been investigated in the proof of Lemma \ref{5.8}. 
We shall prove the following result on the spectrum of $\mathscr{U}_{\delta}^{\ep}(T_{a}, 0)$ for $0 < |\delta| \ll 1$.

\begin{thm}\label{thm7.1}
There exist positive constants $\ep_{2}$, $\delta_{2}$ and $\Lambda_2$ such that if $0 < \ep \le \ep_{2}$ and $| \delta | \le \delta_{2}$, then 
\[
\sigma \left( \mathscr{U}_{\delta}^{\ep} \left( T_a, 0 \right) \right) \subset \left\{ 1, e^{\frac{a}{a^{\ep} } T_{a} \lambda^{\ep}_{\delta} } \right\} \cup \left\{ \mu \in \mathbb{C} ; |\mu|\le  e^{ - T_{a} \Lambda_2 } \right \} 
.\]
Here $1$ and $e^{T_{a}\lambda^{\ep}_{\delta} }$ are simple eigenvalues of $\mathscr{U}_{\delta}^{\ep} \left( T_{a}, 0 \right)$ and $\lambda^{\ep}_{\delta}$ satisfies
\[
\lambda_{\delta}^{\ep} =  2 \delta^{2} \tilde{\eta}_{0}^{\ep} \, \re[ K z^{\ep}_{0} ]_{ +, \ep } + O(\delta^3)
\]
uniformly for $0 < \ep \le \ep_{2}$ and $|\delta| \le \delta_{2}$. 
The eigenspace of the eigenvalue $1$ is spanned by $\delt u^{\ep}_{\delta}$.
\end{thm} 

To prove Theorem \ref{thm7.1}, we first show the following lemma.

\begin{lem}\label{lem7.2}
There exist positive constants $\Lambda_2$ and $r_0$ such that for all $0<r\le r_0$ there exists a positive constant $\delta_{2} = \delta_{2}(r)$ such that for $0<\ep\le\ep_{2}$ and $|\delta|\le\delta_{2}$, it holds that
\[
\rho \left( \mathscr{U}_{\delta}^{\ep} \left( T_a, 0 \right) \right) \supset \{ \mu ; |\mu - 1 | \ge r \} \cap \left\{ \mu ; |\mu| \ge e^{-T_{a} \Lambda_2 } \right\},
\]
and $\sigma\left(\mathscr{U}_{\delta}^{\ep}\left(T_{a},0\right)\right)\cap\{\mu;|\mu-1|<r\}$ consists of eigenvalues and the spectral projection associated with the set of these eigenvalues is a finite rank operator of rank $2$.
\end{lem}

To prove Lemma \ref{lem7.2}, we observe that $\mathscr{U}^{\ep}_{\delta} (t, \tau)$ is written as 
\[
\mathscr{U}^{\ep}_{\delta} (t, \tau) = \mathscr{V}_{\omega^{\ep}_{\delta}}^{\ep} (t - \tau) - \frac{a \delta}{a^{\ep} (1 + \omega^{\ep}_{\delta} ) } \int_\tau^t \mathscr{V}_{\omega^{\ep}_{\delta} }^{\ep} (t - s) M_{\delta}^{\ep} (s)\mathscr{U}_{\delta}^{\ep} (s, \tau) \,ds.
\]
See \eqref{7.3} below. Using this formula, we shall investigate the spectrum of $\mathscr{U}^{\ep}_{\delta} ( T_{a}, 0 )$ by a perturbation argument. We thus introduce an operator $\mathscr{S}_{\omega, \delta}^{\ep}$ on $\mathfrak{B}(X^1, \mathcal{Y} (0,T_a))$ defined by 
\[
(\mathscr{S}_{\omega, \delta}^{\ep} \mathscr{W} v)(t) = - \frac{a}{ a^{\ep} (1 + \omega) } \int_0^t \mathscr{V}_{\omega}^{\ep} (t - s)M_{\delta}^{\ep} (s) \mathscr{W} v (s)\,ds
\]
for $\mathscr{W}\in \mathfrak{B}(X^1, \mathcal{Y} (0,T_a))$ and $v \in X^1$. 
The following estimates hold for $\mathscr{S}_{\omega, \delta}^{\ep}$. 

\begin{lem}\label{lem7.3} 
{\rm (i)} Let $|\omega| \le \frac{1}{4}$. Then the operator $\mathscr{S}_{\omega, \delta}^{\ep}$ is a bounded operator on $\mathfrak{B}(X^1, \mathcal{Y} (0,T_a))$ and satisfies 
\[
\|\mathscr{S}_{\omega, \delta}^{\ep} \mathscr{W} \|_{\mathfrak{B}(X^1, \mathcal{Y} (0,T_a))}  \le C\| \mathscr{W} \|_{\mathfrak{B}(X^1, \mathcal{Y} (0,T_a))} 
\]
and 
\[
\|\mathscr{S}_{\omega, \delta}^{\ep} \|_{\mathfrak{B}(\mathfrak{B}(X^1, \mathcal{Y} (0,T_a)) )} \le C
\]
uniformly for $0 < \ep \le \ep_{2}$, $|\omega| \le \frac{1}{4}$ and $|\delta|\le \delta_{2}$. 

\vspace{1ex}   
{\rm (ii)} There exists a positive constant $\delta_{2}$ such that if $|\delta| \le \delta_{2}$, then $I - \delta \mathscr{S}_{\omega, \delta}^{\ep}$ has the bounded inverse $(I - \delta \mathscr{S}_{\omega, \delta}^{\ep} )^{-1} \in \mathfrak{B}(\mathfrak{B}(X^1, \mathcal{Y}  (0,T_a)) )$  
with estimate 
\[
\|(I - \delta \mathscr{S}_{\omega, \delta}^{\ep} )^{-1}\|_{\mathfrak{B}(\mathfrak{B}(X^1, \mathcal{Y} (0,T_a)) )}\le 2 
\]
uniformly for $0 < \ep \le \ep_{2}$, $|\omega| \le \frac{1}{4}$ and $|\delta|\le \delta_{2}$. 
\end{lem}

\vspace{1ex}
\noindent 
{\bf Proof.} Applying Lemmas \ref{lem5.1} and \ref{lem6.2}, we obtain  
\begin{align*} 
|||\mathscr{S}_{\omega, \delta}^{\ep} \mathscr{W}v |||_{\ep, \mathcal{Y} (0,T_a)} 
&  \le \frac{1}{2}|||\mathscr{S}_{\omega, \delta}^{\ep} \mathscr{W}v |||_{\ep, \mathcal{Y} (0,T_a)} \\
& \quad + C\left\{|||\mathscr{W} v |||_{\ep, \mathcal{X} (0, T_a)} +  ||| \hat u_\delta^\ep |||_{\ep, \mathcal{Y}_a} |||\mathscr{W} (\cdot) v |||_{\ep, \mathcal{Y} (0, T_a)} \right\} \\
& \ \le \frac{1}{2}|||\mathscr{S}_{\omega, \delta}^{\ep} \mathscr{W}v |||_{\ep, \mathcal{Y}  (0,T_a)}  +  C \|\mathscr{W} \|_{\mathfrak{B}(X^1, \mathcal{Y} (0,T_a))} ||| v |||_{\ep, X^1},
\end{align*}
and hence,
\[
 ||| \mathscr{S}_{\omega, \delta}^{\ep} \mathscr{W}v |||_{\ep, \mathcal{Y} (0,T_a)} \le C\|\mathscr{W} \|_{\mathfrak{B}(X^1, \mathcal{Y} (0,T_a))} ||| v |||_{\ep, X^1} 
\]
uniformly for $0< \ep \le \ep_{2}$, $|\omega| \le \frac{1}{4}$ and $|\delta| \le \delta_{2}$.  
This implies that 
\[
\|\mathscr{S}_{\omega, \delta}^{\ep} \mathscr{W} \|_{\mathfrak{B}(X^1, \mathcal{Y} (0,T_a))}  \le C\| \mathscr{W} \|_{\mathfrak{B}(X^1, \mathcal{Y} (0,T_a))},  
\]
and hence, 
\[
\|\mathscr{S}_{\omega, \delta}^{\ep} \|_{\mathfrak{B}(\mathfrak{B}(X^1, \mathcal{Y} (0,T_a)) )} \le C
\]
uniformly for $0 < \ep \le \ep_{2}$, $|\omega| \le \frac{1}{4}$ and $|\delta|\le \delta_{2}$.  
It then follows that if $|\delta|\le\frac{1}{2C}$, then there exists $(I - \delta \mathscr{S}_{\omega, \delta}^{\ep} )^{-1} \in \mathfrak{B}(\mathfrak{B}(X^1, \mathcal{Y} (0,T_a)) )$ with estimate $\|(I - \delta \mathscr{S}_{\omega, \delta}^{\ep} )^{-1}\|_{\mathfrak{B}(\mathfrak{B}(X^1, \mathcal{Y} (0,T_a)) )}\le 2$. This completes the proof. \hfill$\square$ 
 
\vspace{1ex}
We next give a proof of Lemma \ref{lem7.2}. 

\vspace{1ex}
\noindent
{\bf Proof of Lemma \ref{lem7.2}}. 
Let $u(t)$ be a solution of \eqref{7.1}--\eqref{7.2} with $F=0$. Since $u(t) = \mathscr{U}_{\delta}^{\ep}(t,\tau)v$ and $\delt u + \frac{a}{a^{\ep} } L_{\omega^{\ep}_{\delta} }^{\ep} u = - \frac{a \delta}{a^{\ep} (1 + \omega^{\ep}_{\delta}) } M_{\delta}^{\ep}(t) u$, the solution operator $\mathscr{U}_{\delta}^{\ep}(t,\tau)$ for \eqref{7.1} satisfies
\begin{equation}\label{7.3}
\begin{split} 
\mathscr{U}_{\delta}^{\ep} (t,\tau) v &= \mathscr{V}_{\omega^{\ep}_{\delta} }^{\ep}(t - \tau) v - \frac{a \delta}{a^{\ep} (1 + \omega^{\ep}_{\delta}) } \int_\tau^t \mathscr{V}_{ \omega^{\ep}_{\delta}}^{\ep} (t - s)M_{\delta}^{\ep} (s)u(s)\,ds \\
&= \mathscr{V}_{\omega^{\ep}_{\delta}}^{\ep} (t - \tau) v - \frac{a \delta}{a^{\ep} (1 + \omega^{\ep}_{\delta}) } \int_\tau^t \mathscr{V}_{\omega^{\ep}_{\delta} }^{\ep} (t - s)M_{\delta}^{\ep} (s)\mathscr{U}_{\delta}^{\ep} (s, \tau) v \,ds.
\end{split}
\end{equation}

We set $\mathscr{W}_{\delta}^{\ep} (t, s) = \mathscr{U}_{\delta}^{\ep} (t,s) - \mathscr{V}_{\omega^{\ep}_{\delta}}^{\ep} (t - s)$. It then follows that 
\begin{equation}\label{7.4}
\begin{split}
\mathscr{W}_{\delta}^{\ep} (t,\tau) v & = - \frac{a \delta}{ a^{\ep} (1 + \omega^{\ep}_{\delta}) } \int_\tau^t \mathscr{V}_{\omega^{\ep}_{\delta}}^{\ep} (t - s)M_{\delta}^{\ep} (s)\mathscr{U}_{\delta}^{\ep} (s, \tau) v\,ds \\
&  = - \frac{a \delta}{ a^{\ep} (1 + \omega^{\ep}_{\delta}) } \int_\tau^t \mathscr{V}_{\omega^{\ep}_{\delta} }^{\ep} (t - s)M_{\delta}^{\ep} (s)\mathscr{V}_{\omega^{\ep}_{\delta} }^{\ep} (s-\tau) v\,ds \\
& \quad \quad - \frac{a \delta}{ a^{\ep} (1 + \omega^{\ep}_{\delta}) } \int_\tau^t \mathscr{V}_{\omega^{\ep}_{\delta} }^{\ep} (t - s) M_{\delta}^{\ep} (s)\mathscr{W}_{\delta}^{\ep} (s, \tau) v\,ds.
\end{split}
\end{equation}

We claim that 
\begin{equation}\label{7.5}
\|\mathscr{W}_{\delta}^{\ep} (T_a,0)\|_{\mathfrak{B}(X^1)} \le 2 C |\delta|
\end{equation}
uniformly for $0 < \ep \le \ep_{2}$ and $|\delta| \le \delta_{2}$.  
Indeed, we see from \eqref{7.4} that $\mathscr{W}_{\delta}^{\ep} (t,0)$ is written as 
\[
\mathscr{W}_{\delta}^{\ep} (t,0) v = - \delta(\mathscr{S}^{\ep}_{\omega^{\ep}_{\delta},  \delta} \mathscr{V}^{\ep}_{\omega^{\ep}_{\delta} } (\cdot) v)(t) - \delta ( \mathscr{S}^{\ep}_{\omega^{\ep}_{\delta}, \delta } \mathscr{W}_{\delta}^{\ep} (\cdot,0) v)(t).
\] 
This implies
\[
(I + \delta \mathscr{S}^{\ep}_{\omega^{\ep}_{\delta}, \delta} ) \mathscr{W}_{\delta}^{\ep} (\cdot,0) = - \delta \mathscr{S}^{\ep}_{\omega^{\ep}_{\delta}, \delta} \mathscr{V}_{\omega^{\ep}_{\delta} }^{\ep} (\cdot) .
\]

By Lemmas \ref{lem5.1} and \ref{lem7.3}, we have  
\[
||| \mathscr{S}^{\ep}_{\omega^{\ep}_{\delta}, \delta} \mathscr{V}_{\omega^{\ep}_{\delta} }^{\ep} (\cdot) v |||_{\ep,\mathcal{Y} (0,T_a)} \le C ||| v |||_{\ep, X^1}
\]
uniformly for $0 <\ep \le \ep_{2}$ and $|\delta| \le\delta_{2}$. We also find from Lemma \ref{lem7.3} that 
\[
\mathscr{W}_{\delta}^{\ep} (\cdot,0) = - \delta(I + \delta \mathscr{S}^{\ep}_{\omega^{\ep}_{\delta}, \delta} )^{-1}\mathscr{S}^{\ep}_{\omega^{\ep}_{\delta}, \delta} \mathscr{V}_{\omega^{\ep}_{\delta} }^{\ep} (\cdot)
\]
and
\[
||| \mathscr{W}_{\delta}^{\ep} (\cdot,0) v |||_{\ep, \mathcal{Y} (0,T_a)} \le 2 |\delta| |||\mathscr{S}^{\ep}_{\omega^{\ep}_{\delta}, \delta} \mathscr{V}_{\omega^{\ep}_{\delta} }^{\ep} (\cdot) v |||_{\ep, \mathcal{Y} (0,T_a)} \le 2C |\delta| ||| v |||_{\ep, X^1}.
\]
In particular, we have $|||\mathscr{W}_{\delta}^{\ep} (T_a,0) v |||_{\ep, X^1}\le 2 C |\delta| ||| v |||_{\ep, X^1}$, which yields $\|\mathscr{W}_{\delta}^{\ep} (T_a,0)\|_{\mathfrak{B}(X^1)} \le 2 C |\delta|$ uniformly for $0 < \ep \le \ep_{2}$ and $|\delta| \le \delta_{2}$. This proves \eqref{7.5}. 
 
We next consider the resolvent of $\mathscr{U}_{\delta}^{\ep} (T_{a},0)$. For a fixed positive constant $r$ we set 
\[
\Sigma_r = \{\mu\in \mathbb{C}; |\mu-1|\ge r,\,|\mu|\ge\mathrm{e}^{-\frac{3}{4}\kappa_1 T_a}\}. 
\] 
We see from Lemma \ref{rem5.9} that if $|\omega^{\ep}_{\delta}| \le C\min\{r, 1\}$, then 
\[
\|(\mu-\mathscr{V}_{\omega^{\ep}_{\delta}}^{\ep} (T_a))^{-1}\|_{\mathfrak{B}(X^1)} \le C\left(\frac{1}{r} + \frac{1}{|\mu|}\right).
\] 
Since $\omega^{\ep}_{\delta} = \tilde{\omega}^{\ep}_{\delta} \delta^{2}$, there exists a positive constant $\delta_{2} = O(r)$ $(r \to 0)$  such that if $|\delta| \le \delta_{2}$, then  
\[
\|\mathscr{W}_{\delta}^{\ep} (T_a, 0)(\mu-\mathscr{V}_{\omega^{\ep}_{\delta}}^{\ep} (T_a))^{-1}\|_{\mathfrak{B}(X^1)} \le\frac{1}{2}. 
\]
We thus conclude that $\Sigma_r \subset \rho(\mathscr{U}_{\delta}^{\ep} (T_a, 0))$ and   
\begin{align*}
& (\mu - \mathscr{U}_{\delta}^{\ep} (T,0))^{-1} \\ 
& \quad = (\mu - \mathscr{V}_{\omega^{\ep}_{\delta}}^{\ep} (T_a))^{-1}(I - \mathscr{W}_{\delta}^{\ep} (T_a, 0)(\mu - \mathscr{V}_{\omega^{\ep}_{\delta}}^{\ep} (T_a))^{-1})^{-1} 
\end{align*}
for $\mu \in \Sigma_{r}$. Furthermore, it holds 
\[
\|(\mu - \mathscr{U}_{\delta}^{\ep} (T,0))^{-1}\|_{\mathfrak{B}(X^1)}  \le 2C\left(\frac{1}{r} + \frac{1}{|\mu|}\right)
\]
for $\mu \in \Sigma_r$ uniformly in $0 < \ep \le \ep_{2}$ and $|\delta| \le \delta_{2}$. 
This completes the proof. 
\hfill$\square$

We next investigate a part of $\sigma( \mathscr{U}_{\delta}^{\ep} \left(T_{a},0\right) )$ in $\{ \mu \in \mathbb{C};|\mu-1|<r \}$. To specify it, we consider the spectrum of the linearized operator $B_{\delta}^{\ep}$ around the bifurcating time periodic solution $u_{\delta}^{\ep}$ of \eqref{4.5} which is given by
\begin{eqnarray}\label{7.6}
B_{\delta}^{\ep} = B^{\ep}( \omega^{\ep}_{\delta} ) + \delta M_{\delta}^{\ep},  
\end{eqnarray}
where $M_{\delta}^{\ep}$ is the operator defined by 
\[
(M_{\delta}^{\ep} u)(t) = M_{\delta}^{\ep}(t) u(t) \quad \mbox{\rm for $t \in \mathbb{T}_{T_{a} }$ } 
\]
with $M_{\delta}^{\ep}(t)$ given in \eqref{7.1}. Note that $B^{\ep}_{0} = B^{\ep}$. 

\begin{thm}\label{thm7.4}
There exists a positive constant $r_0$ such that 
\[
\sigma(- B_{\delta}^{\ep} )\cap\{\lambda\in\mathbb{C};\,|\lambda|\le r_0\} = \{0,\lambda_{\delta}^{\ep} \},
\]
where $0$ and $\lambda_{\delta}^{\ep}$ are simple eigenvalues of $-B_{\delta}^{\ep} $. Furthermore, $\lambda_{\delta}^{\ep}$ satisfies
\[
\lambda_{\delta}^{\ep} =  2 \delta^{2} \tilde{\eta}_{0}^{\ep} \, \re[K z_{0}^{\ep}]_{+,\ep} + O(\delta^3)
\]
uniformly for $0 < \ep \le \ep_{2}$ and $|\delta| \le \delta_{2}$. 
\end{thm}

To give a proof of Theorem \ref{thm7.4}, we introduce notation. We change the bases of $N(B^{\ep}_{0})$ from $\{z_+^{\ep}, z_-^{\ep}\}$ to $\{z^{\ep}_{0}, z^{\ep}_{1} \}$, where 
\[
z^{\ep}_{0} = 2\re z_+^{\ep}, \quad z^{\ep}_{1} = 2\im z_+^{\ep},
\]
and likewise, the dual bases from $\{z_+^{\ep*}, z_-^{\ep*}\}$ to $\{z^{\ep *}_{0}, z^{\ep *}_{1} \}$, where
\[
z^{\ep*}_{0} = \re{ z^{\ep*}_{+} }, \quad z^{\ep*}_{1} = \im{ z^{\ep*}_{+} }.
\]
We set
\[
\llbracket u\rrbracket_{j, \ep} = \langle u, z^{\ep*}_{j} \rangle_\ep \quad (j = 0, 1) 
\]
and define the operators $\tilde{ \mathcal{P} }^{\ep}$ and $\tilde{ \mathcal{Q} }^{\ep}$ by 
\[
\tilde{\mathcal{P}}^{\ep} u = \llbracket u\rrbracket_{0, \ep} z^{\ep}_{0} + \llbracket u\rrbracket_{1, \ep} z^{\ep}_{1}
\]
and 
\[
\tilde{ \mathcal{Q} }^{\ep} = I - \tilde{ \mathcal{P} }^{\ep}, 
\] 
respectively. It then follows that $\tilde{ \mathcal{P} }^{\ep}$ is an eigenprojection for the eigenvalue $0$ of $-B^{\ep}_{0}$. 
Observe also that 
\[
\llbracket \delt u \rrbracket_{0, \ep} =  a \llbracket u \rrbracket_{1, \ep}, \quad 
\llbracket \delt u \rrbracket_{1, \ep} = -a \llbracket u \rrbracket_{0, \ep}. 
\]
Furthermore, it holds that 
\[
\mbox{ \rm $\llbracket u\rrbracket_{0, \ep} = \re [u]_{+, \ep}$ if $u$ is a real valued function. }
\]

We see from Lemma \ref{lem5.8} that if $\re \lambda \ge -\frac{3}{4}\kappa_{1} T_{a}$, then the problem 
\begin{equation}\label{7.7}
(\lambda + B^{\ep} ( \omega ) ) u = F
\end{equation}
has a unique solution $u \in \tilde{\mathcal{Q} }^{\ep} \mathcal{Y}_{a}$ for any given $F \in \tilde{\mathcal{Q} }^{\ep} \mathcal{X}_{a}$ with estimate 
\[
||| u |||_{\ep, \mathcal{Y}_{a} } \le C ||| F |||_{\ep, \mathcal{X}_{a} }
\]
uniformly for $0 < \ep \le \ep_{2}$, $|\omega| \le \frac{1}{4}$ and $\lambda$ with $\re \lambda \ge -\frac{3}{4}\kappa_{1} T_{a}$ and $|\im\lambda|\le \frac{a}{2}$. 
We denote the solution operator for problem \eqref{7.7} by $\mathcal{S}^{\ep} (\lambda, \omega)$, namely, the operator $\mathcal{S}^{\ep} (\lambda, \omega): \tilde{ \mathcal{Q} }^{\ep} \mathcal{X}_{a} \to \tilde{ \mathcal{Q} }^{\ep} \mathcal{Y}_{a}$ is defined by   
\[
u = \mathcal{S}^{\ep}(\lambda, \omega) F,
\]
where $u \in \tilde{ \mathcal{Q} }^{\ep} \mathcal{Y}_{a}$ is the unique solution of \eqref{7.7}.  

We now consider the resolvent problem 
\begin{equation}\label{7.8}
(\lambda + B^{\ep}_{\delta} ) u = F.
\end{equation} 
We shall employ the Lyapunov-Schmidt method to investigate problem \eqref{7.8}. 

We decompose $u$ in \eqref{7.8} into its $\tilde{ \mathcal{P} }^{\ep}$ and $\tilde{ \mathcal{Q} }^{\ep}$ parts as 
\[
u = \zeta_{0} z^{\ep}_{0} + \zeta_{1} z^{\ep}_{1} + U,
\]
where $\zeta_{j} = \llbracket u \rrbracket_{j, \ep}$ $(j = 0, 1)$ and $U = \tilde{ \mathcal{Q} }^{\ep} u \in \tilde{ \mathcal{Q} }^{\ep} \mathcal{Y}_{a}$.  
Applying $\tilde{ \mathcal{P} }^{\ep}$ and $\tilde{ \mathcal{Q} }^{\ep}$ to \eqref{7.8}, we see that \eqref{7.8} is reduced to 
\begin{align} 
\lambda \zeta_{0} + a^{\ep} \omega^\ep_\delta \zeta_{1} + \eta^{\ep}_\delta \llbracket K u \rrbracket_{0, \ep} + \delta \llbracket M( \hat{u}^{\ep}_{\delta}, u ) \rrbracket_{0, \ep} & = \llbracket F \rrbracket_{0, \ep}, \label{7.9} \\[1ex]
\lambda \zeta_{1} - a^{\ep} \omega^\ep_\delta \zeta_{0} + \eta^{\ep}_\delta \llbracket K u \rrbracket_{1, \ep} + \delta \llbracket M( \hat{u}^{\ep}_{\delta}, u ) \rrbracket_{1, \ep} & = \llbracket F \rrbracket_{1, \ep}, \label{7.10} \\[1ex]
(\lambda + \tilde{\mathcal{Q} }^{\ep} B^{\ep}_{\delta} ) U +\delta \tilde{ \mathcal{Q} }^{\ep} \sum_{j = 0, 1} \zeta_j M^{\ep}_{\delta} z^{\ep}_{j} & = \tilde{ \mathcal{Q} }^{\ep} F, \label{7.11}
\end{align}  
where $u = \zeta_{0} z^{\ep}_{0} + \zeta_{1} z^{\ep}_{1} + U$ with $\b{\zeta} = \trans(\zeta_{0}, \zeta_{1} ) \in \mathbb{C}^{2}$ and $U \in \tilde{Q}^{\ep} \mathcal{Y}_{a}$. . 

We shall reduce \eqref{7.9}--\eqref{7.11} to a two-dimensional problem by solving \eqref{7.11} for $U$ in terms of $\zeta_j\,(j = 0,1)$ and $F$ and then substituting $U$ into \eqref{7.9} and \eqref{7.10}. To this end, we next consider the problem 
\begin{equation}\label{7.12}
(\lambda + \tilde{\mathcal{Q}}^\ep B^{\ep}_{\delta} )u = \tilde{\mathcal{Q}}^\ep F.
\end{equation}

\begin{prop}\label{prop7.5}
There exist positive constants $\ep_{2}$ and $\delta_{2}$ such that if $0<\ep\le\ep_{2},\,|\delta|\le\delta_{2}$, $\re\lambda\ge-\frac{3}{4}\kappa_1 T_a$ and $|\im \lambda| \le \frac{a}{2}$, then for any given $F\in\mathcal{X}_a$, problem \eqref{7.12} has a unique solution $U\in\tilde{\mathcal{Q}}^\ep\mathcal{Y}_a$ and $U$ satisfies the estimate
\[
\opnorm{U}_{\ep,\mathcal{Y}_a} \le C\opnorm{F}_{\ep,\mathcal{X}_a}
\]
uniformly for $0<\ep\le\ep_{2},\,|\delta|\le\delta_{2}$ and $\lambda$ with $\re\lambda\ge-\frac{3}{4}\kappa_1 T_a$ and $|\im \lambda| \le \frac{a}{2}$.
\end{prop}

\noindent
{\bf Proof.} We regard $\tilde{\mathcal{Q}}^\ep B_\delta^\ep\tilde{\mathcal{Q}}^\ep = B^\ep(\omega^{\ep}_{\delta}) + \delta \tilde{\mathcal{Q}}^\ep M_\delta^\ep$ as a perturbation of $B^\ep(\omega^{\ep}_{\delta})$. By Lemma \ref{lem5.8} and Lemma \ref{6.2}, we see that if $\re\lambda\ge-\frac{3}{4}\kappa_1 T_a$ and $|\im \lambda| \le \frac{a}{2}$, then
\begin{align*}
\opnorm{\mathcal{S}^\ep(\lambda,\omega^{\ep}_{\delta})\tilde{\mathcal{Q}}^\ep\tilde{M}_\delta^\ep u}_{\ep,\mathcal{Y}_a} &\le \frac{1}{2}\opnorm{\mathcal{S}^\ep(\lambda,\omega^{\ep}_{\delta})\tilde{\mathcal{Q}}^\ep\tilde{M}_\delta^\ep u}_{\ep,\mathcal{Y}_a} \\
& \quad \quad+ C\{\opnorm{u}_{\ep,\mathcal{Y}_a} + \opnorm{\hat{u}_\delta^\ep}_{\ep,\mathcal{Y}_a}\opnorm{u}_{\ep,\mathcal{Y}_a}\}
\end{align*}
uniformly for $u \in \tilde{\mathcal{Q}}^\ep\mathcal{Y}_a$, $0<\ep\le\ep_{2}$, $|\delta|\le\delta_{2}$ and $\re\lambda\ge-\frac{3}{4}\kappa_1 T_a$. This implies that there exists a positive constant $\delta_{2}$ such that if $|\delta|\le\delta_{2}$, then $I + \delta\mathcal{S}^\ep(\lambda,\omega^{\ep}_{\delta})\tilde{\mathcal{Q}}^\ep M_\delta^\ep\tilde{\mathcal{Q}}^\ep$ has a bounded inverse $(I + \delta\mathcal{S}^\ep(\lambda,\omega^{\ep}_{\delta})\tilde{\mathcal{Q}}^\ep M_\delta^\ep\tilde{\mathcal{Q}}^\ep)^{-1}$ on $\tilde{\mathcal{Q}}^\ep\mathcal{Y}_a$ with estimate
\begin{equation*}
\opnorm{(I + \delta \mathcal{S}^\ep(\lambda,\omega^{\ep}_{\delta})\tilde{\mathcal{Q}}^\ep M_\delta^\ep\tilde{\mathcal{Q}}^\ep)^{-1}u}_{\ep,\mathcal{Y}_a} \le C\opnorm{u}_{\ep,\mathcal{Y}_a}
\end{equation*}
uniformly for $0<\ep\le\ep_{2},\,|\delta|\le\delta_{2}$ and $\re\lambda\ge-\frac{3}{4}\kappa_1 T_a$.
It then follows that \eqref{7.12} has a unique solution $U = (I + \delta\mathcal{S}^\ep(\lambda,\omega^{\ep}_{\delta})\tilde{\mathcal{Q}}^\ep M_\delta^\ep\tilde{\mathcal{Q}}^\ep)^{-1}\mathcal{S}^\ep(\lambda,\omega^{\ep}_{\delta})F \in \tilde{\mathcal{Q}}^\ep\mathcal{Y}_a$ and $U$ satisfies the estimate
\[
\opnorm{U}_{\ep,\mathcal{Y}_a} \le C\opnorm{\mathcal{S}^\ep(\lambda,\omega^{\ep}_{\delta})\tilde{\mathcal{Q}}^\ep F}_{\ep,\mathcal{Y}_a}\le C\opnorm{F}_{\ep,\mathcal{X}_a}
\]
uniformly for $0<\ep\le\ep_{2},\,|\delta|\le\delta_{2}$, $\re\lambda\ge-\frac{3}{4}\kappa_1 T_a$ and $|\im \lambda| \le \frac{a}{2}$. This completes the proof. 
\hfill$\square$

\vspace{1ex}
Let $0<\ep\le\ep_{2}$ and $|\delta|\le\delta_{2}$. For $\lambda$ satisfying $\re\lambda \ge-\frac{3}{4}\kappa_1 T_a$, we denote the solution operator for \eqref{7.12} by $\tilde{\mathcal{S}}_\delta^\ep(\lambda)$.

In terms of $\tilde{\mathcal{S}}_\delta^\ep(\lambda)$, \eqref{7.11} is written as
\[
U = \tilde{\mathcal{S}}_\delta^\ep(\lambda)\tilde{\mathcal{Q}}^\ep F - \delta\sum_{j=0,1}\zeta_j \tilde{\mathcal{S}}_\delta^\ep(\lambda)\tilde{\mathcal{Q}}^\ep M_\delta^\ep z_j^\ep.
\]
Substituting this into \eqref{7.9} and \eqref{7.10} and using the fact that $\llbracket M(z_0^\ep,z_j^\ep)\rrbracket_{k,\ep}=0,\,(j,k \in \{0,1\})$, we have
\begin{equation}\label{7.13}
\lambda\b{\zeta} + \delta^2\Gamma_\delta^\ep(\lambda)\b{\zeta} = \mathcal{F}^\ep_\delta(\lambda)F.
\end{equation}
Here $\b{\zeta} = \trans(\zeta_0,\zeta_1)\in \mathbb{C}^2$; $\Gamma_\delta^\ep(\lambda)$ is a $2\times2$ matrix given by 
\begin{align*}
\Gamma_\delta^\ep(\lambda)
& =
\begin{pmatrix}
\tilde{\eta}^{\ep}_{\delta} \llbracket Kz_0^\ep \rrbracket_{0,\ep} & a^\ep \tilde{\omega}^\ep_\delta + \tilde{\eta}^\ep_\delta \llbracket Kz_1^\ep \rrbracket_{0,\ep} \\\vspace{2ex}
\quad \quad \quad\quad + \llbracket M(U_\delta^\ep,z_0^\ep) \rrbracket_{0,\ep} & \quad \quad \quad \quad\quad + \llbracket M(U_\delta^\ep,z_1^\ep) \rrbracket_{0,\ep}\\
-a^\ep \tilde{\omega}^\ep_\delta + \tilde{\eta}^\ep_\delta \llbracket Kz_0^\ep \rrbracket_{1,\ep}  & \tilde{\eta}^\ep_\delta \llbracket Kz_1^\ep \rrbracket_{1,\ep} \\
\quad\quad\quad\quad\quad\quad+ \llbracket M(U_\delta^\ep, z_0^\ep) \rrbracket_{1,\ep} & \quad\quad\quad\quad+ \llbracket M(U_\delta^\ep, z_1^\ep) \rrbracket_{1,\ep}
\end{pmatrix}
\\[2ex]
& \quad -
\begin{pmatrix}
\llbracket M_\delta^\ep\tilde{S}_\delta^\ep(\lambda)\tilde{\mathcal{Q}}^\ep M_\delta^\ep z_0^\ep \rrbracket_{0,\ep} & \llbracket M_\delta^\ep\tilde{S}_\delta^\ep(\lambda)\tilde{\mathcal{Q}}^\ep M_\delta^\ep z_1^\ep \rrbracket_{0,\ep} \\[1ex]
\llbracket M_\delta^\ep\tilde{S}_\delta^\ep(\lambda)\tilde{\mathcal{Q}}^\ep M_\delta^\ep z_0^\ep \rrbracket_{1,\ep} & \llbracket M_\delta^\ep\tilde{S}_\delta^\ep(\lambda) \tilde{\mathcal{Q}}^\ep M_\delta^\ep z_1^\ep \rrbracket_{1,\ep} 
\end{pmatrix}
\end{align*}
and
\begin{equation*}
\mathcal{F}_\delta^\ep(\lambda) F
=
\begin{pmatrix}
\llbracket F \rrbracket_{0,\ep} \\[1ex]
\llbracket F \rrbracket_{1,\ep}
\end{pmatrix}
-
\delta
\begin{pmatrix}
\llbracket M_\delta^\ep\tilde{S}_\delta^\ep(\lambda)\tilde{\mathcal{Q}}^\ep F \rrbracket_{0,\ep} \\[1ex]
\llbracket M_\delta^\ep\tilde{S}_\delta^\ep(\lambda)\tilde{\mathcal{Q}}^\ep F \rrbracket_{1,\ep}
\end{pmatrix}.
\end{equation*}

\begin{prop}\label{prop7.6}
There exists a positive constant $c_0$ such that if $0<\ep\le\ep_{2}$ and $|\delta|\le\delta_{2}$ then $\{\lambda\in\mathbb{C};\,|\lambda|>c_0\delta^2, \re\lambda\ge-\frac{3}{4}\kappa_1 T_a, |\im \lambda| \le \frac{a}{2}\} \subset \rho(- B_\delta^\ep)$ and 
\[
\opnorm{(\lambda + B_\delta^\ep)^{-1}F}_{\ep,\mathcal{Y}_a} \le C\left(\frac{1}{|\lambda|-c_0|\delta|^2} + 1\right)\opnorm{F}_{\ep,\mathcal{X}_a}
\]
uniformly for $0<\ep\le\ep_{2},\,|\delta|\le\delta_{2}$ and $\lambda$ with $|\lambda|>c_0|\delta|^2$, $\re\lambda\ge-\frac{3}{4}\kappa_1 T_a$ and $|\im \lambda| \le \frac{a}{2}$.
\end{prop}

\noindent
{\bf Proof.} By using Proposition \ref{prop7.5}, we see that
\[
|\Gamma_\delta^\ep(\lambda)\b{\zeta}| \le c_0|\b{\zeta}|
\]
uniformly for $0<\ep\le\ep_{2},\,|\delta|\le\delta_{2}$ and $\lambda$ with $\re\lambda\ge-\frac{3}{4}\kappa_1 T_a$. It then follows that if $|\lambda| > c_0|\delta|^2$, $\re\lambda \ge -\frac{3}{4}\kappa_1 T_a$ and $|\im \lambda| \le \frac{a}{2}$, then \eqref{7.13} has a unique solution $\b{\zeta}\in\mathbb{C}^2$ with estimate
\[
|\b{\zeta}| \le \frac{C}{|\lambda|-c_0|\delta|^2}\opnorm{F}_{\ep,\mathcal{X}_a},
\]
from which we deduce that $\{\lambda;\,|\lambda| > c_0|\delta|^2, \re\lambda\ge-\frac{3}{4}\kappa_1 T_a, |\im \lambda| \le \frac{a}{2} \} \subset \rho(-B_\delta^\ep)$ and 
\[
(\lambda + B_\delta^\ep)^{-1}F = (z_0^\ep\,\,z_1^\ep)(\lambda + \delta^2\Gamma_\delta^\ep(\lambda))^{-1} \mathcal{F}_\delta^\ep(\lambda) F + \tilde{\mathcal{S}}_\delta^\ep(\lambda)\tilde{\mathcal{Q}}^\ep F,
\]
\[
\opnorm{(\lambda + \tilde{B}_\delta^\ep)^{-1}F}_{\ep,\mathcal{Y}_a} \le C\left(\frac{1}{|\lambda|-c_0|\delta|^2} + 1\right)\opnorm{F}_{\ep,\mathcal{X}_a}.
\]
This completes the proof. \hfill$\square$

\vspace{1ex}
We take a positive constant $\delta_{2}$ so that $c_0\delta_{2}^2 < \frac{3}{4}\kappa_1 T_a$. It then follows from Proposition \ref{prop7.6} that
\[
\sigma(- B_\delta^\ep)\cap\{\lambda\in\mathbb{C}; \re\lambda\ge-\frac{3}{4}\kappa_1 T_a, |\im \lambda| \le \frac{a}{2} \} \subset \{\lambda\in\mathbb{C};\,|\lambda| < c_0\delta^2\}
\]
for $0<\ep\le\ep_{2}$ and $|\delta|\le\delta_{2}$. We shall, therefore, investigate
\[
\sigma(-B_\delta^\ep)\cap\{\lambda\in\mathbb{C};\,|\lambda|\le c_0\delta^2\}
\]
for $0<\ep\le\ep_{2}$ and $|\delta|\le\delta_{2}$. From the argument above, we see that
\begin{equation*}
\lambda \in \sigma(-B_\delta^\ep)\cap\{\lambda\in\mathbb{C};\,|\lambda|\le c_0\delta^2\} \mbox{ if and only if } {\rm det} \left(\lambda I + \delta^2\Gamma_\delta^\ep(\lambda) \right) = 0.
\end{equation*}
We thus consider zeros of ${\rm det} \left(\lambda I + \delta^2\Gamma_\delta^\ep(\lambda) \right)$. To this end, we regard $\Gamma_\delta^\ep(\lambda)$ as a perturbation of $\Gamma_0^\ep(0)$.

\begin{prop}\label{prop7.7}
If $|\lambda| \le c_1|\delta|$, then
\[
|\Gamma_\delta^\ep(\lambda)-\Gamma_0^\ep(0)| \le C|\delta|
\]
uniformly for $0<\ep\le\ep_{2},\,|\delta|\le\delta_{2}$ and $|\lambda| \le c_1|\delta|$.
\end{prop}

\noindent
{\bf Proof.} We first prove 
\begin{equation}\label{7.14}
\opnorm{\tilde{\mathcal{S}}_\delta^\ep(\lambda)\tilde{\mathcal{Q}}^\ep M_\delta^\ep z_j^\ep - \mathcal{S}^\ep(0,0)\tilde{\mathcal{Q}}^\ep M_0^\ep z_j^\ep}_{\ep,\mathcal{Y}_a} \le C|\delta|\,\,(j=0,1)
\end{equation}
uniformly for $0<\ep\le\ep_{2},\,|\delta|\le \delta_{2}$ and $|\lambda| \le c_1|\delta|$. To show \eqref{7.14}, we write
\[
\tilde{\mathcal{S}}_\delta^\ep(\lambda)\tilde{\mathcal{Q}}^\ep M_\delta^\ep z_j^\ep - \mathcal{S}^\ep(0,0)\tilde{\mathcal{Q}}^\ep M_0^\ep z_j^\ep = I_1 + I_2 + I_3 + I_4,
\]
where
\begin{align*}
& I_1 = \tilde{\mathcal{S}}_\delta^\ep(\lambda)\tilde{\mathcal{Q}}^\ep( M_\delta^\ep z_j^\ep - M_0^\ep z_j^\ep), \quad \quad   \  
I_2 = (\tilde{\mathcal{S}}_\delta^\ep(\lambda) - \tilde{\mathcal{S}}_\delta^\ep(0))\tilde{\mathcal{Q}}^\ep M_0^\ep z_j^\ep, \\
& I_3 = (\tilde{\mathcal{S}}_\delta^\ep(0) - \mathcal{S}^\ep(0,\omega^{\ep}_{\delta}))\tilde{\mathcal{Q}}^\ep M_0^\ep z_j^\ep, \quad   
I_4 = (\mathcal{S}^\ep(0,\omega^{\ep}_{\delta}) - \mathcal{S}^\ep(0,0))\tilde{\mathcal{Q}}^\ep M_0^\ep z_j^\ep.
\end{align*}

As for $I_1$, since $M_\delta^\ep z_j^\ep - M_0^\ep z_j^\ep = \delta \tilde{\eta}_\delta^\ep K z_j^\ep + \delta M(U^{\ep}_{\delta},z_j^\ep)$, we have $\opnorm{I_1}_{\ep,\mathcal{Y}_a} \le C|\delta|$. 
To estimate $I_2$, we note that
\[
\tilde{\mathcal{S}}_\delta^\ep(\lambda) = \tilde{\mathcal{S}}_\delta^\ep(0)\sum_{N=0}^\infty(-1)^N\lambda^N\tilde{\mathcal{S}}_\delta^\ep(0)^N.
\]
This, together with Proposition \ref{prop7.5}, implies that
\[
\opnorm{I_2}_{\ep,\mathcal{Y}_a} \le C|\lambda|\opnorm{\tilde{\mathcal{S}}_\delta^\ep(0) M_0^\ep z_j^\ep}_{\ep,\mathcal{Y}_a} \le C|\delta|.
\] 
As for $I_3$, we see from the proof of Proposition \ref{prop7.5} that 
\[
\opnorm{I_3}_{\ep,\mathcal{Y}_a} \le C|\delta|\opnorm{\mathcal{S}^\ep(0,\omega^{\ep}_{\delta})\tilde{\mathcal{Q}}^\ep M_0^\ep z_j^\ep}_{\ep,\mathcal{Y}_a} \le C|\delta|.
\]
Since $(\mathcal{S}^\ep(0,\omega^{\ep}_{\delta}) - \mathcal{S}^\ep(0,0))\tilde{\mathcal{Q}}^\ep = (B^{\ep}(\omega^{\ep}_{\delta})\tilde{\mathcal{Q}}^\ep)^{-1} - (B^{\ep}\tilde{\mathcal{Q}}^\ep)^{-1}$ and $M_0^\ep z_j^\ep = M(z^{\ep}_{0}, z^{\ep}_{j})$, we estimate $I_4$ as in the proof of Theorem \ref{thm4.2} (i) to obtain $\opnorm{I_4}_{\ep,\mathcal{Y}_a} \le C|\delta|$. The estimate \eqref{7.14} is thus proved.

It follows from \eqref{7.14} that
\[
|\llbracket M_\delta^\ep\tilde{\mathcal{S}}_\delta^\ep(\lambda)\tilde{\mathcal{Q}}^\ep M_\delta^\ep z_j^\ep \rrbracket_{k,\ep} - \llbracket M_\delta^\ep\mathcal{S}^\ep(0,0)\tilde{\mathcal{Q}}^\ep M_0^\ep z_j \rrbracket_{k,\ep} | \le C|\delta|, \quad (j, k \in\{0,1\})
\]
uniformly for $0<\ep\le\ep_{2},\,|\delta|\le\delta_{2}$ and $\lambda\in\Sigma$. Since $\tilde{\eta}^\ep_\delta  = \tilde{\eta}_0^\ep + O(\delta)$ and $\tilde{\omega}^\ep_\delta = \tilde{\omega}_0^\ep + O(\delta)$, we conclude that $|\Gamma_\delta^\ep(\lambda) - \Gamma_0^\ep(0)| \le C|\delta|$ uniformly for $0<\ep\le\ep_{2}$, $|\delta| \le\delta_{2}$ and $|\lambda| \le c_1|\delta|$. This completes the proof. \hfill$\square$

\vspace{1ex}
We are now in a position to prove Theorem \ref{thm7.4}. We set 
\[
D_\delta^\ep(\tilde{\lambda}) = {\rm det} \left(\tilde{\lambda} + \Gamma_\delta^\ep(\delta^2\tilde{\lambda}) \right).
\]
If $D_\delta^\ep(\tilde{\lambda})\neq0$, then $\lambda = \delta^2\tilde{\lambda} \in \rho(-B_\delta^\ep)$ and 
\[
(\lambda + B_\delta^\ep)^{-1}F = \begin{pmatrix}z_0^\ep & z_1^\ep\end{pmatrix}(\lambda + \delta^2 \Gamma_\delta^\ep(\lambda))^{-1}\mathcal{F}_\delta^\ep(\lambda)F + \tilde{\mathcal{S}}_\delta^\ep(\lambda)\tilde{\mathcal{Q}}^\ep F.
\]
If $D_\delta^\ep(\tilde{\lambda}) = 0$, then $\lambda = \delta^2\tilde{\lambda} \in \sigma(-B_\delta^\ep)$. Therefore, we investigate zeros of $D_\delta^\ep(\tilde{\lambda})$ in $\{\tilde{\lambda} \in\mathbb{C};\,|\tilde{\lambda}| \le \tilde{c}_0\}$, where $\tilde{c}_0 = \frac{c_0}{\delta_{2}}$. We shall prove the following proposition. 

\begin{prop}\label{prop7.8}
There exists a positive constant $\delta_{2}$ such that if $0<\ep\le\ep_{2}$ and $|\delta| \le\delta_{2}$ then $\tilde{D}_\delta^\ep(\tilde{\lambda})$ has two zeros $\tilde{\lambda} = 0$ and $\tilde{\lambda}_\delta^\ep$ in $\{\tilde{\lambda}\in\mathbb{C};\,|\tilde{\lambda}| \le\tilde{c}_0\}$ and both zeros are of order $1$. Here $\tilde{\lambda}_\delta^\ep$ satisfies $\tilde{\lambda}_\delta^\ep = \tilde{\lambda}_0^\ep + O(\delta),\,\tilde{\lambda}_0^\ep = 2\tilde{\eta}_0^\ep\real[Kz_0^\ep]_{+,\ep}$. 
\end{prop}

It follows from Proposition \ref{prop7.8} that
\[
\sigma(-B_\delta^\ep)\cap\{\lambda\in\mathbb{C};\,|\lambda| \le c_0\delta^2\} = \{0,\lambda_\delta^\ep\}, 
\]
where $0$ and $\lambda_\delta^\ep = \delta^2\tilde{\lambda}_\delta^\ep$ are simple eigenvalues of $-B_\delta^\ep$. This proves Theorem \ref{thm7.4}. Furthermore, we shall see from the proof of Proposition \ref{prop7.8} below that the eigenspaces for the eigenvalues $0$ and $\lambda_\delta^\ep$ are spanned by $\delt u_\delta^\ep$ and $z_0^\ep + O(\delta)$, respectively.

\vspace{1ex}
Let us prove Proposition \ref{prop7.8}.

\vspace{1ex}
\noindent
{\bf Proof of Proposition \ref{prop7.8}.} We first observe that
\begin{equation}\label{7.15}
\Gamma_\delta^\ep(0)
=
\begin{pmatrix}
\gamma_{11,\delta}^\ep(0) & 0\\
\gamma_{21,\delta}^\ep(0) & 0
\end{pmatrix},
\end{equation}
where
\begin{align*}
\gamma_{11,\delta}^\ep(\lambda) & = \tilde{\eta}^\ep_\delta \llbracket Kz_0^\ep \rrbracket_{0,\ep} + \llbracket M(U_\delta^\ep,z_0^\ep) 
\rrbracket_{0,\ep} -\llbracket M_\delta^\ep\tilde{S}_\delta^\ep(\lambda)\tilde{\mathcal{Q}}^\ep M_\delta^\ep z_0^\ep \rrbracket_{0,\ep},
\\[1ex]
\gamma_{21,\delta}^\ep(\lambda) & = -a^\ep \tilde{\omega}^\ep_\delta + \tilde{\eta}^\ep _\delta \llbracket Kz_0^\ep \rrbracket_{1,\ep} + \llbracket M(U_\delta^\ep,z_0^\ep) \rrbracket_{1,\ep} - \llbracket M_\delta^\ep\tilde{\mathcal{S}}_\delta^\ep(\lambda)\tilde{\mathcal{Q}}^\ep M_\delta^\ep z^\ep_0 \rrbracket_{1,\ep}.
\end{align*}
Indeed, differentiating the equation 
\begin{equation}\label{7.16}
B^\ep(\omega^{\ep}_{\delta})u_\delta^\ep + \eta^\ep_\delta K u_\delta^\ep + N(u_\delta^\ep) = 0
\end{equation}
in $t$, we have
\[
B_\delta^\ep \delt U_\delta^\ep + \delta\frac{a^\ep}{a}\tilde{\omega}^\ep_\delta \delt^2 z_0^\ep + M_\delta^\ep\delt z_0^\ep = 0.
\]
This, together with $\delt z_0^\ep = -az_1^\ep$ and Proposition \ref{prop7.6}, gives \eqref{7.15}. It then follows that $\tilde{\lambda} = 0$ is a zero of $D_\delta^\ep(\tilde{\lambda})$.

We next claim that 
\begin{equation}\label{7.17}
D_0^\ep(\tilde{\lambda}) = \tilde{\lambda}(\tilde{\lambda} + \tilde{\lambda}_0^\ep).
\end{equation}
This can be verified by using \eqref{7.16}. Indeed, we see from \eqref{7.16} that
\begin{equation}\label{7.18}
B^\ep(\omega^{\ep}_{\delta})U_\delta^\ep + \frac{a^\ep}{a}\delta\tilde{\omega}^\ep_\delta \delt z_0^\ep + \delta\tilde{\eta}^\ep_\delta K(z_0^\ep + \delta U_\delta^\ep) + N(z_0^\ep + \delta U_\delta^\ep) = 0.
\end{equation}
Noting that $\llbracket M(z_0^\ep,z_0^\ep) \rrbracket_{j,\ep} = 0$ $(j=0,1)$, we have
\[
\tilde{\eta}^\ep_\delta \llbracket Kz_0^\ep \rrbracket_{0,\ep} + \delta\tilde{\eta}^\ep_\delta \llbracket KU_\delta^\ep\rrbracket_{0,\ep} + \llbracket M(U_\delta^\ep,z_0^\ep)\rrbracket_{0,\ep} + \frac{\delta}{2}\llbracket M(U_\delta^\ep,U_\delta^\ep)\rrbracket_{0,\ep} = 0
\]
and 
\[
U_\delta^\ep = -\mathcal{S}^\ep(0,\omega^{\ep}_{\delta})\tilde{\mathcal{Q}}^\ep\{\delta \tilde{\eta}^\ep_\delta K(z_0^\ep + \delta U_\delta^\ep) + N(z_0^\ep + \delta U_\delta^\ep)\}.
\]
Letting $\delta\rightarrow0$, we obtain
\begin{equation}\label{7.19}
\tilde{\eta}_0^\ep\llbracket Kz_0^\ep \rrbracket_{0,\ep} + \llbracket M(U_0^\ep,z_0^\ep) \rrbracket_{0,\ep} = 0
\end{equation}
and 
\begin{equation}\label{7.20}
U_0^\ep = -\frac{1}{2}\tilde{\mathcal{S}}_0^\ep(0)\tilde{\mathcal{Q}}^\ep M(z_0^\ep,z_0^\ep).
\end{equation}
We deduce from \eqref{7.19} and \eqref{7.20} that
\[
\gamma_{11,0}^\ep(0) = -\llbracket M(z_0^\ep,\tilde{\mathcal{S}}_0^\ep(0)M(z_0^\ep,z_0^\ep))\rrbracket_{0,\ep} = -2\tilde{\eta}_0^\ep\llbracket Kz_0^\ep\rrbracket_{0,\ep}=-\tilde{\lambda}_0^\ep,
\]
and hence,
\[
\Gamma_0^\ep(0)
=
\begin{pmatrix}
-\tilde{\lambda}_0^\ep & 0\\
\gamma_0^\ep(0) & 0
\end{pmatrix}.
\]
This implies \eqref{7.17}.

We take $\delta_{2}$ so small that $\tilde{c}_0 = \frac{c_0}{\delta_{2}} \le 2|\tilde{\lambda}_0^\ep|$. By Proposition \ref{prop7.6}, we obtain
\[
|D_\delta^\ep(\tilde{\lambda}) - D_0^\ep(\tilde{\lambda})| \le c_1|\delta|
\]
for $0<\ep\le\ep_{2},\,|\delta|\le\delta_{2}$ and $|\tilde{\lambda}| \le \tilde{c}_0$. Therefore, there exists a positive constant $\delta_{2}$ such that if $|\delta| \le \delta_{2}$, then $|D_0^\ep(\tilde{\lambda})|>|D_\delta^\ep(\tilde{\lambda})-D_0^\ep(\tilde{\lambda})|$ on $|\tilde{\lambda}| = \tilde{c}_0$. Applying the Rouch\'{e} theorem, we find that $D_\delta^\ep(\tilde{\lambda})$ has two zeros in $\{\tilde{\lambda}\in\mathbb{C};\,|\tilde{\lambda}|\le \tilde{c}_0\}$; one of them is $\tilde{\lambda} = 0$. Furthermore, we see that, by taking $\delta_{2}$ smaller if necessary, that there exists a positive constant $c_3$ such that if $|\delta| \le\delta_{2}$, then $|D_0^\ep(\tilde{\lambda})| > |D_\delta^\ep(\tilde{\lambda})-D_0^\ep(\tilde{\lambda})|$ on $|\tilde{\lambda}-\tilde{\lambda}_0^\ep| = c_3|\delta|$. The Rouch\'{e} theorem then implies that $D_\delta^\ep(\tilde{\lambda})$ has a zero $\tilde{\lambda}_\delta^\ep$ of order $1$ in $\{\tilde{\lambda}\in\mathbb{C};\,|\tilde{\lambda}-\tilde{\lambda}^\ep_{0}|<c_3|\delta|\}$. This completes the proof. \hfill$\square$

\vspace{1ex}
\noindent
{\bf Proof of Theorem \ref{thm7.1}} Let $\mu_0$ and $\mu_1$ be eigenvalues of $\mathscr{U}_{\delta}^{\ep} \left(T_{a},0\right)$ with $|\mu_j - 1|<r$ $(j = 0,1)$. We know that one of $\mu_j$'s, say $\mu_0$, is equal to $1$. On the other hand, as in the proof of Theorem \ref{thm5.6}, i.e., the proof of \cite[Theorem 4.2]{Hsia-Kagei-Nishida-Teramoto1}, we see from \eqref{7.2'}--\eqref{7.2'''} that if $\mathrm{e}^{ \frac{a}{a^{\ep} } T_{a}\lambda}\in\rho\left(\mathscr{U}_{\delta}^{\ep} \left(T_{a},0\right)\right)$, then $\lambda\in\rho(-B_{\delta}^{\ep})$. Therefore, since $\sigma\left(\mathscr{U}_{\delta}^{\ep} \left(T_{a},0\right)\right)\cap\{\mu;|\mu-1|<r\} = \{\mu_0,\mu_1\}$, by Theorem \ref{thm7.4}, we conclude that $\mu_0 = 1,\,\mu_1 = \mathrm{e}^{ \frac{a}{a^{\ep} }  T_{a}\lambda^{\ep}_{\delta}}$. This completes the proof.\hfill$\square$

\vspace{2ex}
\noindent
{\bf Acknowledgements.}
Y. Kagei was partly supported by JSPS KAKENHI
Grant Numbers 16H03947, 16H06339 and 20H00118. T. Nishida is partly supported by JSPS KAKENHI Grant Number 20K03699.

\end{document}